\documentclass[]{interact}
\usepackage{graphicx}
\usepackage{color}
\usepackage{endfloat}
\usepackage{amsmath}
\usepackage{amssymb}

\newtheorem{thm}{\sc Theorem}[section]
\newtheorem{prop}{\sc Proposition}[section]

\renewcommand\thefigure{\arabic{figure}}
\renewcommand\thetable{\arabic{table}}

\usepackage[captionskip=5ex,farskip=5ex]{subfig}

\usepackage[numbers,sort&compress]{natbib}
\bibpunct[, ]{[}{]}{,}{n}{,}{,}
\renewcommand\bibfont{\fontsize{10}{12}\selectfont}
\makeatletter
\def\NAT@def@citea{\def\@citea{\NAT@separator}}
\makeatother

\begin{document}
%

\baselineskip 0.49cm   

\title{Real-time reconstruction of moving point/dipole wave sources from boundary measurements}
\author{
\name{Takashi Ohe\thanks{Email:ohe@xmath.ous.ac.jp}}
\affil{Department of Applied Mathematics, Faculty of Science, 
Okayama University of Science, 1-1 Ridai-cho, Kita-ku, Okayama 700-0005,
JAPAN}
}

\maketitle

\begin{abstract}
This paper is concerned with a reconstruction method for multiple moving
 point/dipole wave sources.
We assume that the number, locations, and magnitudes/moments of wave
 sources are unknown, and consider the problem to reconstruct these
 parameters from the measurement
 of the wave field on the boundary.
For this problem, we derive algebraic relations between the parameters of wave
 sources and the reciprocity gap functionals, and propose a
real-time reconstruction procedure based on these relations.
We perform some numerical experiments, and show the effectiveness of our
 reconstruction procedure.
\end{abstract}


\begin{amscode}
35L05, 35R30,65M32
\end{amscode}


%
\setcounter{page}{1}
\section{Introduction}

Reconstruction problem for wave sources often arises in
science, engineering and medical fields, and has many important
applications in these areas, e.\@g.\@ reconstruction of seismic source,
identification of noisy sound source in the human body, and so on
\cite{Ando1,Crocker1,Wang1}.
Such kind of problem can be formulated as an \emph{inverse source problem} for
wave equation.
In this paper, we consider a simplified case such that the media is
homogeneous and isotropic, and measurement of the wave field is
given on the whole boundary of the domain.
Then we can formulate the problem as follows:
\par
Let $u$ be a wave field described as a solution of
the following initial-boundary value problem for three dimensional scalar wave equation:
\begin{equation} 
\left\{
\begin{array}{rll}
  \displaystyle 
  \frac{1}{c^2}\partial_t^2 u(t,\boldsymbol{r}) - \Delta u(t,\boldsymbol{r}) &=  F(t,\boldsymbol{r})  & \mbox{ in } (0,T) \times \Omega, \\
  u(0,\boldsymbol{r}) &= 0,  & \mbox{ in } \Omega, \\
  \partial_t u(0,\boldsymbol{r}) &= 0,  & \mbox{ in } \Omega, \\
  u(t,\boldsymbol{r}) &= 0, &  \mbox{ on } (0,T) \times \Gamma,
\end{array}
\right.
\label{eq:mixed_problem_wave_eq}
\end{equation}
where $\boldsymbol{r} = (x,y,z)$, $\Omega \subset {\mathbb R}^3$ is a simply connected bounded
domain with $C^\infty$-boundary $\Gamma = \partial\Omega$, $c>0$
is the given wave propagation speed, $T > 2\cdot {\mathrm{diam}}\, \Omega/c >0$ is a given constant, 
and $F(t,\boldsymbol{r})$ describes the wave source.
In (\ref{eq:mixed_problem_wave_eq}), suppose that the wave source
$F$ is unknown, and
consider the problem to reconstruct unknown $F$ from observation
$\phi$ on $\Gamma$ given by
\begin{equation*}
   \phi(t,\boldsymbol{r}) = \partial_\nu u(t,\boldsymbol{r}), \mbox{ on } [0,T]\times\Gamma,
\end{equation*}
where $\partial_\nu u$ denotes the outward normal derivative of $u$ on $\Gamma$.
\par
Many researchers discussed the reconstruction of wave sources 
in theoretical and numerical points of view, e.\@g.\@ \cite{Anikonov1,Bruckner1,Isakov1,Isakov2,Rossing1}.
Unfortunately, the solution of this inverse problem is not unique in
general, and hence we usually restrict the source term to some ideal models.
In this paper, we assume that the source term is described by multiple moving
point sources
\begin{equation}
  F(t,\boldsymbol{r}) = \sum_{k=1}^K q_k(t) \delta(\boldsymbol{r} - \boldsymbol{p}_k(t)),
  \label{eq:point_source}
\end{equation}
or moving dipole sources
\begin{equation}
  F(t,\boldsymbol{r}) = - \sum_{k=1}^K \boldsymbol{m}_k(t)\cdot\nabla \delta(\boldsymbol{r} - \boldsymbol{p}_k(t)).
  \label{eq:dipole_source}
\end{equation}
In (\ref{eq:point_source}), $K$ denotes the number of point sources,
and $\boldsymbol{p}_k(t) \equiv (p_{k,x}(t),\,p_{k,y}(t),\,p_{k,z}(t)) \in D$ and $q_k(t) \in {\mathbb R}$ denote the location and
magnitude of $k$-th point source at $t$, respectively, where $D$ is a
compact subset of $\Omega$.
Note that we need not specify $D$ explicitly.
The symbol $\delta$ describes the Dirac's delta distribution which is
understood as a linear functional
\begin{equation*}
   \langle \delta(\cdot-\boldsymbol{p}), v(\cdot)\rangle_\Omega =
    v(\boldsymbol{p}), \qquad \boldsymbol{p} \in \Omega,
\end{equation*}
for any $v \in C(\overline{\Omega})$.
Also in (\ref{eq:dipole_source}), $K$ denotes the number of dipole sources,
and $\boldsymbol{p}_k(t) \in D$ and $\boldsymbol{m}_k(t) \equiv (m_{k,x}(t),\,m_{k,y}(t),\,m_{k,z}(t))\in {\mathbb R}^3$ denote the location
and dipole moment of $k$-th dipole source at $t$, respectively.
The gradient of $\delta$ is understood as
\begin{align*}
   \langle \nabla \delta(\cdot-\boldsymbol{p}), v(\cdot)\rangle_\Omega 
 =&  - \langle \delta(\cdot-\boldsymbol{p}),
 \nabla v(\cdot)\rangle_\Omega  \nonumber \\
 =& - \nabla v(\boldsymbol{p}), \qquad \boldsymbol{p} \in \Omega
\end{align*}
for any $v \in C^1(\Omega) \cap C(\overline{\Omega})$.
Hence, we consider the solution $u$ of (\ref{eq:mixed_problem_wave_eq}) in a weak sense, 
i.\@e.\@ $u \in C^1((0,T];H^1_0(\Omega))$ that satisfies 
\begin{align}
  &  \displaystyle  \frac{1}{c^2}\int_\Omega \partial_t u(T,\boldsymbol{r}) v(T,\boldsymbol{r})  dV(\boldsymbol{r})
    - \frac{1}{c^2}\int_\Omega u(T,\boldsymbol{r})         \partial_t
    v(T, \boldsymbol{r}) dV(\boldsymbol{r}) \nonumber \\
  & \displaystyle - \int_0^T \int_\Gamma \partial_\nu u(t, \boldsymbol{r}) v(t, \boldsymbol{r})
      dS(\boldsymbol{r}) dt \nonumber \\
  & +\displaystyle \int_0^T \int_\Omega u(t, \boldsymbol{r}) 
       \left( \frac{1}{c^2}\partial_t^2	 v(t, \boldsymbol{r}) - \Delta v(t, \boldsymbol{r}) \right)
    dV(\boldsymbol{r})dt \nonumber \\
 = & \ \  {\mathcal F}(v)
  \label{eq:weak_form_mixed_problem}
\end{align}
for any $v \in C^\infty([0,T]\times\overline{\Omega})$.
The right hand side term ${\mathcal F}(v)$ is expressed by
\begin{align}
   {\mathcal F}(v) 
=& \displaystyle 
    \int_0^T \sum_{k=1}^K q_k(t) \langle \delta( \cdot -
\boldsymbol{p}_k(t)), v(t,\cdot)\rangle_\Omega dt \nonumber \\
=& \displaystyle
    \sum_{k=1}^K \int_0^T q_k(t) v(t,\boldsymbol{p}_k(t)) dt,
\label{eq:F_for_point_source}
\end{align}
for moving point sources (\ref{eq:point_source}), and
\begin{align}
   {\mathcal F}(v) 
=& \displaystyle
  -\int_0^T \sum_{k=1}^K \boldsymbol{m}_k(t) \cdot \langle \nabla \delta( \cdot -
\boldsymbol{p}_k(t)), v(t,\cdot)\rangle_\Omega dt \nonumber \\
=& \displaystyle \sum_{k=1}^K \int_0^T \boldsymbol{m}_k(t) \cdot \nabla
 v(t,\boldsymbol{p}_k(t)) dt.
\label{eq:F_for_dipole_source}
\end{align}
for moving dipole sources (\ref{eq:dipole_source}).
\par

%
Reconstruction methods of unknown sources can be categorised into two
 types, one is based on an optimisation technique like the least
squares method\cite{Bruckner1,Komornik1,Komornik2,Ohnaka1,Rashedi1}, 
and the other is based on an algebraic idea \cite{Badia2,Badia3,Inui1,Ohe1}.
The former method has high
versatility, but it usually requires large computational cost because we
need to solve partial differential equation iteratively.
In the latter method, the problem is usually reduced into algebraic equations of parameters of
unknown sources. 
This kind of method has a merit on the computational cost since we need not
 to solve partial differential equation iteratively instead of a
 drawback on relatively low versatility.
In this paper, we focus our interest on the latter type of methods.
\par
The first work on algebraic reconstruction methods for wave sources was done by El Badia and Ha Duong in 2001\cite{Badia2,Badia3}.
They assumed that several point sources are fixed in $D$,
i.\@e.\@ $\boldsymbol{p}_k(t) \equiv \boldsymbol{p}_k \in D$, 
and derived a reconstruction procedure from observations $\phi$
 on the whole of $\Gamma$ and the final state
$u(T,\boldsymbol{r})$ and $\partial_t u(T,\boldsymbol{r})$ in $\Omega$.
The key ideas of their method are {\em the reciprocity gap functional}
\cite{Andrieux1} with respect to the space variable $\boldsymbol{r}$
and the Fourier transform with respect to the time variable $t$.
Since their method needs to compute the inverse Fourier transform in the
reconstruction procedure, 
we can not reconstruct unknown source is simultaneously on observing data.
In 2011, Inui, Ohnaka and the author extended  the idea of
\cite{Badia3}, and
propose a new algebraic procedure in which one can reconstruct the
locations and magnitudes of several point sources 
almost simultaneously on observations \cite{Ohe1}.
 We can also find some researches which apply the reciprocity gap
functional to localise a point wave source\cite{Dogan1,Dogan2}.
For the reconstruction of moving wave source, only a few papers are
published in present(e.\@g.\@ \cite{Hu1,Nakaguchi1,Wang2, Wang3}),
and in these papers
unknown source have been assumed as a single point source.
The purpose of this paper is to extend our result \cite{Ohe1} to the
case of multiple moving point and dipole sources, and give a
 real-time reconstruction procedure, where the phrase 'real-time' means simultaneously on observing data.
\par
The contents of this paper are as follows.
Section 2 shows some regularity results for the solution of (\ref{eq:mixed_problem_wave_eq}),
especially, for observation data $\phi$. 
In section 3, we describe the detail of our reconstruction method;
First, we define the reciprocity gap
functional for the wave equation.
Next, we give an appropriate set of functions, 
and derive algebraic relations between the reciprocity gap functionals
for these functions and parameters of unknown sources.
Based on these relations,
we propose a real-time reconstruction procedure for unknown sources.
In section 4, we perform some numerical experiments and discuss the
effectiveness of our reconstruction procedure. 
\par
Throughout this paper, we use the following notations for some functional
spaces.
The Lebesgue space of the square integrable functions on $\Omega$ is
denoted by $L^2(\Omega)$, and the set of all functions $u$ of which the weak
derivative $D^\alpha u \in L^2(\Omega)$ for every multiindex
$|\alpha| \leq m$ by $H^m(\Omega)$.
$H^m_0(\Omega)$ denotes the closure of ${\mathcal D}(\Omega)$ in
$H^m(\Omega)$ where ${\mathcal D}(\Omega)$ is the set of
all infinitely differentiable functions with the compact support in $\Omega$.
We denote the dual space of $H^m_0(\Omega)$ by $H^{-m}(\Omega)$.
\par

\section{Some regularity results for observation data}

Before discussion of the reconstruction of wave sources, we show a regularity
result on the solution of initial-boundary value problem (\ref{eq:mixed_problem_wave_eq}) for moving point and
dipole sources,  especially, the regularity of the observation data
$\phi =\partial_\nu u$ on $\Gamma$.
\par
\begin{prop}
Let $\ell$ be a non-negative integer.
For moving point sources (\ref{eq:point_source}), assume that $\boldsymbol{p}_k \in C^{\ell+2}([0,T]; D),\ 
q_k \in C^{\ell+1}([0,T];\mathbb{R})$,  $|d_t{\boldsymbol{p}}_k(t)| < c$ and $q_k(0) =
 d_t{q}_k(0) = d^2_t{q}_k(0)=\cdots = d_t^{(\ell +1)}q_k(0) = 0$, where $d_t$ denotes the
 derivative with respect to $t$.
Then, the solution $u$ of (\ref{eq:mixed_problem_wave_eq}) satisfies $u
 \in C([0,T]; L^2(\Omega))$ and $\partial_t u \in C([0,T]; L^2(\Omega))$.
Specifically, the restriction $u$ on $[0,T]\times (\Omega \backslash D)$
 satisfies
\begin{align}
  u|_{[0,T]\times(\Omega\backslash D)} \in 
& C([0,T];H^{\ell+1}(\Omega\backslash D)), \label{eq:regularity_1}\\
 \nabla u|_{[0,T]\times(\Omega\backslash D)} \in
& C([0,T];H^{\ell}(\Omega\backslash D)), \\
 \nabla \partial_t^m u|_{[0,T]\times(\Omega\backslash D)} \in
& C([0,T];H^{\ell-m}(\Omega\backslash D)), \quad 0 \leq m \leq \ell,
\end{align}
and the normal derivative $\partial_\nu u$ on $\Gamma$ satisfies
\begin{equation}
 \partial_\nu u \in H^{\ell}([0,T]\times\Gamma).
  \label{eq:regularity_4} 
\end{equation}
\par
 For moving dipole sources (\ref{eq:dipole_source}), assume that
 $\boldsymbol{p}_k \in C^{\ell+2}([0,T]; D),\ \boldsymbol{m}_k \in
 C^{\ell+2}([0,T];{\mathbb R}^3)$, $|d_t{\boldsymbol{p}}_k(t)| < c$, and $\boldsymbol{m}_k(0) =
 d_t{\boldsymbol{m}}_k(0) = d_t^2{\boldsymbol{m}}_k(0) = \cdots = d_t^{(\ell+2)}\boldsymbol{m}_k(0) = \boldsymbol{0}$.
Then, the solution $u$ satisfies $u \in C([0,T]; H^{-1}(\Omega))$, and $\partial_t u \in C([0,T]; H^{-1}(\Omega))$.
Specifically, the restriction $u$ on $[0,T]\times(\Omega \backslash D)$ and the
 normal derivative $\partial_\nu u$ have the same regularity properties
 (\ref{eq:regularity_1})-(\ref{eq:regularity_4}) as for moving point sources.
\end{prop}

\noindent
\begin{proof}
First, we consider the case where the source term $F$ is expressed by moving point
sources (\ref{eq:point_source}).
Let $u_{\mathrm N}$ be the solution of the following initial value
problem for the wave equation in the whole space ${\mathbb R}^3$:
\begin{equation*}
 \left\{
 \begin{array}{rll}
 \displaystyle \frac{1}{c^2}\partial_t^2 u_{\mathrm N} - \Delta u_{\mathrm N}&= \displaystyle
  \sum_{k=1}^K q_k(t)\delta(\boldsymbol{r}-\boldsymbol{p}_k(t)),& \qquad \mbox{in
  } (0,T) \times {\mathbb R}^3, \\
 u_{\mathrm N}(0,\boldsymbol{r}) &= 0, & \qquad \mbox{ in } {\mathbb R}^3, \\
 \partial_t u_{\mathrm N}(0,\boldsymbol{r}) &= 0, & \qquad \mbox{ in } {\mathbb
  R}^3.
  \end{array}
\right.
\end{equation*}
Then, $u_{\mathrm N}$ is explicitly expressed by the following extended Li\'{e}nard-Wiechert retarded potential\cite{Jackson1}:
\begin{equation}
  u_{\mathrm N}(t,\boldsymbol{r}) = \sum_{k=1}^K
   \frac{1}{4\pi}\cdot\frac{q_k(s_k(t, \boldsymbol{r}))}{\left|\boldsymbol{r} -
						\boldsymbol{p}_k(s_k(t,\boldsymbol{r}))\right|
   \cdot h_k(s_k(t,\boldsymbol{r}), \boldsymbol{r})},
   \label{eq:special_sol_point_source}
\end{equation}
where $s_k(t, \boldsymbol{r})$ is determined as a solution $s$ of the equation
\begin{equation*}
   t = s +  \frac{\left|\boldsymbol{r} - \boldsymbol{p}_k(s)\right|}{c},
\end{equation*}
for each $t, \boldsymbol{r}$ and $k$,  and 
\begin{equation*}
  h_k(s,\boldsymbol{r}) = 1 - \frac{d_t{\boldsymbol{p}}_k(s)\cdot(\boldsymbol{r} -
 \boldsymbol{p}_k(s))}{c\left|\boldsymbol{r}-\boldsymbol{p}_k(s) \right|}.
\end{equation*}
Since $\boldsymbol{p}_k\in C^{\ell+2}([0,T];D)$ and $|d_t{\boldsymbol{p}}_k(t)|< c$, $s_k(t, \boldsymbol{r})$ is uniquely determined and
$s_k \in C^{\ell+1}(([0,T] \times \Omega) \backslash
\bigcup_{t \in (0,T)}
\{(t,\boldsymbol{p}_k(t))\}) \cap C([0,T]\times\Omega)$, 
and hence we have $h_k \in C^{\ell +1}(([0,T] \times \Omega) \backslash \bigcup_{t \in (0,T)}
\{(t,\boldsymbol{p}_k(t))\})$.
For the derivation of (\ref{eq:special_sol_point_source}), see e.\@g.\@ Appendix A in \cite{Nakaguchi1}.
The restriction of $u_{\mathrm N}$ on $[0,T]\times\Omega$
satisfies $u_{\mathrm{N}}|_{[0,T]\times\Omega} \in
C([0,T];L^2(\Omega))$ since $s_k(t,\boldsymbol{p}_k(t)) = t$ and $|\boldsymbol{r} -
\boldsymbol{p}_k(s_k(t,\boldsymbol{r}))| =
 O(|\boldsymbol{r} - \boldsymbol{p}_k(t)|)$ as $\boldsymbol{r} \rightarrow \boldsymbol{p}_k(t)$.
Moreover, the restriction to $[0,T]\times({\mathbb{R}^3}\backslash D)$ satisfies
$u_{\mathrm N}|_{[0,T]\times({\mathbb{R}^3}\backslash D)}\in C^{\ell+1}([0,T]\times({\mathbb{R}^3}\backslash D))$
and $\nabla u_{\mathrm N}|_{[0,T]\times({\mathbb{R}^3}\backslash D)} \in
C^{\ell}([0,T]\times({\mathbb{R}^3}\backslash D))$.
Hence, $u_{\mathrm N}|_{[0,T]\times\Gamma} \in C^{\ell+1}([0,T]\times\Gamma) \subset
H^{\ell+1}([0,T]\times\Gamma)$ and $\bm \nu \cdot \nabla u_{\mathrm N}|_{[0,T]\times\Gamma} \in C^{\ell}([0,T]\times\Gamma) \subset
H^{\ell}([0,T]\times\Gamma)$ since $\Gamma$ is in $C^{\infty}$-class,
where $\bm \nu$ is the
unit vector normal to $\Gamma$.
\par
Next, let $u_{\mathrm H}$ be the solution of initial-boundary value problem of the
homogeneous wave equation with inhomogeneous Dirichlet condition on $\Gamma$:
\begin{equation}
 \left\{
 \begin{array}{rll}
 \displaystyle
 \frac{1}{c^2} \partial_t^2 u_{\mathrm H} - \Delta u_{\mathrm H}&= 0, & \qquad \mbox{ in
  } (0,T) \times \Omega, \\
 u_{\mathrm H}(0,\boldsymbol{r}) &= 0, & \qquad \mbox{ in } \Omega, \\
 \partial_t u_{\mathrm H}(0,\boldsymbol{r}) &= 0, & \qquad \mbox{ in } \Omega, \\
 u_{\mathrm H} &= -u_{\mathrm N}|_{[0,T]\times\Gamma}, &  \qquad \mbox{ on }
  [0,T]\times\Gamma.
\end{array}
\right.
  \label{eq:homogeneous_wave_eq_0p}
\end{equation}
Owing to Remark 2.10 in \cite{Lasiecka1}, 
the solution $u_{\mathrm H}$ of (\ref{eq:homogeneous_wave_eq_0p}) exists in
$C([0,T];H^{\ell+1}(\Omega))$ since $u_{\mathrm N}|_{[0,T]\times\Gamma}$ satisfies all compatibility condition
$\partial_t^{k} u_{\mathrm N}(0,\boldsymbol{r}) = 0$ for
$k=0,1,2,\cdots,\ell+1$ on $\Gamma$.
Moreover, $\partial^{m}_t u_{\mathrm H} \in C([0,T];H^{\ell-m+1}(\Omega))$ and $\partial_\nu u_{\mathrm H} \in
H^{\ell}([0,T]\times\Gamma)$.
\par
Let $u = u_{\mathrm H} + u_{\mathrm N}$, then $u$ becomes the
solution of initial-boundary value problem (\ref{eq:mixed_problem_wave_eq}) in
$C([0,T];L^2(\Omega))$,  and satisfies (\ref{eq:regularity_1})-(\ref{eq:regularity_4}).
\par
Next, we consider the case where the source term is expressed by dipole sources model (\ref{eq:dipole_source}).
As same as for point sources model, let us consider the solution
$u_{\mathrm N}$  of the following initial value
problem in $\mathbb{R}^3$: 
\begin{equation*}
 \left\{
 \begin{array}{rll}
 \displaystyle
 \frac{1}{c^2} \partial_t^2 u_{\mathrm N} - \Delta u_{\mathrm N}&= -\displaystyle
  \sum_{k=1}^K \boldsymbol{m}_k(t)\cdot\nabla\delta(\boldsymbol{r}-\boldsymbol{p}_k(t)), &
  \mbox{ in } (0,T) \times {\mathbb R}^3, \\
 u_{\mathrm N}(0,\boldsymbol{r}) &= 0, & \mbox{ in } {\mathbb R}^3, \\
 \partial_t u_{\mathrm N}(0,\boldsymbol{r}) &= 0, & \mbox{ in } {\mathbb
  R}^3.
  \end{array}
 \right.
\end{equation*}
Then, $u_{\mathrm N}$ is explicitly expressed as
\begin{align}
{u_{\mathrm N}(t, \boldsymbol{r})} =   
&  \sum_{k=1}^K
   \frac{1}{4\pi}\cdot
   \frac
   {\boldsymbol{m}_k(s_k(t,\boldsymbol{r}))\cdot(\boldsymbol{r}-\boldsymbol{p}_k(s_k(t, \boldsymbol{r}))}
   {\left|\boldsymbol{r} -\boldsymbol{p}_k(s_k(t,\boldsymbol{r}))\right|^3
   \cdot h_k(s_k(t,\boldsymbol{r}), \boldsymbol{r})} 
 \nonumber \\
&+  \sum_{k=1}^K
   \frac{1}{4\pi c}\cdot
  \partial_t \left(
   \frac
   { \boldsymbol{m}_k(s_k(t,\boldsymbol{r}))\cdot(\boldsymbol{r}-\boldsymbol{p}_k(s_k(t, \boldsymbol{r})))}
   {\left|\boldsymbol{r} -\boldsymbol{p}_k(s_k(t,\boldsymbol{r}))\right|^2
   \cdot (h_k(s_k(t,\boldsymbol{r}), \boldsymbol{r}))}  \right)
\label{eq:special_sol_dipole_source}
\end{align}
We can derive (\ref{eq:special_sol_dipole_source}) using a similar
discussion as the one for moving point sources, but we omit it.
%
Under the assumption for $\boldsymbol{p}_k$ and $\boldsymbol{m}_k$, one can see
that $u_{\mathrm{N}} \in C^{\ell+1}([0,T];H^{-1}(\Omega))$.
Also for restrictions $u_{\mathrm N}|_{[0,T]\times({\mathbb{R}^3 \backslash D})}$, 
 $u_{\mathrm N}|_{[0,T]\times\Gamma}$ and $\bm \nu \cdot \nabla u_{\mathrm
d}|_{[0,T]\times\Gamma}$, we can show the same regularities as for moving point
sources.
Then, using the same discussion, we can derive 
that the solution $u$ of (\ref{eq:mixed_problem_wave_eq}) is in $C([0,T];H^{-1}(\Omega))$, and (\ref{eq:regularity_1})-(\ref{eq:regularity_4}).
\par
\end{proof}

\section{Reconstruction method}
\subsection{Reciprocity gap functional}
In this section, we present a reconstruction method for moving point and
dipole sources from observation $\phi=\partial_\nu u$ on $\Gamma$.
The key idea of our method is {\emph{the reciprocity gap functional}} that
is defined on the subspace of test functions $v$ in (\ref{eq:weak_form_mixed_problem}).
This idea is widely applied to
various inverse problems, e.\@g.\@ inverse conductivity problems and
inverse scattering problems \cite{Athanasiadis1,
Andrieux1,Chung1,Colton1,Badia1,Haddar1,Inui1,Nara1}.
First, we show a definition of the reciprocity gap functional for scalar wave equations.
\par
Let ${\mathcal W} \subset C^\infty([0,T]\times\overline{\Omega};{\mathbb
C})$ be a set of complex-valued
functions $v$ that satisfy the homogeneous wave equation with vanishing
condition at $t=T$: 
\begin{equation*}
  \left\{
  \begin{array}{rll}
     \displaystyle
     \frac{1}{c^2}\partial_t^2 v - \Delta v &= 0, &\qquad \mbox{ in } (0,T)\times\Omega, \\
     v(T,\boldsymbol{r}) &= 0, & \qquad \mbox{ in } \Omega, \\
     \partial_t v(T,\boldsymbol{r}) &= 0, & \qquad \mbox{ in } \Omega.
  \end{array}
  \right.
\end{equation*}
For given observation data $\phi \in L^2((0,T)\times\Gamma)$, we define
the {\emph{reciprocity gap functional}} $\mathcal R_\phi$ on ${\mathcal W}$  as follows: 
\begin{equation}
  {\mathcal R}_\phi(v) \equiv - \int_0^T \int_\Gamma \phi (t, \boldsymbol{r}) v(t,
   \boldsymbol{r})  dS(\boldsymbol{r}) dt, \qquad v \in {\mathcal W}.
   \label{eq:definition_reciprocity_gap}
\end{equation}
Then, since $\phi = \partial_\nu u$ and $u$ satisfies the weak form
(\ref{eq:weak_form_mixed_problem}), 
we obtain
\begin{align}
  {\mathcal R}_\phi(v)  
= & \,\displaystyle
    {\mathcal F}(v) -\frac{1}{c^2} \int_\Omega \partial_t u(T,\boldsymbol{r})
 v(T,\boldsymbol{r}) dV(\boldsymbol{r}) 
    +\frac{1}{c^2} \int_\Omega  u(T,\boldsymbol{r})
 \partial_t v(T,\boldsymbol{r}) dV(\boldsymbol{r}) \nonumber \\
 &  \displaystyle
    -\int_0^T \int_\Omega u(t, \boldsymbol{r}) 
       \left( \frac{1}{c^2}\partial_t^2	 v(t, \boldsymbol{r}) - \Delta v(t, \boldsymbol{r}) \right)
    dV(\boldsymbol{r})dt \nonumber \\
= & \,{\mathcal F}(v).
   \label{eq:relation_reciprocity_gap}
\end{align}
\par
The equation (\ref{eq:relation_reciprocity_gap})  gives the
relation between the reciprocity gap functional ${\mathcal R}_\phi(v)$
and the source term $F$, and suggests that we can reconstruct the
unknown source term from ${\mathcal R}_\phi(v)$ for appropriate choice
of functions $v$.
In the following subsections, we give our choice of functions $v \in \cal
W$,  and propose reconstruction procedures of moving point sources
and dipole sources.
\par
\subsection{Reconstruction of moving point sources}
Throughout this subsection, we assume that $\boldsymbol{p}_k \in C^{6}([0,T]; D),\ 
q_k \in C^5([0,T];\mathbb{R})$ and $q_k(0) =
 d_t{q}_k(0) = d_t^2{q}_k(0) = \cdots = d_t^5 q_k(0) = 0$.
Then, from Proposition 2.1 and since $[0,T] \times \Gamma$ is a three dimensional smooth
manifold, observation data $\phi$ is in $H^4((0,T)\times\Gamma) \subset C^2((0,T)\times\Gamma)$.
For the reconstruction of moving point sources, we choose the following
five sequences of functions in $\mathcal W$:
\begin{align}
  f_{n,\varepsilon}(t, \boldsymbol{r};\tau) =
 & \displaystyle (x+\mathrm{i}y)^n
   \eta_\varepsilon \left(t - \left(\tau - \frac{z}{c}
			      \right)\right),\quad n = 0,1,2,\cdots,  \label{eq:def_f_n_e}\\
  g_{n,\varepsilon}(t, \boldsymbol{r};\tau) =
& \displaystyle -\partial_t
   f_{n,\varepsilon}(t, \boldsymbol{r};\tau), \quad n = 0,1,2,\cdots, \label{eq:def_g_n_e}\\
  h_{n,\varepsilon}(t, \boldsymbol{r};\tau) =
& \displaystyle z (\partial_x -
  \mathrm{i}\partial_y)f_{n,\varepsilon}(t, \boldsymbol{r};\tau) - (x -
  \mathrm{i}y)\partial_z f_{n,\varepsilon}(t, \boldsymbol{r};\tau),
  \nonumber \\
   & \hspace{10ex} n = 1,2,\cdots, \label{eq:def_h_n_e}\\
  i_{n,\varepsilon}(t, \boldsymbol{r};\tau) =
& \displaystyle \partial_t^2
   f_{n,\varepsilon}(t, \boldsymbol{r};\tau), \quad n = 0,1,2,\cdots, \label{eq:def_i_n_e}\\
  j_{n,\varepsilon}(t, \boldsymbol{r};\tau) =
& \displaystyle -\partial_t 
   h_{n,\varepsilon}(t, \boldsymbol{r};\tau), \quad n = 1,2,\cdots, \label{eq:def_j_n_e}
\end{align}
where $\tau \in I_\varepsilon
\equiv \left[\max\left\{0,\sup_{\boldsymbol{r}\in\Omega}z/c
+ \varepsilon\right\},\, \min\left\{T,T+
\inf_{\boldsymbol{r}\in\Omega}z/c - \varepsilon\right\}\right] \subset [0,T]$, $0 <
\varepsilon \ll 1 $, and $\eta_\varepsilon \in
C^{\infty}({\mathbb R};{\mathbb R})$ denotes the standard mollifier
function with the support $[-\varepsilon, \varepsilon]$ (e.\@g.\@ Appendix C in \cite{Evans1}).
We note that sequences $\{f_{n,\varepsilon}\},\
\{g_{n,\varepsilon}\}$ and $\{h_{n,\varepsilon}\}$ have already applied 
to the reconstruction of fixed point sources \cite{Ohe1}.
Supplemental sequences $\{i_{n,\varepsilon}\}$ and
$\{j_{n,\varepsilon}\}$ are used to treat the effect of moving
velocities of sources.
\par
%
Due to the assumptions for $\boldsymbol{p}_k$ and $q_k$, the observation data
$\phi$ is in $C^2([0,T]\times\Gamma)$.
Then reciprocity gap functionals 
${\mathcal R}_\phi(f_{n,\varepsilon}),\ {\mathcal
R}_\phi(g_{n,\varepsilon}),\cdots, {\mathcal
R}_\phi(j_{n,\varepsilon})$
converge as $\varepsilon \rightarrow +0$. 
Now, let 
\begin{align*}
   {\mathcal R}_\phi(f_n)(\tau) & \equiv \displaystyle \lim_{\varepsilon \rightarrow +0}
   {\mathcal R}_\phi(f_{n,\varepsilon}(\cdot,\cdot;\tau)), \\
   {\mathcal R}_\phi(g_n)(\tau) & \equiv \displaystyle \lim_{\varepsilon \rightarrow +0}
   {\mathcal R}_\phi(g_{n,\varepsilon}(\cdot,\cdot;\tau)), \\
   {\mathcal R}_\phi(h_n)(\tau) &\equiv \displaystyle \lim_{\varepsilon \rightarrow +0}
   {\mathcal R}_\phi(h_{n,\varepsilon}(\cdot,\cdot;\tau)), \\
   {\mathcal R}_\phi(i_n)(\tau) &\equiv \displaystyle \lim_{\varepsilon \rightarrow +0}
   {\mathcal R}_\phi(i_{n,\varepsilon}(\cdot,\cdot;\tau)), \\
   {\mathcal R}_\phi(j_n)(\tau) &\equiv \displaystyle \lim_{\varepsilon \rightarrow +0}
   {\mathcal R}_\phi(j_{n,\varepsilon}(\cdot,\cdot;\tau)).
\end{align*}
Then we establish
\begin{align}
   {\mathcal R}_\phi(f_n)(\tau) 
=& \displaystyle -\lim_{\varepsilon \rightarrow +0} \int_0^T \int_\Gamma
  \phi(t, \boldsymbol{r}) f_{n,\varepsilon}(t, \boldsymbol{r};\tau) dS(\boldsymbol{r})dt
  \nonumber \\
=& \displaystyle -\lim_{\varepsilon \rightarrow +0} \int_\Gamma (x+
 \mathrm{i}y)^n \int_0^T \phi(t, \boldsymbol{r}) \eta_\varepsilon\left(t-
							   \left(\tau -
							    \frac{z}{c}\right)
							  \right) dt dS(\boldsymbol{r})
  \nonumber \\
=& \displaystyle - \int_\Gamma (x + \mathrm{i}y)^n \phi\left( \tau -
					    \frac{z}{c}, \boldsymbol{r} \right)
 dS(\boldsymbol{r}),
 \label{eq:expression_of_R_f_n}
\end{align}
and similarly,
\begin{align}
 {\mathcal R}_\phi(g_n) (\tau) 
  =&  -\int_\Gamma (x + \mathrm{i}y)^n \partial_t \phi\left(\tau -
							\frac{z}{c},
							\boldsymbol{r} \right)
 dS(\boldsymbol{r}) \nonumber \\
 =& \displaystyle -d_\tau \int_\Gamma (x + \mathrm{i}y)^n \phi\left( \tau
 - \frac{z}{c}, \boldsymbol{r} \right) dS(\boldsymbol{r}) \nonumber \\
 =& d_\tau {\mathcal R}_\phi(f_n) (\tau), 
  \label{eq:expression_of_R_g_n} \\
 {\mathcal R}_\phi(h_n)(\tau) 
=& \displaystyle 
  - \int_\Gamma 2nz(x + \mathrm{i}y)^{n-1} \phi\left(\tau -
					    \frac{z}{c}, \boldsymbol{r} \right) dS(\boldsymbol{r})
  \nonumber \\
 & \displaystyle \hspace{5ex}  - \frac{1}{c}\int_\Gamma (x-\mathrm{i}y)(x + \mathrm{i}y)^{n} \partial_t \phi\left(\tau -
					    \frac{z}{c}, \boldsymbol{r} \right)
 dS(\boldsymbol{r}) \nonumber \\
 =& \displaystyle 
  - \int_\Gamma 2nz(x + \mathrm{i}y)^{n-1} \phi\left(\tau - \frac{z}{c}, \boldsymbol{r} \right) dS(\boldsymbol{r})  \nonumber \\
  & \displaystyle \hspace{5ex}  
 - \frac{1}{c} d_\tau \int_\Gamma (x-\mathrm{i}y)(x + \mathrm{i}y)^{n} \phi\left(\tau -  \frac{z}{c}, \boldsymbol{r} \right) dS(\boldsymbol{r}), 
  \label{eq:expression_of_R_h_n} \\
 {\mathcal R}_\phi(i_n)(\tau) 
  =& - \int_\Gamma (x + \mathrm{i}y)^n \partial_t^2 \phi\left(\tau - \frac{z}{c}, \boldsymbol{r} \right)
  dS(\boldsymbol{r})  \nonumber \\
 =& d_\tau^2 {\mathcal R}_\phi(f_n) (\tau), 
   \label{eq:expression_of_R_i_n} \\
 {\mathcal R}_\phi(j_n)(\tau) 
=& \displaystyle
  - \int_\Gamma 2nz(x + \mathrm{i}y)^{n-1} \partial_t \phi\left(\tau -
					    \frac{z}{c}, \boldsymbol{r} \right) dS(\boldsymbol{r})
  \nonumber \\
 & \displaystyle \hspace{5ex}  - \frac{1}{c}\int_\Gamma (x-\mathrm{i}y)(x + \mathrm{i}y)^{n} \partial_t^2 \phi\left( \tau -
					    \frac{z}{c},
 \boldsymbol{r}\right)dS(\boldsymbol{r}) 
   \nonumber \\
 =& d_\tau {\mathcal R}_\phi(h_n) (\tau), 
 \label{eq:expression_of_R_j_n}
 \end{align}
since the observation data $\phi$ is in $C^2((0,T)\times\Gamma)$, 
where $d_\tau$ denotes the derivative with respect to $\tau$.
\par
Also since ${\mathcal F}(f_{{n,\varepsilon}}),\ {\mathcal
F}(g_{{n,\varepsilon}}),\cdots,{\mathcal F}(i_{{n,\varepsilon}})$
converge as $\varepsilon \rightarrow +0$, we obtain the following
formulae that express explicit relations to the parameters of moving point sources:
\begin{align}
  {\mathcal F}(f_n)(\tau) 
\equiv& \displaystyle \lim_{\varepsilon \rightarrow  +0}
  {\mathcal F}(f_{n,\varepsilon}(\cdot,\cdot;\tau)) \nonumber \\
=& \displaystyle \lim_{\varepsilon \rightarrow +0}
  \sum_{k=1}^K \int_0^T q_k(t)
  f_{n,\varepsilon}(t,\boldsymbol{p}_k(t);\tau) \nonumber \\
=& \displaystyle \lim_{\varepsilon \rightarrow +0}
  \sum_{k=1}^K \int_0^T q_k(t)
  (p_{k,x}(t) + \mathrm{i} p_{k,y}(t))^n \eta_{\varepsilon}\left(t - \tau +
							\frac{p_{k,z}(t)}{c}\right)
  dt \nonumber \\
=& \displaystyle \sum_{k=1}^K  q_k(t_k(\tau))
 \xi_k(t_k(\tau))\cdot(p_{k,xy}(t_k(\tau)))^n,
\label{eq:relation_RGf_parameter_point_source}
\\[2ex]
 {\mathcal F}(g_n) (\tau)
\equiv& \displaystyle \lim_{\varepsilon \rightarrow  +0}
  {\mathcal F}(g_{n,\varepsilon}(\cdot,\cdot;\tau)) \nonumber \\
 =&   
   \displaystyle \sum_{k=1}^K  d_\tau(q_k \xi_k) \cdot ({p_{k,xy}})^n 
    + n \sum_{k=1}^K  q_k \xi_k \cdot d_\tau(p_{k,xy})  \cdot
    (p_{k,xy})^{n-1},  
    \label{eq:relation_RGg_parameter_point_source} 
\\[2ex]
 {\mathcal F}(h_n) (\tau)
\equiv& \displaystyle \lim_{\varepsilon \rightarrow  +0}
  {\mathcal F}(h_{n,\varepsilon}(\cdot,\cdot;\tau)) \nonumber \\
   =
&  
   \displaystyle 
    2n \sum_{k=1}^K  q_k \xi_k \cdot p_{k,z} \cdot (p_{k,xy})^{n-1} 
   + \hat{R}_{h_n}  
\label{eq:relation_RGh_parameter_point_source} 
\end{align}
\begin{align}
%
 {\mathcal F}(i_n) (\tau)
\equiv& \displaystyle \lim_{\varepsilon \rightarrow  +0}
  {\mathcal F}(i_{n,\varepsilon}(\cdot,\cdot;\tau)) \nonumber \\
  =
& 
  \displaystyle \sum_{k=1}^K  
    d_\tau^2(q_k\xi_k) \cdot (p_{k,xy})^n 
 + n \sum_{k=1}^K  q_k \xi_k \cdot d_\tau^2(p_{k,xy})\cdot(p_{k,xy})^{n-1} \nonumber \\
 &  
  + \hat{R}_{i_n}  
\label{eq:relation_RGi_parameter_point_source} \\[2ex]
{\mathcal F}(j_n) (\tau)
\equiv& \displaystyle \lim_{\varepsilon \rightarrow  +0}
  {\mathcal F}(j_{n,\varepsilon}(\cdot,\cdot;\tau)) \nonumber \\
    =
&   \displaystyle 2n \sum_{k=1}^K  q_k \xi_k \cdot d_\tau(p_{k,z})
\cdot (p_{k,xy})^{n-1} 
     + \hat{R}_{j_n},
\label{eq:relation_RGj_parameter_point_source} 
\end{align}
where 
\begin{itemize}
\item $ p_{k,xy}(t_k(\tau)) =  p_{k,x}(t_k(\tau)) +
 \mathrm{i} p_{k,y}(t_k(\tau))$,
\item $t_k(\tau)$ is the unique solution $t$ of the equation
\begin{equation}
   t - \tau + \frac{p_{k,z}(t)}{c} = 0,
   \label{eq:def_t_k_tau}
\end{equation}
for each $k$ and $\tau$,
\item $\xi_k(\tau)$ is the derivative of $t_k(\tau)$ and can be expressed by
\begin{equation}
   \xi_k(\tau) = \frac{dt_k}{d\tau}(\tau) =\left(1 - \frac{d_\tau
							      ({p}_{k,z}(t_k(\tau)))}{c}\right),
 \label{eq:def_xi}
\end{equation}
\item $\hat{R}_{h_n},\ \hat{R}_{i_n}$ and $\hat{R}_{j_n}$ are defined by
\begin{align}
  \hat{R}_{h_n}
=& 
   \displaystyle 
     \frac{1}{c} \sum_{k=1}^K
                  \left\{d_\tau(q_k\xi_k) \cdot \overline{p_{k,xy}} + q_k \xi_k \cdot
		   d_\tau(\overline{p_{k,xy}}) \right\} (p_{k,xy})^{n}
    \nonumber \\
 & 
    \displaystyle 
    + \frac{n}{c} \sum_{k=1}^K
      q_k \xi_k \cdot d_\tau(p_{k,xy}) \cdot \overline{p_{k,xy}}\cdot (p_{k,xy})^{n-1} 
   \label{eq:expresion_of_Rhn} 
\\[2ex]
%
 \hat{R}_{i_n} 
=&  \displaystyle 2n \sum_{k=1}^Kd_\tau(q_k \xi_k) \cdot
 d_\tau(p_{k,xy})  \cdot (p_{k,xy})^{n-1}  \nonumber \\
 & + \displaystyle n(n-1) \sum_{k=1}^K q_k \xi_k \cdot (d_\tau(p_{k,xy}))^2 \cdot (p_{k,xy})^{n-2}, 
   \label{eq:expresion_of_Rin}
\end{align}
\begin{align}
 \hat{R}_{j_n} 
=&  
  \displaystyle 2n \sum_{k=1}^K 
   \left\{d_\tau(q_k \xi_k) \cdot p_{k,xy} + (n-1) q_k\xi_k
 \cdot d_\tau(p_{k,xy})   \right\}\cdot p_{k,z}\cdot
 (p_{k,xy})^{n-2} \nonumber \\
&  \displaystyle + \frac{1}{c} \sum_{k=1}^K 
  d_\tau^2(q_k\xi_k) \cdot \overline{p_{k,xy}}\cdot
 (p_{k,xy})^n \nonumber \\
&  \displaystyle + \frac{1}{c} \sum_{k=1}^K 
   \left\{2d_\tau(q_k\xi_k)\cdot d_\tau(\overline{p_{k,xy}}) +
    q_k\xi_k \cdot d_\tau^2 (\overline{p_{k,xy}})\right\} \cdot(p_{k,xy})^n
  \nonumber \\
&  
  \displaystyle + \frac{n}{c} \sum_{k=1}^K 
  \left\{2d_\tau(q_k\xi_k) \cdot d_\tau(p_{k,xy})+
           q_k\xi_k \cdot d_\tau^2 (p_{k,xy}) \right\} \cdot
 \overline{p_{k,xy}} \cdot (p_{k,xy})^{n-1} \nonumber \\
& \displaystyle + \frac{n}{c} \sum_{k=1}^K 
 2q_k\xi_k\cdot d_\tau(p_{k,xy}) \cdot d_\tau(\overline{p_{k,xy}})\cdot (p_{k,xy})^{n-1}
   \nonumber \\
&  
  \displaystyle + \frac{n(n-1)}{c} \sum_{k=1}^K 
     q_k \xi_k \cdot (d_\tau(p_{k,xy}))^2 \cdot \overline{p_{k,xy}}\cdot
     (p_{k,xy})^{n-2}.
   \label{eq:expresion_of_Rjn}
\end{align}
\end{itemize}
 In (\ref{eq:relation_RGg_parameter_point_source})-(\ref{eq:relation_RGj_parameter_point_source}), we omit the argument
 $(t_k(\tau))$ on $q_k,\ p_{k,xy},\ p_{k,z}$ and their derivatives, and
 the argument $(\tau)$ on $\xi_k$ to simplify the expression,
e.\@g.\@ $q_k(t_k(\tau))$ is simplified as $q_k$.
 In the following, we use these notations if we do not need to
 show these arguments explicitly.
Derivations of
 (\ref{eq:relation_RGf_parameter_point_source})-(\ref{eq:relation_RGj_parameter_point_source}) 
 are given in Appendix A.
We note that equations
 (\ref{eq:relation_RGf_parameter_point_source})-(\ref{eq:relation_RGh_parameter_point_source})
are similar to the results for fixed point sources \cite{Ohe1}, but
some supplemental terms arise due to the effect of moving velocities of sources, 
e.\@g.\@ $\xi_k(\tau)$ and ${\mathrm d}_\tau
 (p_{k,xy}(t_k(\tau)))$.
\vspace{2ex}
\par
\newpage
\noindent
\textbf{Note 3.1.} 
Since functions $f_{n,\varepsilon},\ g_{n,\varepsilon},\cdots,\
j_{n,\varepsilon}$ 
have propagation property along only $z$-axis, the reciprocity gap
functionals of these functions have `retarded' property only on
$z$-coordinate.
Also we can obtain an another expression of the perturbation  factor
$\xi_k(\tau)$ due to the $z$-component of the velocity of point sources:
\begin{equation}
  \xi_k(\tau)  \displaystyle = \left(1 + \frac{d_t
 {p}_{k,z}(t_k(\tau))}{c}\right)^{-1}.
  \label{eq:another_expression_of_xi} 
\end{equation}
Equation (\ref{eq:another_expression_of_xi}) shows that 
the perturbation  factor $\xi_k$ effects similarly as the
factor $1/h$ in the extended Li\'{e}nard-Wiechert retarded potential
(\ref{eq:special_sol_point_source}) and
(\ref{eq:special_sol_dipole_source}).
\vspace{2ex}
\par
Using expressions (\ref{eq:relation_RGf_parameter_point_source})-(\ref{eq:relation_RGj_parameter_point_source}), we obtain the following
reconstruction theorem for moving point sources:
\par
\begin{thm}
Let $I_0 \equiv \left[\max\{0,\sup_{\boldsymbol{r} \in \Omega}z/c\},
 \min\{T + \inf_{\boldsymbol{r} \in \Omega}z/c,T\}\right] \subset[0,T]$.
For each $\tau \in I_0$, let $K(\tau)$ be the number of point sources such that $q_k(t_k(\tau)) \neq 0$.
Assume that $K(\tau) \leq K_M$ for given $K_M$,  and $p_{j,xy}(t_j(\tau)) \neq p_{k,xy}(t_k(\tau))$ if $j \neq k$.
Then, we can identify $K(\tau)$ from the reciprocity gap functionals ${\mathcal R}_\phi(f_n)(\tau),\ n=0,1,2,\cdots,2K_M$.
Also we can uniquely reconstruct $\boldsymbol{p}_k(t_k(\tau))$ and
 $q_k(t_k(\tau))$,  $k=1,2,\cdots,K(\tau)$ from 
\begin{itemize}
\item ${\mathcal R}_\phi(f_n)(\tau),\ n=0,1,2,\cdots,2K(\tau)-1$, 
\item ${\mathcal R}_\phi(g_n)(\tau),\ n=0,1,2,\cdots,2K(\tau)-1$,
\item ${\mathcal R}_\phi(h_n)(\tau),\ n=1,2,3,\cdots,K(\tau)$,
\item ${\mathcal R}_\phi(i_n)(\tau),\ n=0,1,2,\cdots,2K(\tau)-1$,
\item ${\mathcal R}_\phi(j_n)(\tau),\ n=1,2,3,\cdots,K(\tau)$.
\end{itemize}
\end{thm}
\vspace{2ex}
\par
\noindent
\begin{proof}
We show the proof of the theorem in the following five steps which
 describe the procedure of the reconstruction process.
\vspace{1ex}
\begin{description}
\item[Step 1.] Identify $K(\tau)$ from ${\mathcal R}_\phi(f_n)(\tau),\ n=0,1,2,\cdots, 2K_M$.
\item[Step 2.] Reconstruct $p_{k,xy}(t_k(\tau)),\ k=1,2,\cdots,K(\tau)$, and identify perturbed magnitudes $q_k(t_k(\tau))\xi_k(\tau)$
	   from ${\mathcal R}_\phi(f_n)(\tau),\ n=0,1,2,\cdots,2K(\tau)-1$.
\item[Step 3.] Reconstruct $p_{k,z}(t_k(\tau)),\ k=1,2,\cdots,K(\tau)$ from ${\mathcal R}_\phi(g_n)(\tau),\ n=0,1,2,\cdots,2K(\tau)-1$ and ${\mathcal
	   R}_\phi(h_n)(\tau),\ n=1,2,3,\cdots,K(\tau)$.
\item[Step 4.] Identify $d_\tau(p_{k,z}(t_k(\tau))),\ k=1,2,\cdots,K(\tau)$ from ${\mathcal R}_\phi(i_n)(\tau),\ n=1,2,3,\cdots,2K(\tau)$ and ${\mathcal
	   R}_\phi(j_n)(\tau),\ n=1,2,3,\cdots,K(\tau)$.
\item[Step 5.] Compute $\xi_k(\tau)$ using 
	   $d_\tau({p}_{k,z}(t_k(\tau)))$ for each $k$, and
	   reconstruct the magnitude $q_k(t_k(\tau))$ 
	   from its perturbed value $q_k(t_k(\tau))\xi_k(\tau)$.
\end{description}
\vspace{1ex}
We show the detail of each step bellow.
\vspace{1ex}
\par
\noindent
\textbf{Step 1.} Define $L \times L$ Hankel matrix
\begin{equation*}
  H_{L,\mu}(\tau) = \left(
        \begin{array}{cccc}
          {\mathcal R}_\phi(f_{\mu})(\tau) &
	  {\mathcal R}_\phi(f_{\mu+1})(\tau) & 
           \cdots & 
          {\mathcal R}_\phi(f_{\mu+L-1})(\tau) \\
          {\mathcal R}_\phi(f_{\mu+1})(\tau) & 
          {\mathcal R}_\phi(f_{\mu+2})(\tau) &
	   \cdots & 
          {\mathcal R}_\phi(f_{\mu+L})(\tau) \\
          \vdots & \vdots &  \ddots & \vdots \\
          {\mathcal R}_\phi(f_{\mu+L-1})(\tau) & 
          {\mathcal R}_\phi(f_{\mu+L})(\tau) &
	  \cdots & 
          {\mathcal R}_\phi(f_{\mu+2L-2})(\tau)
     \end{array}
\right).
\end{equation*}
Then, from (\ref{eq:relation_RGf_parameter_point_source}) and using
corollary 3 in \cite{Nara1}, we
can determine $K(\tau)$ by
\begin{equation}
  K(\tau) = \max\left\{ L \, |\, \det H_{L,0}(\tau) \neq 0\right\}.
  \label{eq:determine_K}
\end{equation}
Equation (\ref{eq:determine_K}) shows that $\det H_{L,0}(\tau)$ vanishes for $L \geq K(\tau)+1$.
Hence we can determine $K(\tau)$ from ${\mathcal R}_\phi(f_n)(\tau)$ for
$n=0,1,2,\cdots,2K_M$ if $K(\tau) \leq K_M$
\vspace{1ex}
\par
\noindent
\textbf{Step 2.}
From the definition of $H_{L,\mu}(\tau)$ and using Theorem 2 in \cite{Badia1}, we can reconstruct $p_{k,xy}(t_k(\tau)),\
k=1,2,\cdots,K(\tau)$  as eigenvalues of the matrix $(H_{K(\tau),0}(\tau))^{-1} H_{K(\tau),1}(\tau)$.
For this computation, we need $\mathcal{R}_\phi(f_n)(\tau)$ for $n=0,1,2,\cdots,2K(\tau)-1$.
\par
Now, let us define $K(\tau) \times K(\tau)$-matrix
\begin{equation*}
 V(\tau) = \left(\begin{array}{cccc}
 1 & 1 & \cdots & 1 \\
 p_{1,xy}(t_k(\tau))   & p_{2,xy}(t_k(\tau))   & \cdots & p_{K(\tau),xy}(t_k(\tau))   \\ 
 (p_{1,xy}(t_k(\tau)))^2 & (p_{2,xy}(t_k(\tau)))^2 & \cdots & (p_{K(\tau),xy}(t_k(\tau)))^2 \\
   \vdots   & \vdots     & \ddots & \vdots           \\
 (p_{1,xy}(t_k(\tau)))^{K(\tau)-1} & (p_{2,xy}(t_k(\tau)))^{K(\tau)-1} & \cdots &
  (p_{K(\tau),xy}(t_k(\tau)))^{K(\tau)-1} 
 \end{array} 
 \right),
\end{equation*}
and ${\bm \zeta}(\tau),{\bm \beta}(\tau) \in {\mathbb{C}}^{K(\tau)}$  by
\begin{equation*}
  {\bm \zeta}(\tau) = \left( \begin{array}{c}
     q_1(t_k(\tau)) \xi_1(\tau) \\ q_2(t_k(\tau)) \xi_2(\tau) \\ \vdots \\ q_{K(\tau)}(t_k(\tau)) \xi_{K(\tau)}(\tau)
     \end{array}
  \right),
\quad
 {\bm \beta}(\tau) = \left( \begin{array}{c}
     {\mathcal R}_\phi(f_0)(\tau) \\ {\mathcal R}_\phi(f_1)(\tau) \\
	      \vdots \\ {\mathcal R}_\phi(f_{K(\tau)-1})(\tau) 
\end{array}
\right).
\end{equation*}
Then, equations (\ref{eq:relation_RGf_parameter_point_source}) for
$n=0,1,2,\cdots,K(\tau)-1$ are rewritten by
\begin{equation}
   V(\tau) {\bm \zeta}(\tau) = {\bm \beta}(\tau).
   \label{eq:eq_for_perturbed_lambda}
\end{equation}
Since $V(\tau)$ is a Vandermonde type matrix and from the assumption
that $p_{j,xy}(t_j(\tau)) \neq
p_{k,xy}(t_k(\tau))$ for $j \neq k$, we derive $\det V(\tau) \neq 0$ .
Then, the equation (\ref{eq:eq_for_perturbed_lambda}) is
uniquely solvable, 
and we can identify the perturbed magnitudes
$q_k(t_k(\tau))\,\xi_k(\tau),\ k=1,2,\cdots,K(\tau)$.
\par
\vspace{2ex}
\par
\noindent
\textbf{Step 3.}
Owing to equation (\ref{eq:relation_RGh_parameter_point_source}), if we know values of
${\mathcal R}_\phi(h_n)(\tau)$ and $\hat{R}_{h_n}$ for
$n=1,2,\cdots,K(\tau)$, we can identify $q_k\,\xi_k \cdot p_{k,z},\
k=1,2,\cdots,K(\tau)$ as the unique solution of the system of linear
equation:
\begin{equation}
 \sum_{k=1}^{K(\tau)}(q_k \xi_k \cdot p_{k,z} ) \cdot  (p_{k,xy})^{n-1} 
= \frac{1}{2n}\left({\mathcal R}(h_n)(\tau) -
					 \hat{R}_{h_n}\right), 
\quad  n=1,2,\cdots,K(\tau).
  \label{eq:equation_in_Step3}
\end{equation}
Then, dividing each solution $q_k \xi_k \cdot p_{k,z}$ of (\ref{eq:equation_in_Step3}) by the perturbed
magnitude $q_k \xi_k$ that is identified in step 2,
we obtain $p_{k,z}$.
Hence we consider the estimation of $\hat{R}_{h_n}$.
\par
In the expression (\ref{eq:expresion_of_Rhn}) of $\hat{R}_{h_n}$,
 $d_\tau(q_k \xi_k)$ and $d_\tau(p_{k,xy})$ remain unknowns.
We apply (\ref{eq:relation_RGg_parameter_point_source})
 to identify these unknowns. 
Let us define ${\bm \psi}_k,\ {\bm \psi}'_k \in {\mathbb{C}}^{2K(\tau)}$ by
\begin{equation*}
{\bm \psi}_k = \left(
      \begin{array}{c}
      1 \\ 
      p_{k,xy} \\ 
     (p_{k,xy})^2 \\
     (p_{k,xy})^3 \\ 
      \vdots \\ 
     (p_{k,xy})^{2K(\tau)-1}
      \end{array}
\right), \qquad
{\boldsymbol{\psi}}'_k = \left(
      \begin{array}{c}
      0 \\ 1 \\ 
      2p_{k,xy} \\ 
      3(p_{k,xy})^2 \\ 
      \vdots \\ 
      (2K(\tau)-1)(p_{k,xy})^{2K(\tau)-2}
      \end{array}
\right)
\end{equation*}
and the matrix $\tilde{V}(\tau) \in {\mathbb C}^{2K(\tau) \times 2K(\tau)}$ by
\begin{equation}
  \tilde{V}(\tau) = \left(
        \begin{array}{cccccccc}
        {\bm \psi}_1  & {\bm \psi}_2  & \cdots & {\bm \psi}_{K(\tau)} &
        {\bm \psi}'_1 & {\bm \psi}'_2 & \cdots & {\bm \psi}'_{K(\tau)}
\end{array}
\right).
\label{eq:def_tildeV}
\end{equation}
Also define ${\bm \iota}(\tau),\ {\bm \gamma}(\tau) \in {\mathbb{C}}^{2K(\tau)}$
 by
\begin{equation*}
   \hspace{-8ex}
  {\bm \iota}{(\tau)} = \left( \begin{array}{c}
    d_\tau(q_1\xi_1) \\ d_\tau(q_2\xi_2) \\ \vdots \\ d_\tau(q_{K(\tau)} \xi_{K(\tau)}))
      \\[0.5ex]
     q_1\xi_1 \cdot d_\tau(p_{1,xy})    \\ q_2\xi_2 \cdot d_\tau(p_{2,xy})    \\ 
     \vdots \\ 
     q_{K(\tau)}\xi_{K(\tau)} \cdot d_\tau(p_{K(\tau),xy})
     \end{array}
  \right),
\quad
  {\bm \gamma}{(\tau)} = \left( \begin{array}{c}
     {\mathcal R}_\phi(g_0)(\tau) \\ {\mathcal R}_\phi(g_1)(\tau) \\
	      \vdots \\ {\mathcal R}_\phi(g_{2K(\tau)-1})(\tau) 
\end{array}
\right).
\end{equation*}
Then, equations (\ref{eq:relation_RGg_parameter_point_source}) for
  $n=0,1,2,\cdots,2K(\tau)-1$ are rewritten as
\begin{equation}
   \tilde{V}(\tau) {\bm \iota}(\tau) = {\bm \gamma}(\tau).
   \label{eq:eq_of_LAMBDA}
\end{equation}
The matrix $\tilde{V}(\tau)$ is a generalised Vandermonde-type matrix,
and we establish the following equality:
\begin{equation}
 \det \tilde{V}(\tau) = (-1)^{K(\tau)(K(\tau)-1)/2}
  \prod_{j>k}(p_{j,xy}(t_j(\tau))-p_{k,xy}(t_k(\tau)))^4.
  \label{eq:detV}
\end{equation}
A proof of (\ref{eq:detV}) is given in Appendix B.
From the assumption that $p_{j,xy}(t_j(\tau)) \neq p_{k,xy}(t_k(\tau))$ if
$j \neq k$, we obtain $\det \tilde{V}(\tau) \neq 0$.
Then,  the equation (\ref{eq:eq_of_LAMBDA}) is uniquely solvable, and  we can identify $d_\tau(q_k\,\xi_k)$ and
$d_\tau(p_{k,xy}),\ k=1,2,\cdots,K(\tau)$.
This concludes that we can reconstruct $p_{k,z}(t_k(\tau)),\
k=1,2,\cdots,K(\tau)$.
\par
\vspace{2ex}
\par
\noindent
\textbf{Step 4.}
Owing to (\ref{eq:relation_RGj_parameter_point_source}), we can identify
 $q_k \xi_k \cdot d_\tau(p_{k,z})$ as the unique solution of
 the following system of linear
equation if we know $\mathcal{R}_\phi(j_n)(\tau)$ and $\hat{R}_{j_n}$
for $n=1,2,\cdots,K(\tau)$:
\begin{equation}
 \sum_{k=1}^{K(\tau)}(q_k \xi_k \cdot d_\tau(p_{k,z})) \cdot
  (p_{k,xy})^{n-1} = \displaystyle \frac{1}{2n}\left({\mathcal
						       R}_\phi(j_n)(\tau) -
						       \hat{R}_{j_n}\right), 
   \quad  n=1,2,\cdots,K(\tau).
 \label{eq:equation_for_dpz}
\end{equation}
From the expression (\ref{eq:expresion_of_Rjn}) of $\hat{R}_{j_n}$, we
 can see that $d_\tau^2(q_k \xi_k)$ and $d^2_\tau(p_{k,xy})$ are unknown in the right hand side term of (\ref{eq:equation_for_dpz}).
Using the same idea as Step 3, these unknowns can be identified from
${\mathcal R}_\phi(i_n)(\tau),\ n=0,1,2,\cdots,2K(\tau)-1$ since
 ${\mathcal R}_\phi(i_n)(\tau)$ and $\hat{R}_{i_n}$ are given by
 (\ref{eq:relation_RGi_parameter_point_source}) and
 (\ref{eq:expresion_of_Rin}), respectively.
Then, we have identified all unknowns in the right hand side terms of
(\ref{eq:equation_for_dpz}), and we can reconstruct $d_\tau(p_{k,z}),\
k=1,2,\cdots,K(\tau)$ similarly as Step 3.
\par
\vspace{2ex}
\par
\noindent
\textbf{Step 5.} 
Since we have identified $d_\tau(p_{k,z})$ in Step 4, we can compute the
 perturbation term $\xi_k$ by (\ref{eq:def_xi}), and reconstruct magnitude $q_k$ of each
point source from perturbed magnitude $q_k\,\xi_k$ for $k=1,2,\cdots,K(\tau)$.
\vspace{2ex}
\par
Then, we have reconstructed all parameters of moving point
sources at $t = t_k(\tau),\ k=1,2,\cdots,K(\tau)$, and this completes  the proof of Theorem 3.1.
\par
\end{proof}
\vspace{2ex}
\par
\noindent
\textbf{Note 3.2.} 
 We can use alternative functions that propagate to positive direction with
 respect to $z$-coordinate, i.\@e.\@
\begin{equation*}
  \tilde{f}_{n,\varepsilon}(t, \boldsymbol{r};\tau) =
   (x+\mathrm{i}y)^n
   \eta_\varepsilon \left(t - \left(\tau + \frac{z}{c}
			      \right)\right),\quad n = 0,1,2,\cdots.
\end{equation*}
In this case, $t_k(\tau)$ is replaced by the solution of $t$ of the equation
\begin{equation*}
   t - \tau - \frac{p_{k,z}(t)}{c} = 0,
\end{equation*}
and $\xi_k(\tau)$ is replaced by
\begin{equation*}
   \xi_k(\tau) = \frac{dt_k}{d\tau}(\tau) =\left(1 + \frac{d_\tau
							      ({p}_{k,z}(t_k(\tau)))}{c}\right).
\end{equation*}
\par
\vspace{2ex}
\noindent
\textbf{Note 3.3.} 
In practical cases, it is difficult to apply the condition (\ref{eq:determine_K}) to the
identification of $K(\tau)$
because of noise in the observation data and errors of numerical computation.
In section 4, we propose another heuristic criterion to identify
the number of sources. 
\par
\vspace{2ex}
\par
\noindent
\textbf{Note 3.4.} 
In our reconstruction procedure, we reconstruct the locations and
magnitude of points sources 
at time $t_k(\tau)$ for each $k$ and $\tau$.
We can easily estimate $t_k(\tau)$ for each $k$ and $\tau$ by replacing
$t$ of (\ref{eq:def_t_k_tau}) by $t_k(\tau)$, i.\@e.\@
\begin{displaymath}
     t_k(\tau) = \tau - \frac{p_{k,z}(t_k(\tau))}{c}.
\end{displaymath}
Here, we do not need to solve this equation with respect to $t_k(\tau)$
but only substitute $p_{k,z}(t_k(\tau))$ that was identified in Step 3.
\par
\vspace{2ex}
\par
\noindent
\noindent
\textbf{Note 3.5.} 
Since we have obtained $p_{k,z}(t_k(\tau))$ in Step 3, we may compute 
$d_\tau p_{k,z}(t_k(\tau))$ using numerical differentiation in
stead of Step 4.
(Here, we call \textbf{Step 4'}.)
We compare these two methods in section 4.
\par
\vspace{2ex}
\par
\newpage
\noindent
\textbf{Note 3.6.} 
The merits of our reconstruction procedure are (a) real-timeness
and (b) independence on the initial condition.
Reconstruction methods based on an optimisation procedure such as the
least squares method usually need whole time observation data and the initial
condition because we need to solve the initial-boundary value problem.
On the other hand, the proof of Theorem 3.1 shows that 
our method needs observations in the interval $[\tau -\sup_{\boldsymbol{r}\in\Omega}z/c,\
\tau-\inf_{\boldsymbol{r}\in\Omega}z/c]$
to reconstruct unknown parameters at $t = t_k(\tau)$.
%
This shows the real-timeness and independence on the initial
condition of our reconstruction procedure.
\par
\vspace{2ex}
%
%
\subsection{Reconstruction of moving dipole sources}
Next, we consider the reconstruction of moving dipole sources.
Here, we assume that $\boldsymbol{p}_k \in C^{6}([0,T]; D),\ 
\boldsymbol{m}_k \in C^6([0,T];\mathbb{R}^3)$ and $\boldsymbol{m}_k(0) =
 d_t{\boldsymbol{m}}_k(0) = d_t^2{\boldsymbol{m}}_k(0) = \cdots = d_t^6\boldsymbol{m}_k(0)=
 \boldsymbol{0}$, then the observation data $\phi$ is in $H^4([0,T]\times\Gamma)
 \subset C^2([0,T]\times\Gamma)$.
We also add an assumption that $m_{k,z}(t) \equiv 0$ for all $k$,
i.\@e.\@ the moment of each dipole source is expressed by
$\boldsymbol{m}_k(t)= (m_{k,x}(t),\,m_{k,y}(t), 0)$. 
For the reconstruction of moving dipole sources, we apply the same five sequences of functions $\{f_{n,\varepsilon}\}, \{g_{n,\varepsilon}\}, \cdots,
\{j_{n,\varepsilon}\}$  as for the reconstruction of moving
point sources.
Then, we obtain the following expressions of ${\mathcal F}(f_n)(\tau)$, ${\mathcal F}(g_n)(\tau)$, $\cdots$,
${\mathcal F}(j_n)(\tau)$ using the parameters of moving dipole sources:
\begin{align}
{\mathcal F}(f_n)(\tau) 
=& 
 \displaystyle n \sum_{k=1}^K m_{k,xy} \xi_k \cdot (p_{k,xy})^{n-1}, \quad
 n=1,2,\cdots, \label{eq:relation_RGf_parameter_dipole_source}\\
%
{\mathcal F}(g_n)(\tau)
=& \displaystyle 
    n \sum_{k=1}^K d_\tau(m_{k,xy}\xi_k) \cdot
    (p_{k,xy})^{n-1} \nonumber \\
 &    +n(n-1) \sum_{k=1}^K  m_{k,xy} \xi_k \cdot d_\tau (p_{k,xy})  \cdot
 (p_{k,xy})^{n-2}, \quad n =1,2,\cdots,
 \label{eq:relation_RGg_parameter_dipole_source}
 \\
%
{\mathcal F}(h_n)(\tau)
=& \displaystyle 
    2n(n-1) \sum_{k=1}^K  m_{k,xy} \xi_k \cdot p_{k,z} \cdot
    (p_{k,xy})^{n-2} 
   + \hat{R}_{h_n},
    \quad n=2,3,\cdots,
 \label{eq:relation_RGh_parameter_dipole_source}
\\
%
{\mathcal F}(i_n)(\tau)
=& \displaystyle 
   n \sum_{k=1}^K  
  d_\tau^2 (m_{k,xy} \xi_k)  \cdot (p_{k,xy})^{n-1} 
  + n(n-1) \sum_{k=1}^K m_{k,xy}\xi_k \cdot d^2_\tau (p_{k,xy})
 \cdot (p_{k,xy})^{n-2}
  \nonumber \\
 & + \hat{R}_{i_n},
   \quad n=1,2,\cdots, \label{eq:relation_RGi_parameter_dipole_source}\\
  {\mathcal F}(j_n)(\tau)
=& \displaystyle 
   2n(n-1) \sum_{k=1}^K  m_{k,xy} \xi_k \cdot d_\tau (p_{k,z})  \cdot (p_{k,xy})^{n-2} 
    + \hat{R}_{j_n}, \quad
 n=2,3,\cdots,\label{eq:relation_RGj_parameter_dipole_source}
\end{align}
where 
\begin{itemize}
\item $ m_{k,xy} \equiv  m_{k,x}(t_k(\tau)) +
 \mathrm{i} m_{k,y}(t_k(\tau))$,
\item $\hat{R}_{h_n}$, $\hat{R}_{i_n}$, and $\hat{R}_{j_n}$ are expressed by
\begin{align*}
 \hat{R}_{h_n} =
 & \displaystyle 
    \frac{1}{c} \sum_{k=1}^K 
      d_\tau(\overline{m_{k,xy}}\xi_k) \cdot (p_{k,xy})^{n} \nonumber \\
 & \displaystyle 
    + \frac{n}{c} \sum_{k=1}^K
                   \left\{d_\tau(m_{k,xy}\xi_k) \cdot \overline{p_{k,xy}} 
                         + \overline{m_{k,xy}}\xi_k \cdot d_\tau (p_{k,xy})  
                         + m_{k,xy}\xi_k \cdot d_\tau (\overline{p_{k,xy}})  \right\} (p_{k,xy})^{n-1}
    \nonumber \\
 & \displaystyle 
    + \frac{n(n-1)}{c} \sum_{k=1}^K
      m_{k,xy} \xi_k \cdot d_\tau (p_{k,xy}) \cdot \overline{p_{k,xy}}
      \cdot (p_{k,xy})^{n-2}, \nonumber
  \\
 \hat{R}_{i_n} =
   & \displaystyle 
  2n(n-1) \sum_{k=1}^K  
    d_\tau(m_{k,xy}\xi_k) \cdot d_\tau (p_{k,xy}) \cdot (p_{k,xy})^{n-2} 
    \nonumber \\
&  \displaystyle 
  + n(n-1)(n-2) \sum_{k=1}^K
      m_{k,xy} \xi_k \cdot (d_\tau (p_{k,xy}))^2 \cdot
 (p_{k,xy})^{n-3},
  \nonumber
%
\end{align*}
\begin{align*}
 \hat{R}_{j_n} =
&  
  \displaystyle 2n(n-1) \sum_{k=1}^K 
   \left\{ d_\tau (m_{k,xy} \xi_k) \cdot p_{k,xy} + 
    (n-2) m_{k,xy}\xi_k \cdot d_\tau (p_{k,xy})
   \right\} \cdot p_{k,z}\cdot (p_{k,xy})^{n-3} \\
&  
  \displaystyle + \frac{1}{c} \sum_{k=1}^K 
  d_\tau^2 (\overline{m_{k,xy}}\xi_k)\cdot (p_{k,xy})^n 
  + \frac{n}{c} \sum_{k=1}^K 
  d_\tau^2 (m_{k,xy}\xi_k) \cdot \overline{p_{k,xy}} \cdot (p_{k,xy})^{n-1} \\
&  
  \displaystyle + \frac{n}{c} \sum_{k=1}^K 
   \left\{2 d_\tau
    (\overline{m_{k,xy}}\xi_k)\cdot d_\tau ({p_{k,xy}}) +
    \overline{m_{k,xy}}\xi_k \cdot d_\tau^2 (p_{k,xy})\right\}\cdot
   (p_{k,xy})^{n-1}  \\
&  
  \displaystyle + \frac{n}{c} \sum_{k=1}^K 
   \left\{2 d_\tau (m_{k,xy}\xi_k) \cdot d_\tau (\overline{p_{k,xy}}) +
          m_{k,xy}\xi_k \cdot d_\tau^2 (\overline{p_{k,xy}}) \right\}\cdot
   (p_{k,xy})^{n-1}  \\
&  
  \displaystyle + \frac{n(n-1)}{c} \sum_{k=1}^K 
  \overline{m_{k,xy}}\xi_k\cdot(d_\tau (p_{k,xy}))^2 \cdot(p_{k,xy})^{n-2} \\
&  
  \displaystyle + \frac{n(n-1)}{c} \sum_{k=1}^K 
  \left\{ 2 d_\tau (m_{k,xy}\xi_k) \cdot d_\tau (p_{k,xy})
   + m_{k,xy}\xi_k \cdot d_\tau^2 (p_{k,xy})
  \right\} \cdot\overline{p_{k,xy}} \cdot (p_{k,xy})^{n-2}\\
& \displaystyle + \frac{n(n-1)}{c} \sum_{k=1}^K 
    2 m_{k,xy}\xi_k    \cdot d_\tau
   (\overline{p_{k,xy}}) \cdot d_\tau (p_{k,xy})\cdot
    (p_{k,xy})^{n-2} \\
&  
  \displaystyle + \frac{n(n-1)(n-2)}{c} \sum_{k=1}^K 
     m_{k,xy} \xi_k \cdot \overline{p_{k,xy}} \cdot (d_\tau (p_{k,xy}))^2 \cdot
    ( p_{k,xy})^{n-3}.
\end{align*}
\item $t_k(\tau)$ and $\xi_k(\tau)$ are defined by
      (\ref{eq:def_t_k_tau}) and (\ref{eq:def_xi}), respectively.
\end{itemize}
Derivations of equations (\ref{eq:relation_RGf_parameter_dipole_source})-(\ref{eq:relation_RGj_parameter_dipole_source}) are also given in Appendix A.
\par
Similarly to the case for moving point sources, 
we can establish the following reconstruction theorem for moving dipole sources:
\newpage
\begin{thm}
Let $I_0$ be the same interval as in Theorem 3.1. 
 For each $\tau \in I_0$, let $K(\tau)$ be the number of dipole sources such that $\boldsymbol{m}_k(t_k(\tau)) \neq \boldsymbol{0}$.
Assume that $K(\tau) \leq K_M$ for given $K_M$, $p_{j,xy}(t_j(\tau)) \neq p_{k,xy}(t_k(\tau))$ if $j \neq k$, and $m_{k,z}(t) \equiv 0$.
Then, we can identify $K(\tau)$ from the reciprocity gap functionals
 ${\mathcal R}_\phi(f_n)(\tau),\ n=1,2,\cdots,2K_M+1$.
Also we can uniquely reconstruct  $\boldsymbol{p}_k(t_k(\tau))$ and $\boldsymbol{m}_k(t_k(\tau))$ for
 $k=1,2,\cdots,K(\tau)$ from
\begin{itemize}
\item ${\mathcal R}(f_n)(\tau),\ n=1,2,3,\cdots,2K(\tau)$, 
\item ${\mathcal R}(g_n)(\tau),\ n=1,2,3,\cdots,2K(\tau)$,
\item ${\mathcal R}(h_n)(\tau),\ n=2,3,4,\cdots,K(\tau)+1$,
\item ${\mathcal R}(i_n)(\tau),\ n=1,2,3,\cdots,2K(\tau)$,
\item ${\mathcal R}(j_n)(\tau),\ n=2,3,4,\cdots,K(\tau)+1$.
\end{itemize}
\end{thm}
\vspace{2ex}
\par
\noindent
\begin{proof}
Comparing expressions 
(\ref{eq:relation_RGf_parameter_dipole_source})-(\ref{eq:relation_RGj_parameter_dipole_source})
to
(\ref{eq:relation_RGf_parameter_point_source})-(\ref{eq:relation_RGj_parameter_point_source}),
we can find similar relations between the reciprocity gap functionals
and parameters of dipole sources as for moving point sources.
Each step of the reconstruction procedure in the proof of Theorem
3.1 is modified as follows:
\par
\vspace{3ex}
\noindent
\textbf{Steps 1 and 2.}
Dividing (\ref{eq:relation_RGf_parameter_dipole_source}) by $n$, we have
\begin{equation}
 \frac{1}{n}{\mathcal R}_\phi(f_n)(\tau) = 
 \displaystyle \sum_{k=1}^K m_{k,xy}(t_k(\tau))\, \xi_k(\tau) \cdot (p_{k,xy}(t_k(\tau)))^{n-1}, \quad
 n=1,2,\cdots.
 \label{eq:relation_RGf_parameter_dipole_source_modified}
\end{equation}
The right hand side term of (\ref{eq:relation_RGf_parameter_dipole_source_modified}) has the
same form as (\ref{eq:relation_RGf_parameter_point_source}) 
by replacing $q_k(t_k(\tau))$ by $m_{k,xy}(t_k(\tau))$ and the power $n$ by $n-1$.
Hence, we can identify $K(\tau)$ using $H_{k,1}(\tau)$ computed from  ${\mathcal R}_\phi(f_n)(\tau),\
n=1,2,\cdots,2K_M+1$, and reconstruct $p_{k,xy}(t_k(\tau))$, and
$m_{k,xy}(t_k(\tau))\xi_k(\tau)$ from ${\mathcal R}_\phi(f_n)(\tau),\
n=1,2,3,\cdots,2K(\tau)$ using the same procedure as in Steps 1
and 2 of the proof of Theorem 3.1.
\vspace{1ex}
\par
\noindent
\textbf{Steps 3 and 4.} 
Similarly to Steps 1 and 2, by dividing
 (\ref{eq:relation_RGg_parameter_dipole_source})-(\ref{eq:relation_RGj_parameter_dipole_source})
 by $n$, we can obtain
the same forms as
 (\ref{eq:relation_RGg_parameter_point_source})-(\ref{eq:relation_RGj_parameter_point_source}) 
by replacing $q_k(t_k(\tau))$ by $m_{k,xy}(t_k(\tau))$ and the power $n$ by $n-1$.
Hence, we can identify $p_{k,z}(t_k(\tau))$ 
from ${\mathcal R}_\phi(g_n)(\tau),\ n=1,2,3,\cdots,2K(\tau)$ and
${\mathcal R}_\phi(h_n)(\tau),\ n=2,\cdots,K(\tau)+1$,
 and $d_\tau(p_{k,z}(t_k(\tau)))$ from 
${\mathcal R}_\phi(i_n)(\tau),\ n=1,2,3,\cdots,2K(\tau)$ and ${\mathcal
R}_\phi(j_n)(\tau),\ n=2,\cdots,K(\tau)+1$ 
using the same procedure as in Steps 3
and 4 of the proof of Theorem 3.1.
\vspace{1ex}
\par
\noindent
\textbf{Step 5.} 
Finally, we can compute $\xi_k(\tau)$ with (\ref{eq:def_xi}) using $d_\tau (p_{k,z}(t_k(\tau)))$ 
identified in Step 4, and reconstruct each moment $m_{k,xy}(t_k(\tau))$ from its perturbed value $m_{k,xy}(t_k(\tau))\,\xi_k(\tau)$.
Then, we have reconstructed
all parameters of moving dipole sources.
\par
\end{proof}
\vspace{2ex}
\par
\noindent
\textbf{Note 3.7.} 
If we remove the assumption $m_{k,z}(t) \equiv 0$, the
expression (\ref{eq:relation_RGf_parameter_dipole_source}) becomes
\begin{align}
 {\mathcal R}_\phi(f_n) (\tau)
 =&  \displaystyle n \sum_{k=1}^K m_{k,xy} \xi_k \cdot (p_{k,xy})^{n-1} \nonumber \\
  &  \displaystyle - \frac{1}{c} \sum_{k=1}^K d_\tau(m_{k,z}\xi_k) \cdot (p_{k,xy})^n 
                 - \frac{n}{c} \sum_{k=1}^K  m_{k,z} \xi_k \cdot d_\tau (p_{k,xy})  \cdot (p_{k,xy})^{n-1}.
 \label{eq:include_vertical_compnent}
\end{align}
Then, the idea in Step 1 can not  be applied.
Reconstruction method for such cases remains an open problem.
\par
\section{Numerical Experiments}

We now present some numerical experiments for the
reconstruction procedure proposed in section 3.
Throughout the experiments, we set $\Omega = \{ \boldsymbol{r}\,|\,
|\boldsymbol{r}| \le 2 \} \subset \mathbb{R}^3$, $c=1$ and $T=70$.
To give the observation data $\phi = \partial_\nu u$ on $\Gamma$,  we solve the initial-boundary value problem
(\ref{eq:mixed_problem_wave_eq}) numerically using the following boundary integral expression:
\begin{align}
  \frac{1}{2} \phi(t, \boldsymbol{r})
=& \displaystyle 
   \int_0^t \int_\Gamma \partial_{\nu(\boldsymbol{r})} G(t,\boldsymbol{r};s, \bm \rho)
    \phi(s, \bm \rho) dS(\bm \rho) ds \nonumber \\
 & \displaystyle
 + \int_0^t \int_\Omega \partial_{\nu(\boldsymbol{r})} G(t, \boldsymbol{r};s, \bm \rho)
    F(s,\bm \rho) dV(\bm \rho) ds, \nonumber \\
 =& \displaystyle 
  -\frac{1}{4\pi}\int_\Gamma
  \frac{{\bm \nu}(\boldsymbol{r})\cdot(\boldsymbol{r}-\bm \rho)}{|\boldsymbol{r} -
 \bm \rho|^3} \cdot \phi\left(t -  \frac{|\boldsymbol{r}-\bm \rho|}{c}, \bm \rho\right) dS(\bm \rho) \nonumber \\
 & \displaystyle 
  - \frac{1}{4\pi c}\int_\Gamma
 \frac{{\bm \nu}(\boldsymbol{r})\cdot(\boldsymbol{r}-\bm \rho)}{|\boldsymbol{r} -  \bm \rho|^2} \cdot \partial_t \phi\left(t -
  \frac{|\boldsymbol{r}-\bm \rho|}{c}, \bm \rho\right) dS(\bm \rho) \nonumber\\
 & \displaystyle
  + \partial_{\nu(\boldsymbol{r})} u_{\mathrm N}(t,\boldsymbol{r}), \qquad (t,\boldsymbol{r})\in
 [0,T]\times\Gamma
\label{eq:boundary_integral_equation}
\end{align}
where $G(t, \boldsymbol{r};s, \bm \rho)$ is the fundamental solution of the
three-dimensional scalar wave equation defined by
\begin{displaymath}
   G(t, \boldsymbol{r};s, \bm \rho) =
 \frac{1}{4\pi|\boldsymbol{r}-\bm \rho|}\cdot\delta\left(\frac{|\boldsymbol{r}-\bm \rho|}{c}
 - (t-s) \right),
\end{displaymath}
$\partial_{\nu(\boldsymbol{r})}$ denotes the outward normal derivative on
$\Gamma$ with
respect to the variable $\boldsymbol{r}$, and $u_{\mathrm N}$ is given by
(\ref{eq:special_sol_point_source}) for moving point sources, and
(\ref{eq:special_sol_dipole_source}) for moving dipole sources.
The observation points of $\phi$ are arranged as
\begin{align}
   \boldsymbol{r}_{j,k} =
 & (x(\theta_j,\varphi_k),\ y(\theta_j,\varphi_k),\  z(\theta_j)) 
                \equiv
(2 \sin\theta_j\cos\varphi_k,\ 
                      2 \sin\theta_j\sin\varphi_k,\  
                      2 \cos\theta_j) \in \Gamma, \nonumber \\
                  & \qquad \qquad j = 1,2,\cdots,J,\ k=1,2,\cdots,K,
   \label{eq:arrangement_of_observation}
\end{align}
where $\cos\theta_j$ is the $j$-th collocation point of Gauss-Legendre
quadrature, and $\varphi_k = {2 \pi (k-1)}/{K}$.
In the experiments, we assign $J=18,\ K=36$, and therefore  648
observation points are arranged on $\Gamma$.
We give observation data at $t = t_\ell = \Delta t \times \ell,\
\ell=0,1,2,\cdots,T/\Delta t$ on $[0,T]$ where $\Delta t = 0.1$.
To simulate practical observation conditions, we add Gaussian noise  to
the numerical solution at each time $t=t_\ell$ such that
\begin{displaymath}
 \sqrt{\frac
 {\displaystyle\int_\Gamma |\phi^{\mathrm{obs}}(t_\ell, \boldsymbol{r})
 - \phi(t_\ell,\boldsymbol{r})|^2 dS(\boldsymbol{r})}
 {\displaystyle\int_\Gamma |\phi(t_\ell,\boldsymbol{r})|^2 dS(\boldsymbol{r})}} =
 0.1\%,\ 0.5\%,\ 1.0\%, 5.0\%
\end{displaymath}
where 
$\phi^{\mathrm{obs}}$ denotes the
observation data perturbed by Gaussian noise.
%
The application of our reconstruction procedure is performed for every $\tau = \tau_\ell = \Delta  \tau \times \ell,\ \ell=0,1,2,\cdots, T_e/\Delta\tau$, where  $\Delta\tau=0.1$, $T_e = 60.0$.
We approximate the surface integral on $\Gamma$ 
using the Gauss-Legendre quadrature with respect to $\theta$ and the
  trapezoidal rule with respect to $\varphi$.
To approximate the derivative with respect to $\tau$, we apply the central difference, e.\@g.\@
${\mathcal{R}}_\phi(g_n)(\tau)$ is approximated by
\begin{equation*}
 {\mathcal{R}}_\phi(g_n)(\tau) 
= {d}_\tau {\mathcal{R}}_\phi(f_n)(\tau) 
\sim \frac{{\mathcal{R}}_\phi(f_n)(\tau+\Delta \tau)
-{\mathcal{R}}_\phi(f_n)(\tau-\Delta \tau) }{2\cdot\Delta\tau}.
\end{equation*}
\par
%
%
\subsection{Reconstruction of moving point sources}
We show numerical experiments for the reconstruction of moving
point sources.
We arrange the parameters of three point sources as follows:
\begin{description}
\item[source 1.] 
${\boldsymbol{p}}_1(t) = (0.8,\ -0.3,\ 0.8\cos(0.4t)-0.2)$,
\begin{displaymath}
  q_1(t) = \left\{ \begin{array}{ll}
                   0, & 0 \leq t < 4, \\[1ex]
                   \displaystyle \eta\left(\frac{t-4}{10}\right) \cdot \left(1 + \frac{7}{10}\sin\left(\frac{2\pi t}{7}\right)\right),
                   & 4 \leq t < 14, \\[2ex] 
                   \displaystyle 1 + \frac{7}{10}\sin\left(\frac{2\pi t}{7}\right),
                   & 14 \leq t < 35, \\[2ex]
                   \displaystyle \left(1 -  \eta\left(\frac{t-35}{20}\right)\right)  \cdot 
                      \left(1 + \frac{7}{10}\sin\left(\frac{2\pi t}{7}\right)\right),
		   & 35 \leq t < 55, \\[2ex] 
                   0, & 55 \leq t \leq 70.
                   \end{array}
                   \right.
\end{displaymath}
\item[source 2.] ${\boldsymbol{p}}_2(t) = (p_{2,x}(t),\ p_{2,y}(t),\ p_{2,z}(t))$ where
\begin{align*}
      p_{2,x}(t) & =  \sin(-0.2(t+2.5)) + 0.3, \\
      p_{2,y}(t) & =  0.7\sin(0.4(t+2.5))-0.2, \\
      p_{2,z}(t) & =  0.5\sin(0.5(t+2.5)) + 0.2,
\end{align*}
\begin{displaymath}
  q_2(t) = \left\{ \begin{array}{ll}
                   0, & 0 \leq t < 10, \\[1ex]
                   \displaystyle \eta\left(\frac{t-10}{12}\right) \cdot \sin\left(\frac{2\pi (t-4)}{12}\right), 
                   & 10 \leq t < 22, \\[2ex] 
                   \displaystyle \sin\left(\frac{2\pi (t-4)}{12}\right),
                   & 22 \leq t < 50, \\[2ex] 
    	             \displaystyle \left(1- \eta\left(\frac{t-50}{20}\right)\right)  \cdot   \sin\left(\frac{2\pi (t-4)}{12}\right),
		   & 50 \leq t < 70.
                   \end{array}
                   \right.
\end{displaymath}
\item[source 3.] ${\boldsymbol{p}}_3(t) = (p_{3,x}(t),\ p_{3,y}(t),\ p_{3,z}(t))$ where
\begin{align*}
       p_{3,x}(t) & = 0.9\cos(0.2\pi)\cos(-0.3(t+2.7)) - 0.15\sin(0.2\pi)\sin(-0.3(t+2.7)) -0.5, \\
       p_{3,y}(t) & = 0.9\sin(0.2\pi)\cos(-0.3(t+2.7)) + 0.15\cos(0.2\pi)\sin(-0.3(t+2.7)) +0.6, \\
       p_{3,z}(t) & = 0.8\sin(-0.25(t+2.7)),
\end{align*}
\begin{displaymath}
  q_3(t) = \left\{ \begin{array}{ll}
                   0, & 0 \leq t < 25, \\
                   \displaystyle -\frac{3}{2}\,\eta\left(\frac{t-25}{7}\right),
                   & 25 \leq t < 32, \\
                   -1.5, 
                   & 32 \leq t < 34, \\
                   \displaystyle -\frac{3}{2} \left(1 - \eta\left(\frac{t-34}{10}\right)\right),
		   & 34 \leq t < 44, \\
                   0, & 44 \leq t \leq 70.
                   \end{array}
                   \right.
\end{displaymath}
\end{description}
Here, we set
\begin{displaymath}
    \eta(s) = \left\{\begin{array}{ll}
            0, & s < 0, \\[1ex]
            \displaystyle s - \frac{6\sin(2 \pi s) + \sin^3(2 \pi
            s)}{12\pi},
            & 0 \leq s < 1, \\[1ex]
            1, & 1 \leq s.
     \end{array}\right.
\end{displaymath}
Therefore, $K(\tau)$ changes as
\begin{displaymath}
     K(\tau) = \left\{\begin{array}{ll}
           0, & 0 \leq \tau <3.9,\ \  69.5 \leq \tau \\[1ex]
           1, & 3.9 \leq \tau < 10.2,\ \ 54.0 \leq \tau < 69.5, \\[1ex]
           2, & 10.2 \leq \tau < 24.6,\ \ 44.6 \leq \tau < 54.0, \\[1ex]
           3, & 24.6 \leq \tau < 44.6.
\end{array}\right.
\end{displaymath}
We note that $|\boldsymbol{p}_k(t)| < 1.6$ for any $k$ and $t$, and the range of the moving speeds of sources are $0 \leq
|\dot{\boldsymbol{p}}_1(t)| \leq 0.32$, $0.165 \leq |\dot{\boldsymbol{p}}_2(t)| \leq 0.425$
and $0.045 \leq |\dot{\boldsymbol{p}}_3(t)| \leq 0.333$.
We also note that $q_k \in C^4([0,70]),\ k=1,2,3$.
\par
\begin{figure}[ht]
\centering
\subfloat[Source 1]
{
 \resizebox*{9cm}{!}{\includegraphics{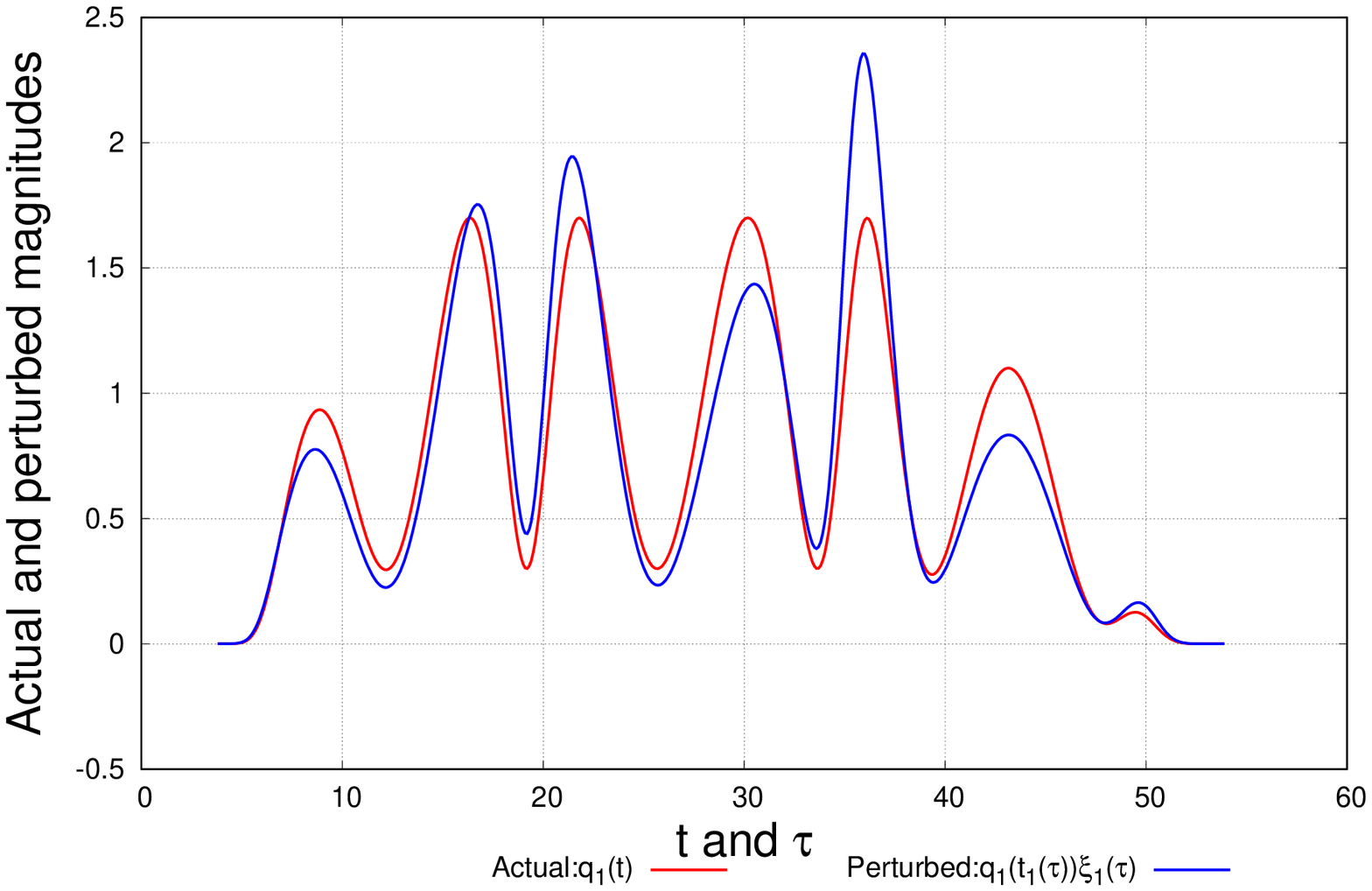}}
 } \\[-2ex] 
\subfloat[Source 2]
{
 \resizebox*{9cm}{!}{\includegraphics{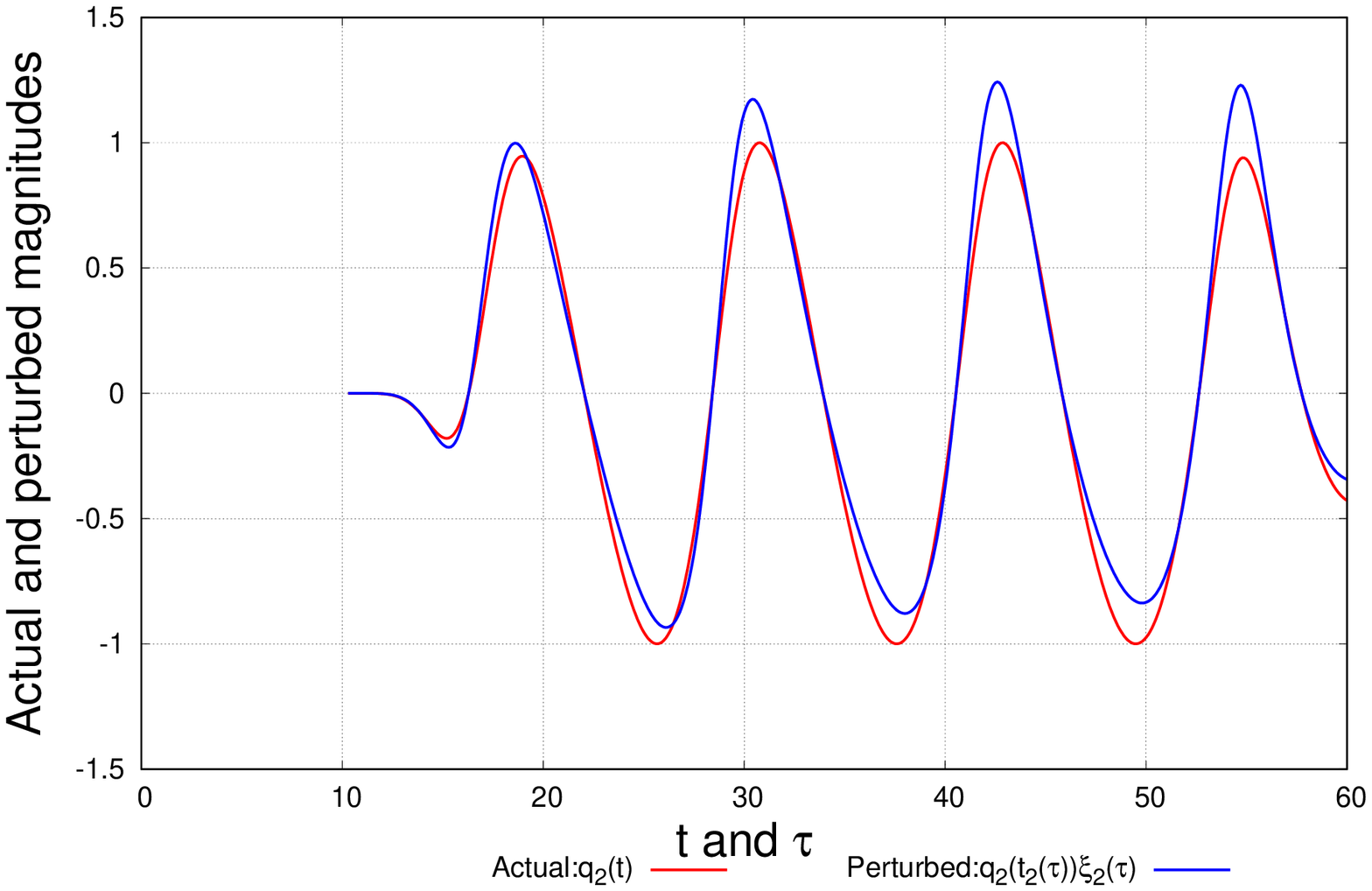}}
 } \\[-2ex] 
\subfloat[Source 3]
{
 \resizebox*{9cm}{!}{\includegraphics{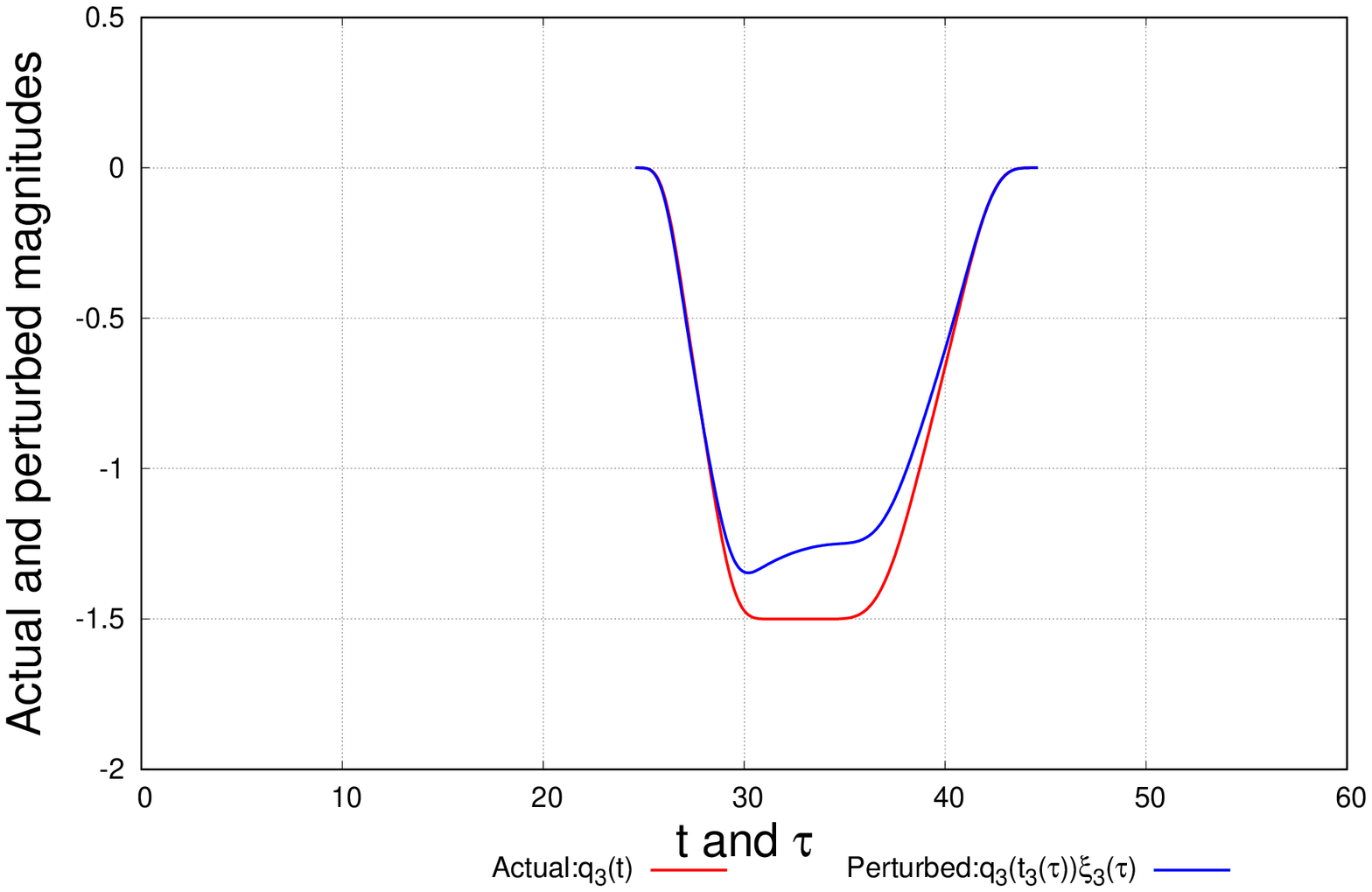}}
 }  
  \caption{Behaviour of the actual and perturbed magnitude of point sources}
  \label{fig:behaviour_actual_perturbed_point_sources}
\end{figure}

Before we show the reconstruction results of moving point source, we
compare the behaviours of magnitudes $q_k(t)$ with respect to $t$
and perturbed magnitudes $q_k(t_k(\tau))\xi_k(t_k(\tau))$ with respect
to $\tau$ in Figure \ref{fig:behaviour_actual_perturbed_point_sources}.
We can find large perturbation on magnitudes due to the $z$-component
of the moving velocities of sources, e.\@g.\@ around $\tau=21,\ 30,\ 35$, and $42$
for source 1 and $30 \leq \tau \leq 37$ for source 3.
One of the purpose of this paper is to give a method to correct these perturbations.
\par
\begin{figure}[ht]
\centering
\subfloat[behaviour of $|\det H_{k,0}(\tau)|$ for each $k$]
{
 \resizebox*{9cm}{!}{\includegraphics{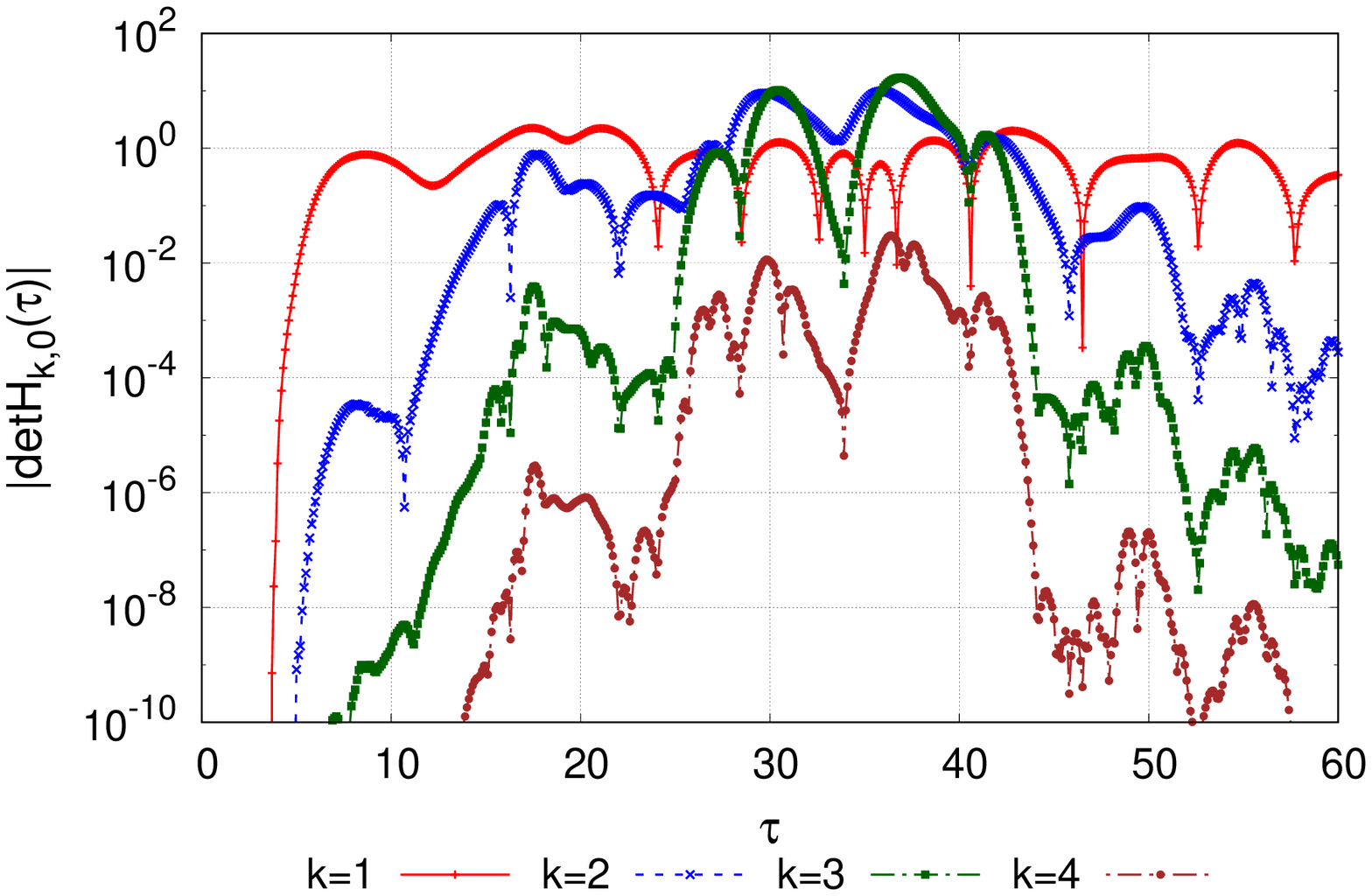}} } \\ 
\subfloat[behaviour of $|\det H_{k,0}(\tau)|$ for fixed $\tau$]
{ \resizebox*{9cm}{!}{\includegraphics{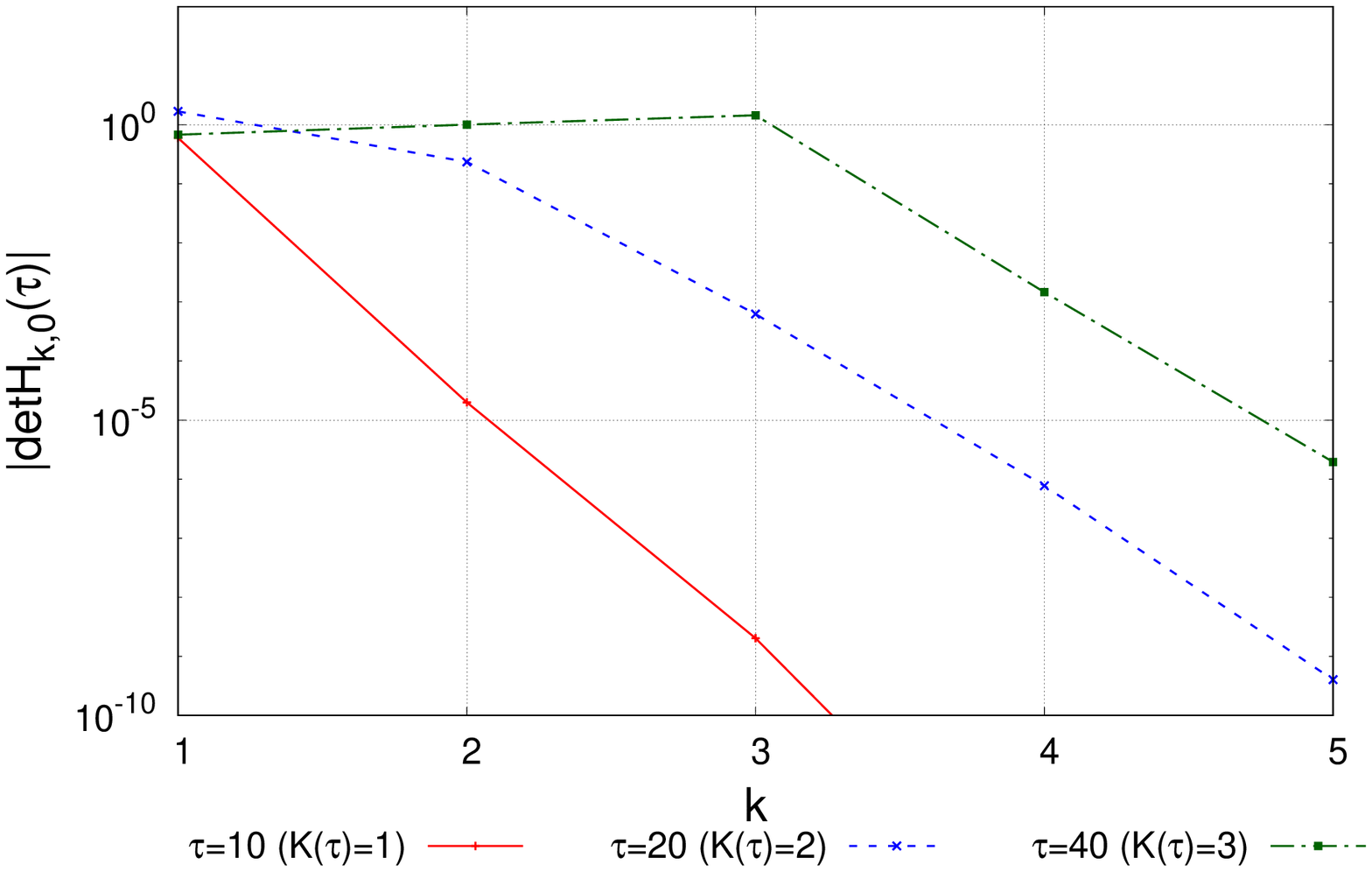}}}
  \caption{Behaviour of $|\det H_{k,0}(\tau)|$ for observation data without noise: point source case.}
  \label{fig:behaviour_detH_point_source}
\end{figure}

Firstly, we discuss the identification of the number $K(\tau)$
of point sources.
Figure \ref{fig:behaviour_detH_point_source}(a) shows the behaviour of $|\det H_{k,0}(\tau)|$ in $0 \leq \tau
\leq T_e$ for $k=1,2,3,4$, and Figure \ref{fig:behaviour_detH_point_source}(b) gives the behaviour of $|\det
H_{k,0}(\tau)|$ at $\tau =10,\ 20$ and $40$ for observation data without noise.
As we mentioned in Note 3.3, the results of Figure \ref{fig:behaviour_detH_point_source} show that we
can not apply the condition (\ref{eq:determine_K}) as it is for the identification of $K(\tau)$.
However,  we can observe that $|\det H_{k,0}(\tau)|$ has relatively large gap
between $k=K(\tau)$ and $K(\tau)+1$. 
From this observation, we apply the following algorithm to identify $K(\tau)$:
\vspace{2ex}
\par
\noindent
\textbf{Algorithm 4.1} 
\begin{description}
\item[Step 0.] Set positive parameters $\varepsilon_0$ and
	   $\varepsilon_G$.
\item[Step 1.] If $|\det H_{1,0}(\tau)| < \varepsilon_0$ and $|\det H_{1,0}(\tau)| > |\det H_{2,0}(\tau)|$, then $K(\tau) = 0$.
\item[Step 2.] If $|\det H_{1,0}(\tau)| > \varepsilon_0$ or $|\det H_{1,0}(\tau)| < |\det H_{2,0}(\tau)|$, then find largest number
	   $k$ such that $|\det H_{k,0}(\tau)|/|\det H_{k-1,0}(\tau)| > \varepsilon_G$, and identify $K(\tau) = k$.
\end{description}
\vspace{2ex}
\par
\begin{figure}[ht]
\centering
   \subfloat[]{
   \resizebox*{9cm}{!}{\includegraphics{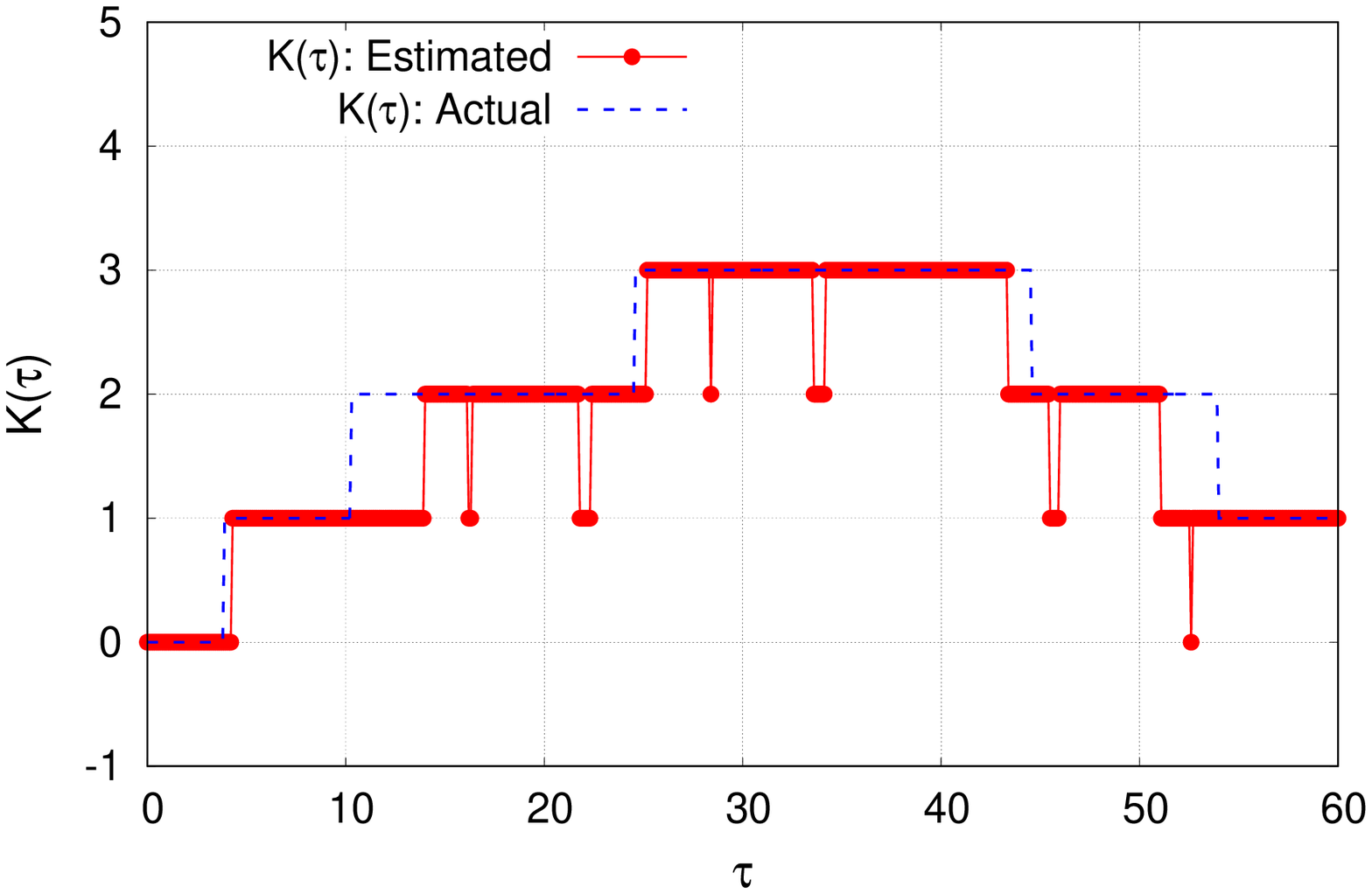}}} \\ 
   \subfloat[]{
   \resizebox*{9cm}{!}{\includegraphics[width=7.5cm]{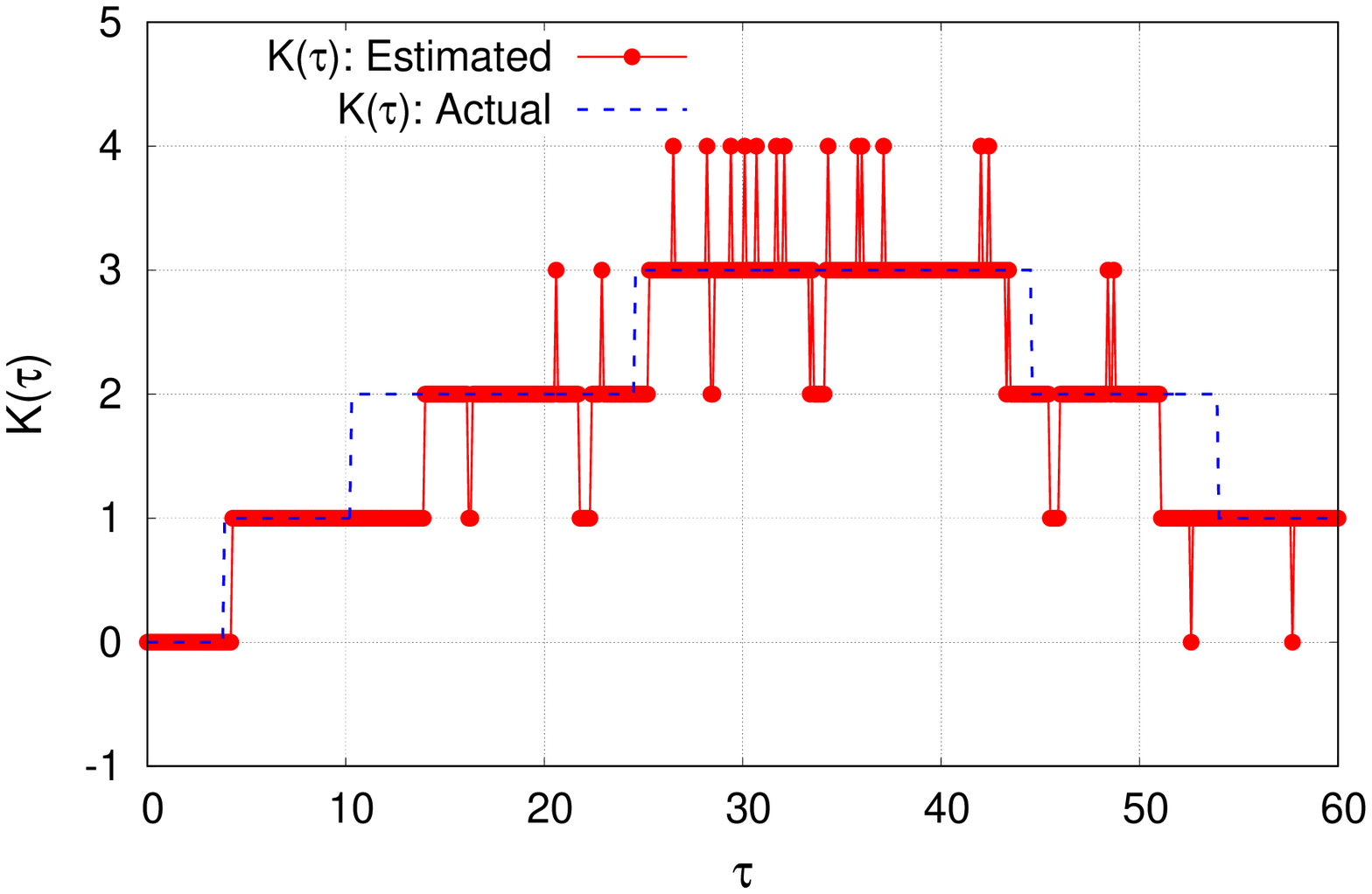}}}
  \caption{Behaviour of estimated number of point sources: (a) for
 observation data without noise, (b) with 0.5\% noise.}
  \label{fig:behaviour_K_point_source}
\end{figure}

\noindent
Figure \ref{fig:behaviour_K_point_source} displays the identification results of $K(\tau)$ from observation data
without noise and with 0.5\% noise.
Here, we set $\varepsilon_0 = 1.0\times 10^{-4}$ and $\varepsilon_G =
2.5 \times 10^{-2}$.
From the result of Figure \ref{fig:behaviour_K_point_source}, we
consider that Algorithm 4.1 works well even for observations with
0.5\% noise.
One can observe that identified $K(\tau)$ is smaller than the actual value around
$\tau=12$ and $\tau=55$.
We discuss a reason of these bad estimates later. 
\par
\par
\begin{figure}[ht]
\centering
  \subfloat[location of source 1]{
  \resizebox*{8cm}{!}{\includegraphics{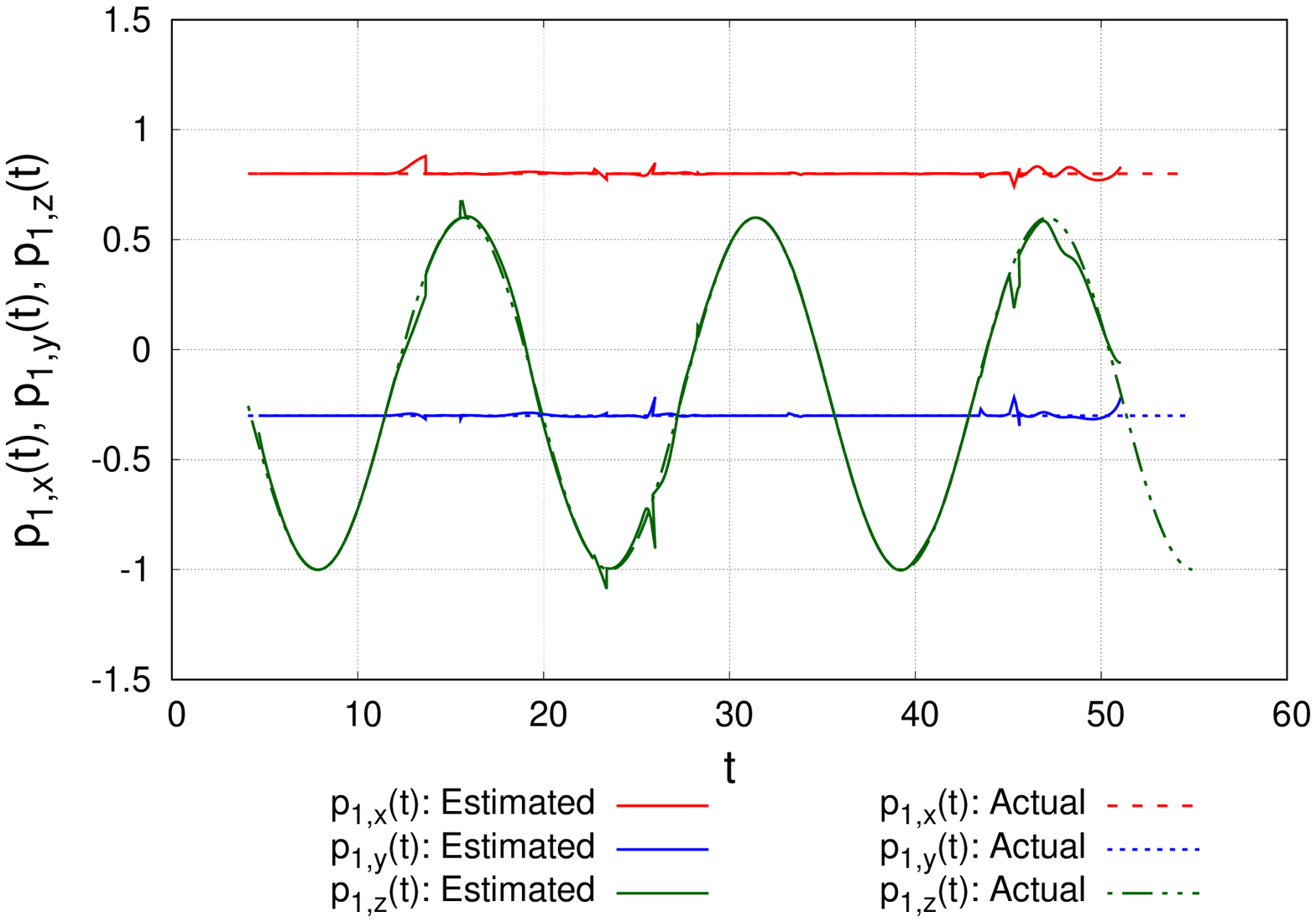}}}
  \subfloat[magnitude of source 1]{
  \resizebox*{8cm}{!}{\includegraphics{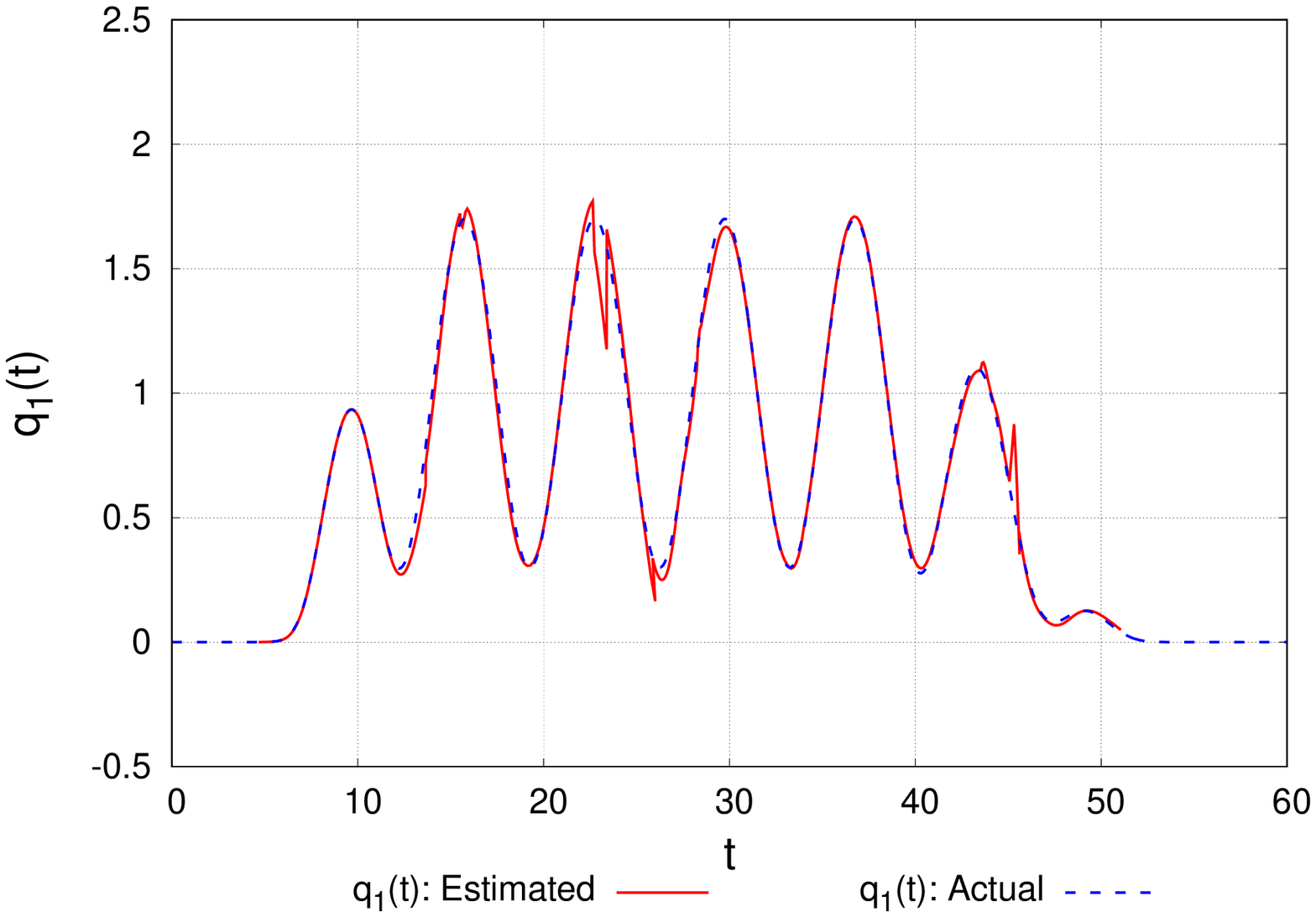}}}
  \\ [-2ex]
  \subfloat[location of source 2]{
  \resizebox*{8cm}{!}{\includegraphics{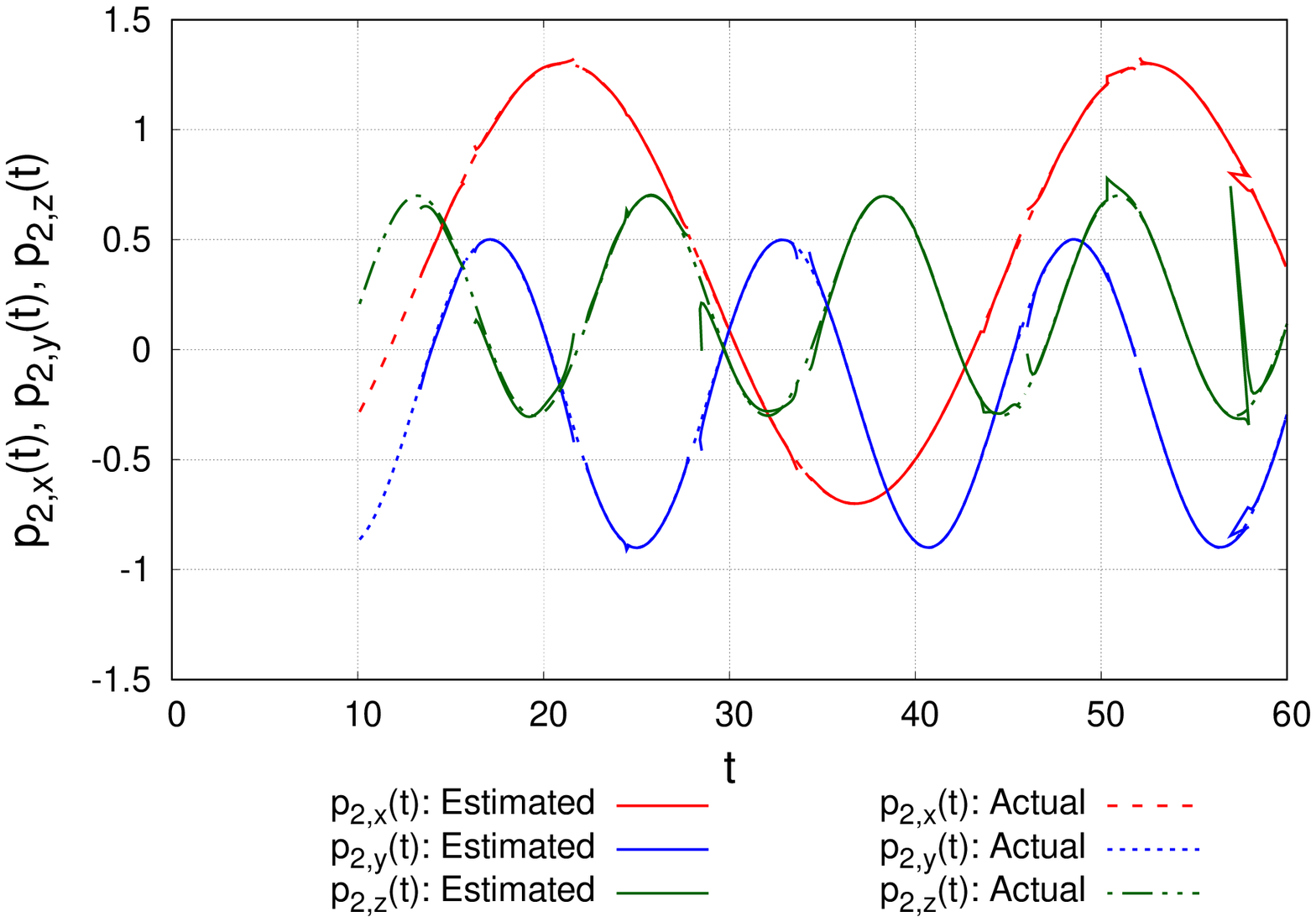}}}
  \subfloat[magnitude of source 2]{
  \resizebox*{8cm}{!}{\includegraphics{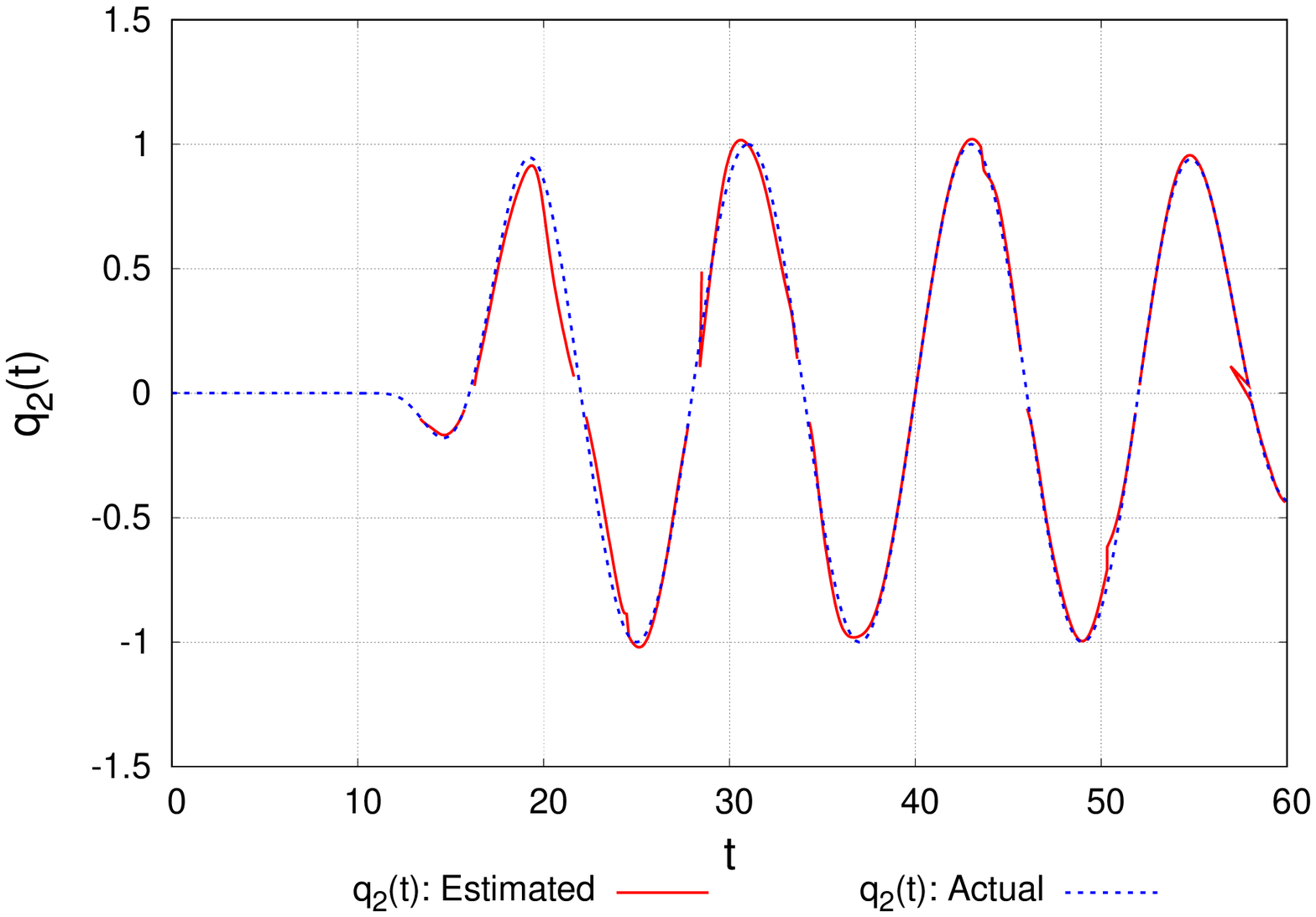}}}
  \\[-2ex]
  \subfloat[location of source 3]{
  \resizebox*{8cm}{!}{\includegraphics{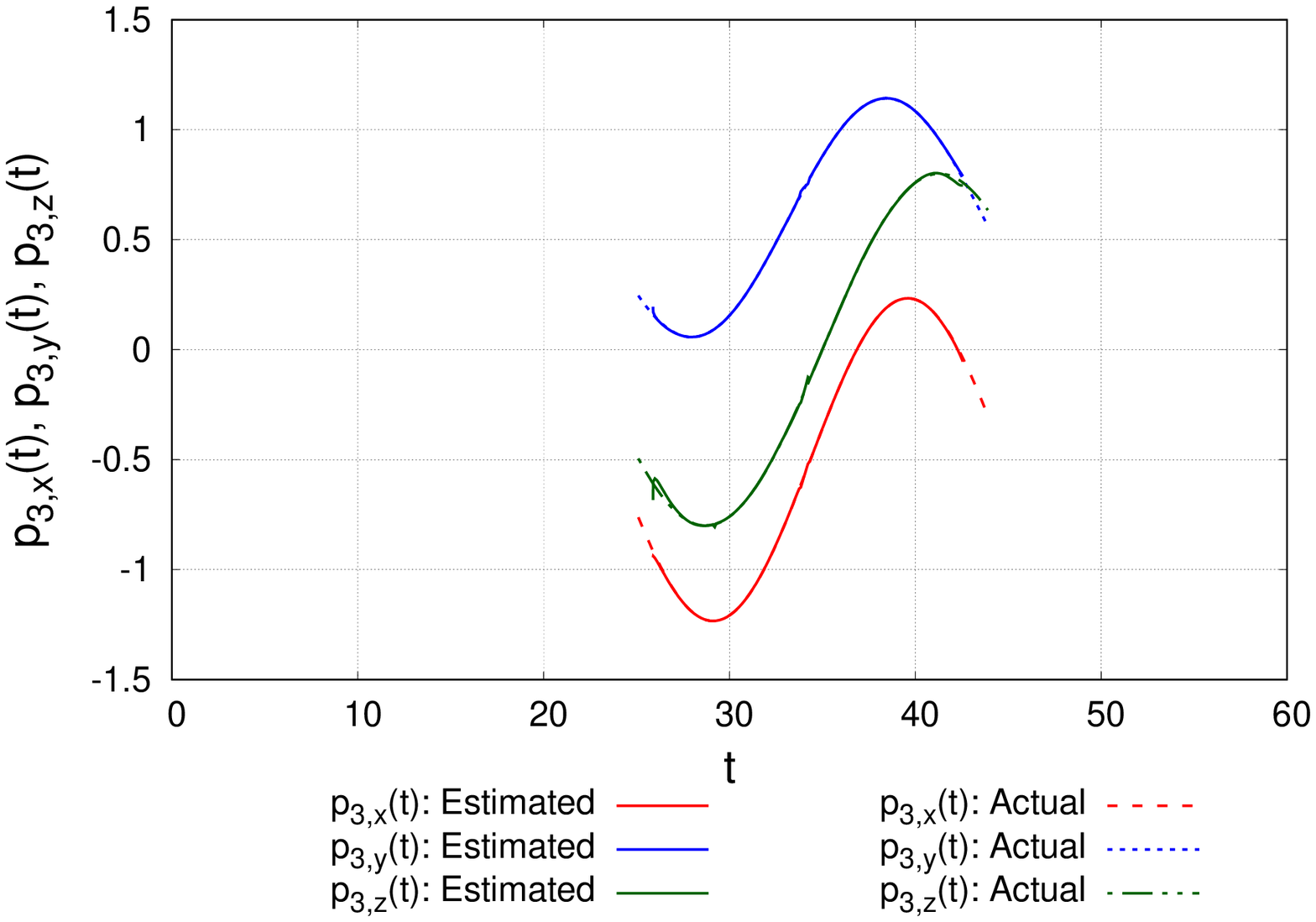}}}
  \subfloat[magnitude of source 3]{
  \resizebox*{8cm}{!}{\includegraphics{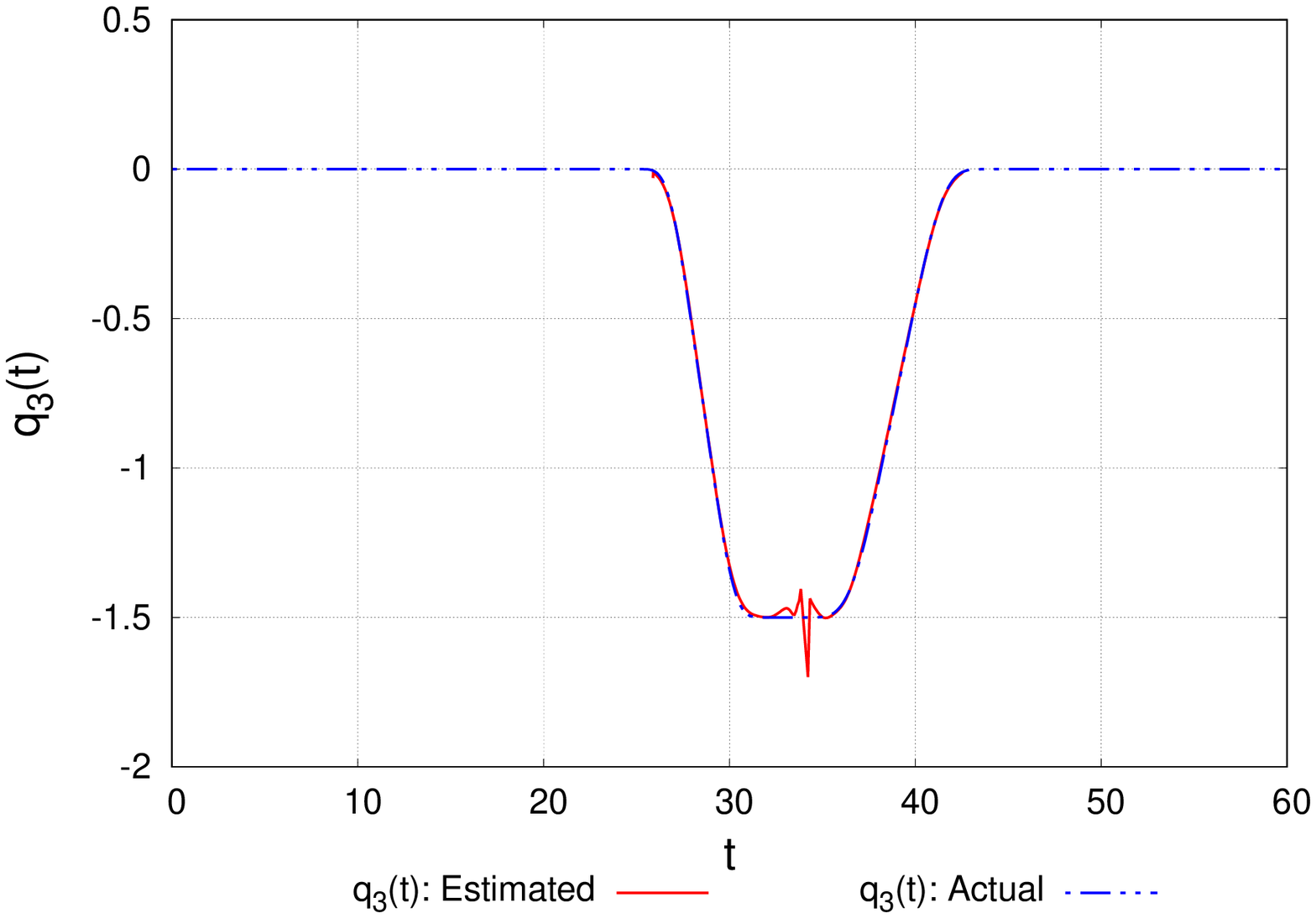}}}
 \caption{Estimated locations and magnitudes of point sources for  observation data without noise.}
\label{fig:reconstruction_results_loc_nl_0.0_point_source}
\end{figure}
\begin{figure}[ht]
\centering
  \subfloat[location of source 1]{
  \resizebox*{8cm}{!}{\includegraphics{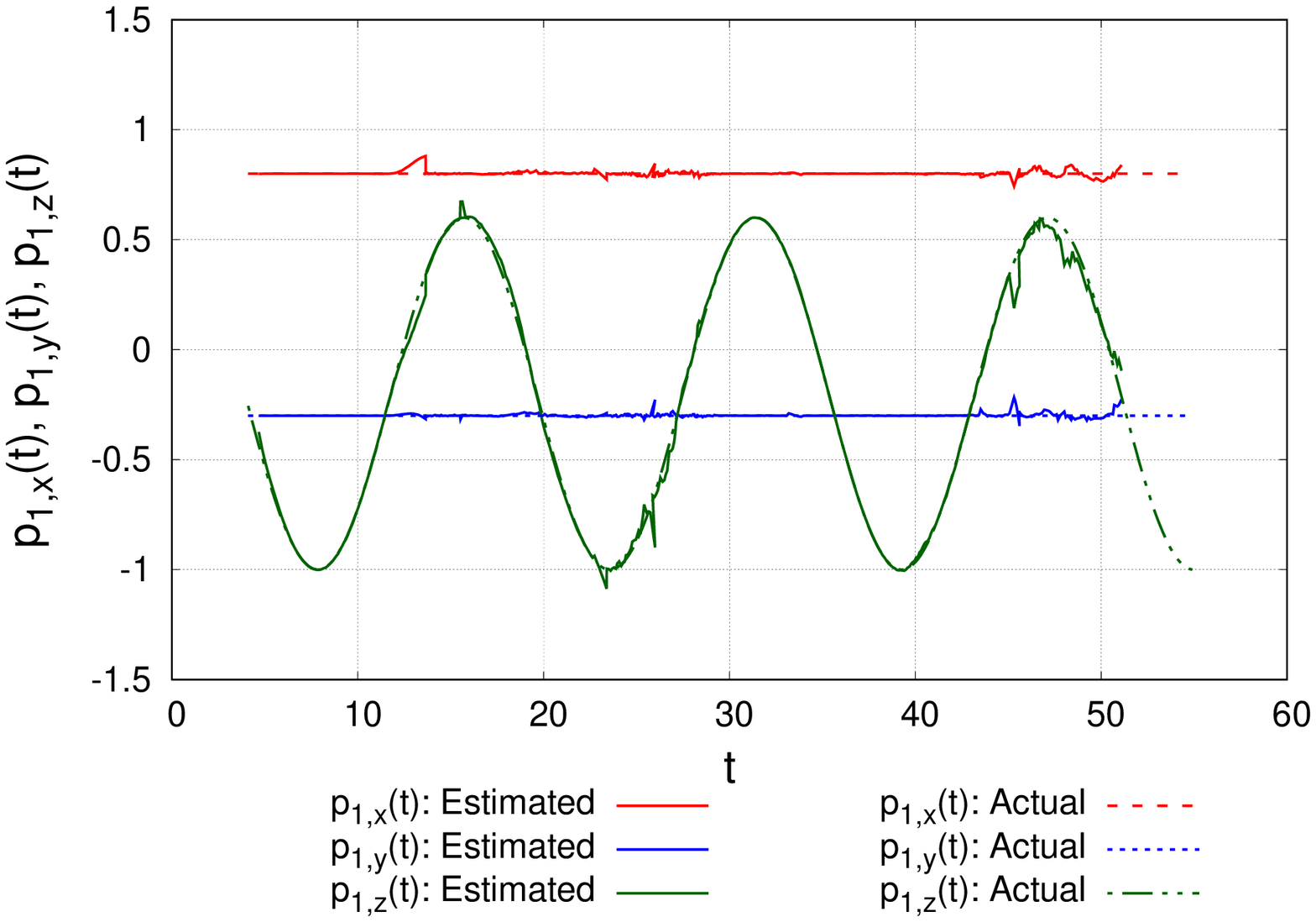}}}
  \subfloat[magnitude of source 1]{
  \resizebox*{8cm}{!}{\includegraphics{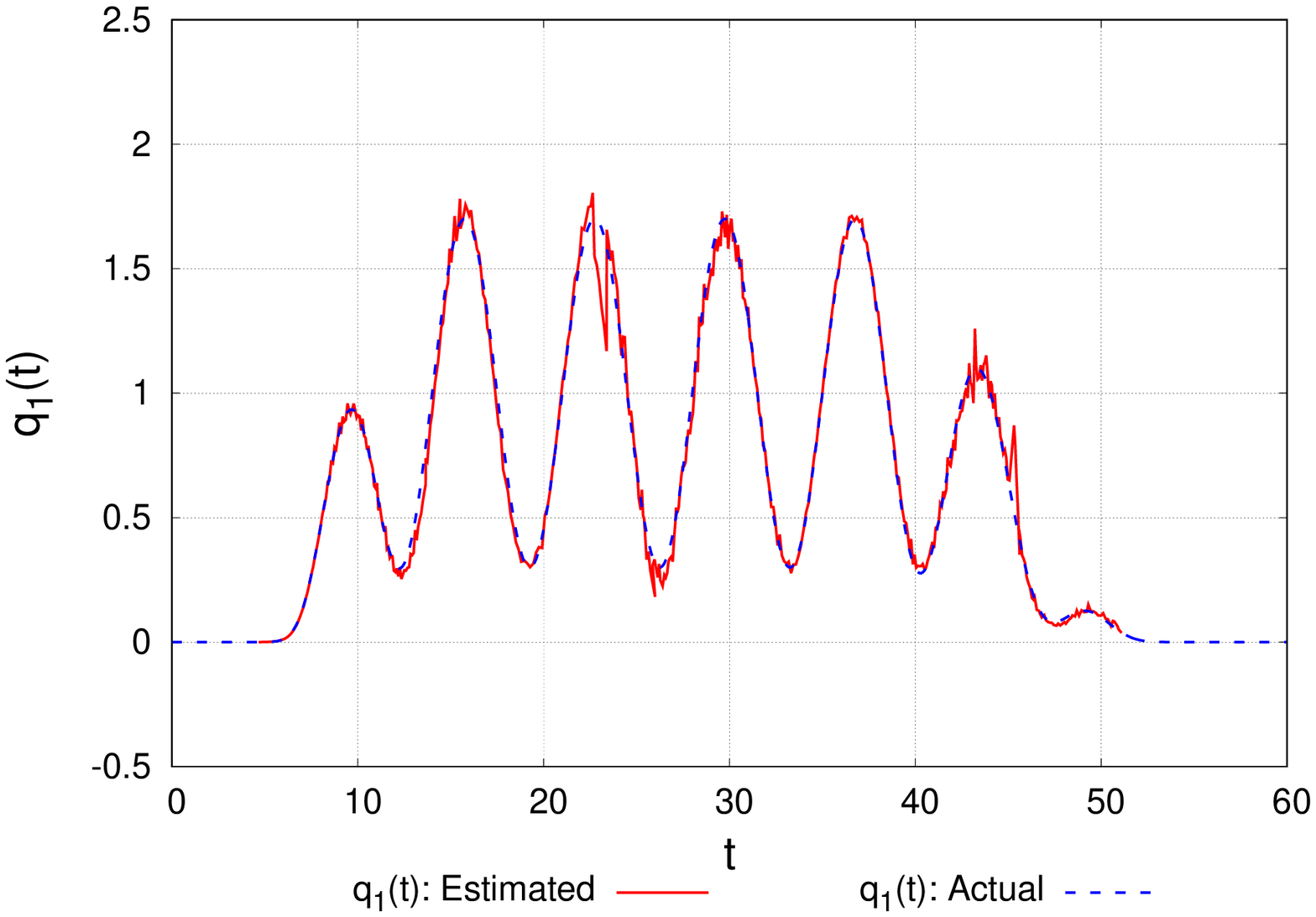}}}
 \\[-2ex]
  \subfloat[location of source 2]{
  \resizebox*{8cm}{!}{\includegraphics{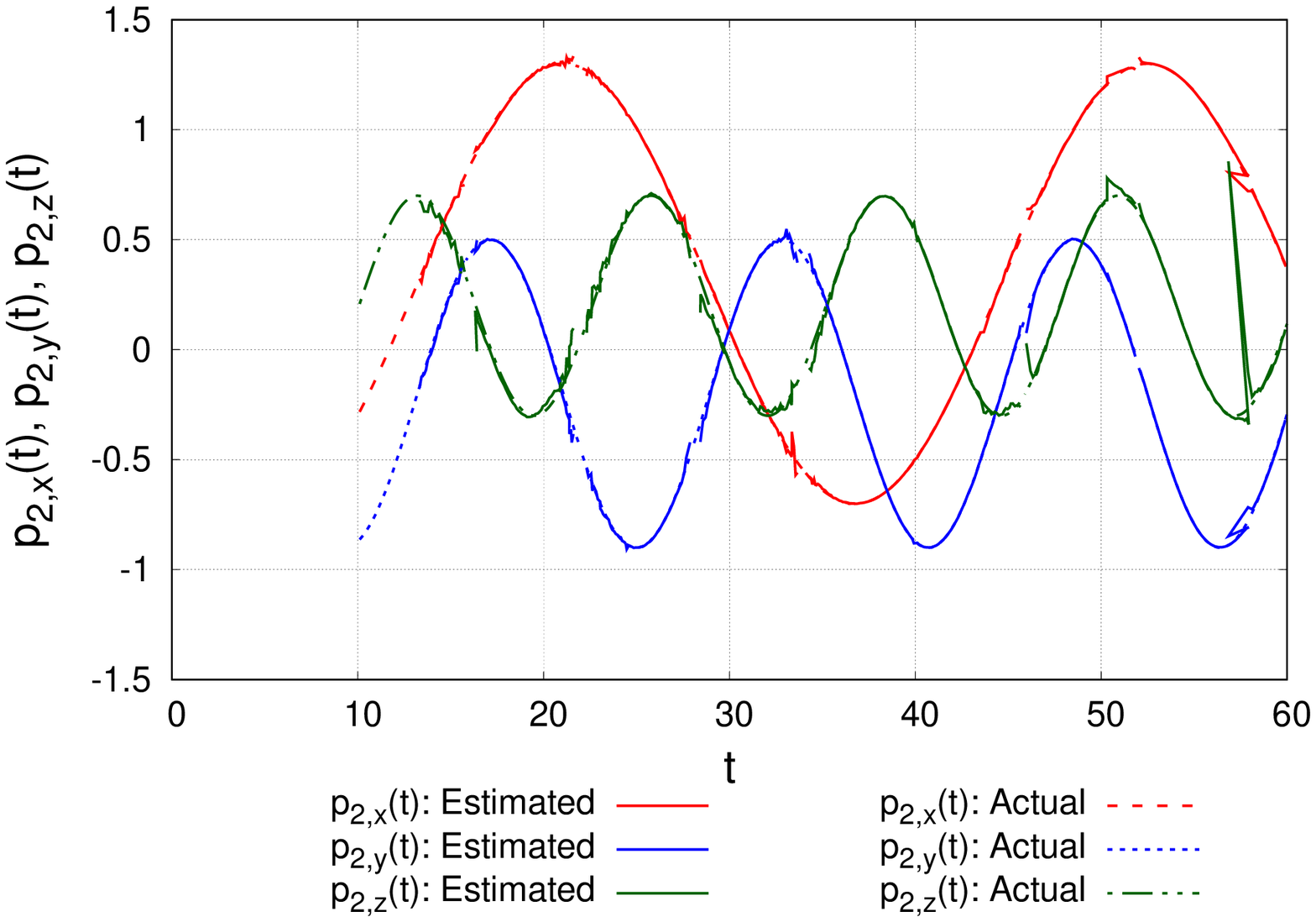}}}
  \subfloat[magnitude of source 2]{
  \resizebox*{8cm}{!}{\includegraphics{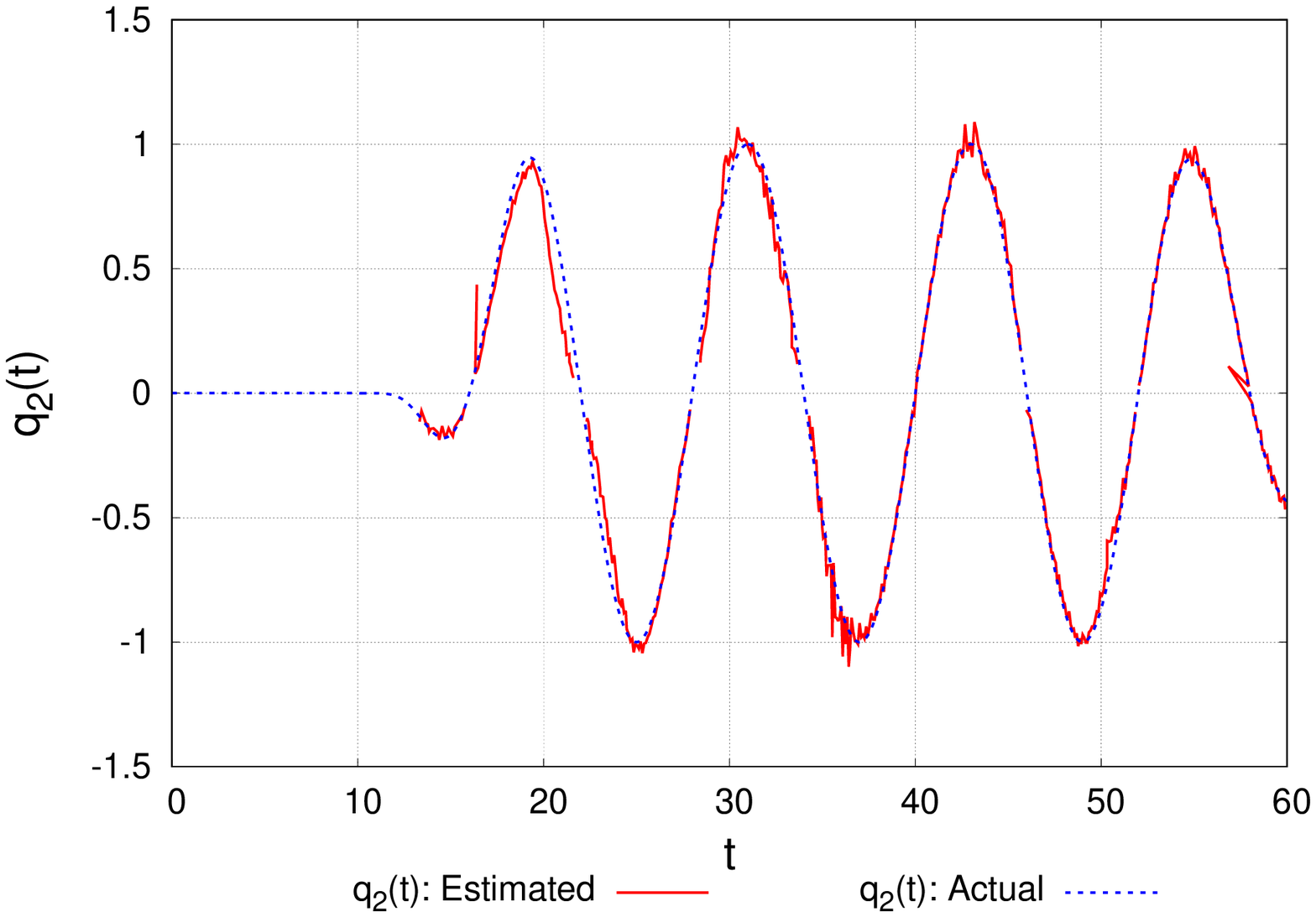}}}
 \\[-2ex]
\subfloat[location of source 3]{
  \resizebox*{8cm}{!}{\includegraphics{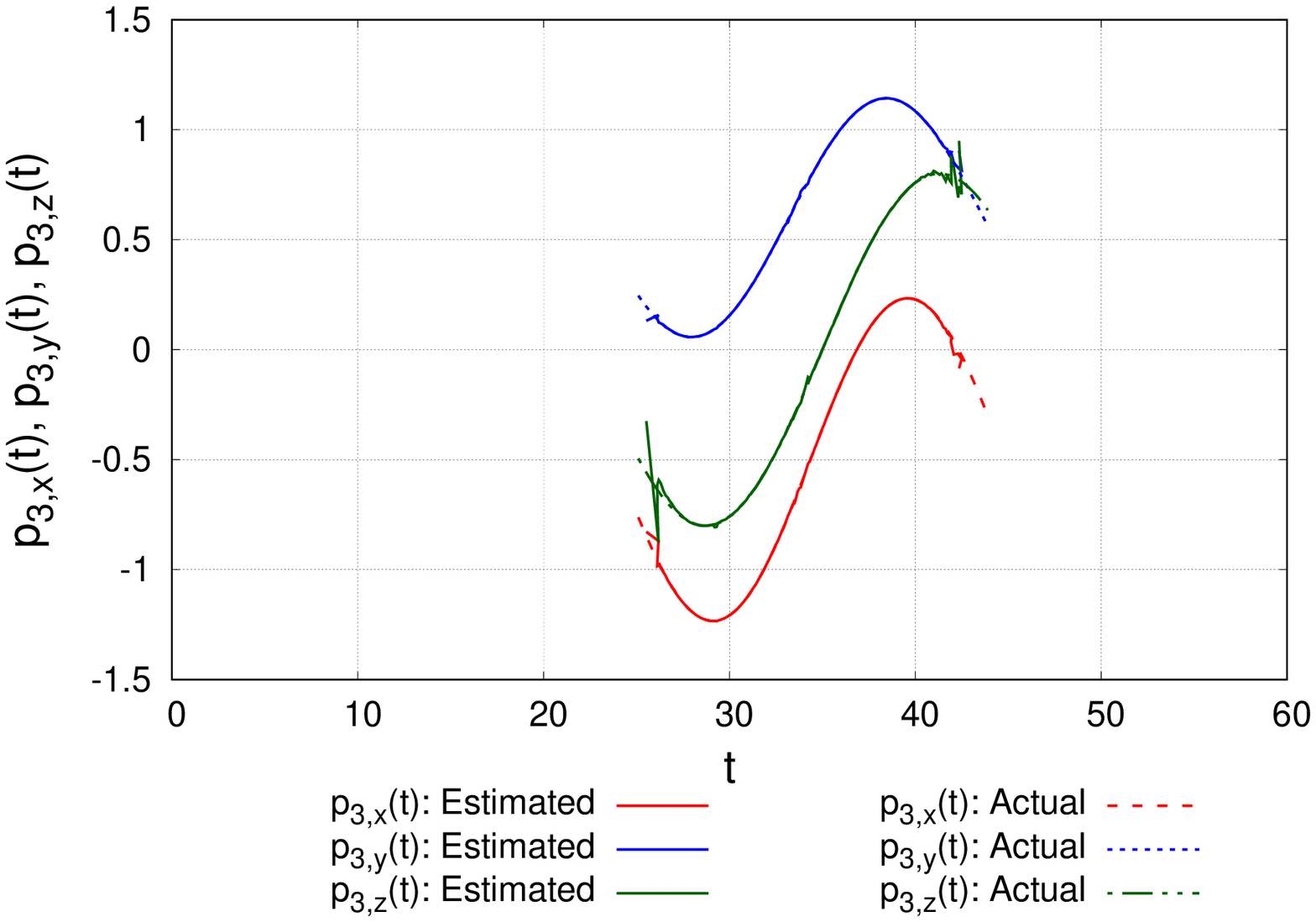}}}
  \subfloat[magnitude of source 3]{
  \resizebox*{8cm}{!}{\includegraphics{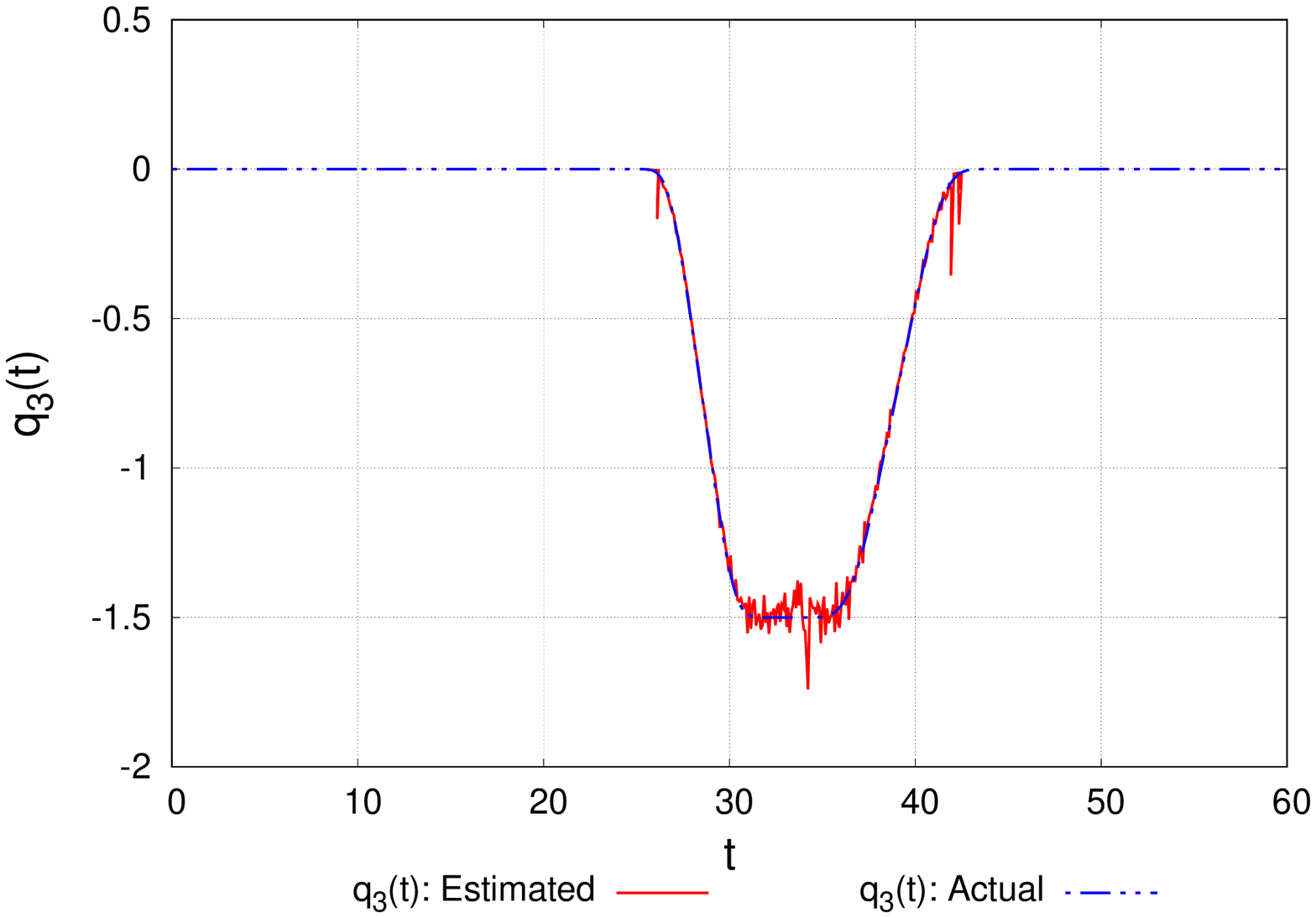}}}
\caption{Estimated locations and magnitudes of point sources for observation data with 0.1\% noise.}
\label{fig:reconstruction_results_loc_nl_0.1_point_source}
\end{figure}
\begin{figure}[ht]
\centering
  \subfloat[location of source 1]{
  \resizebox*{8cm}{!}{\includegraphics{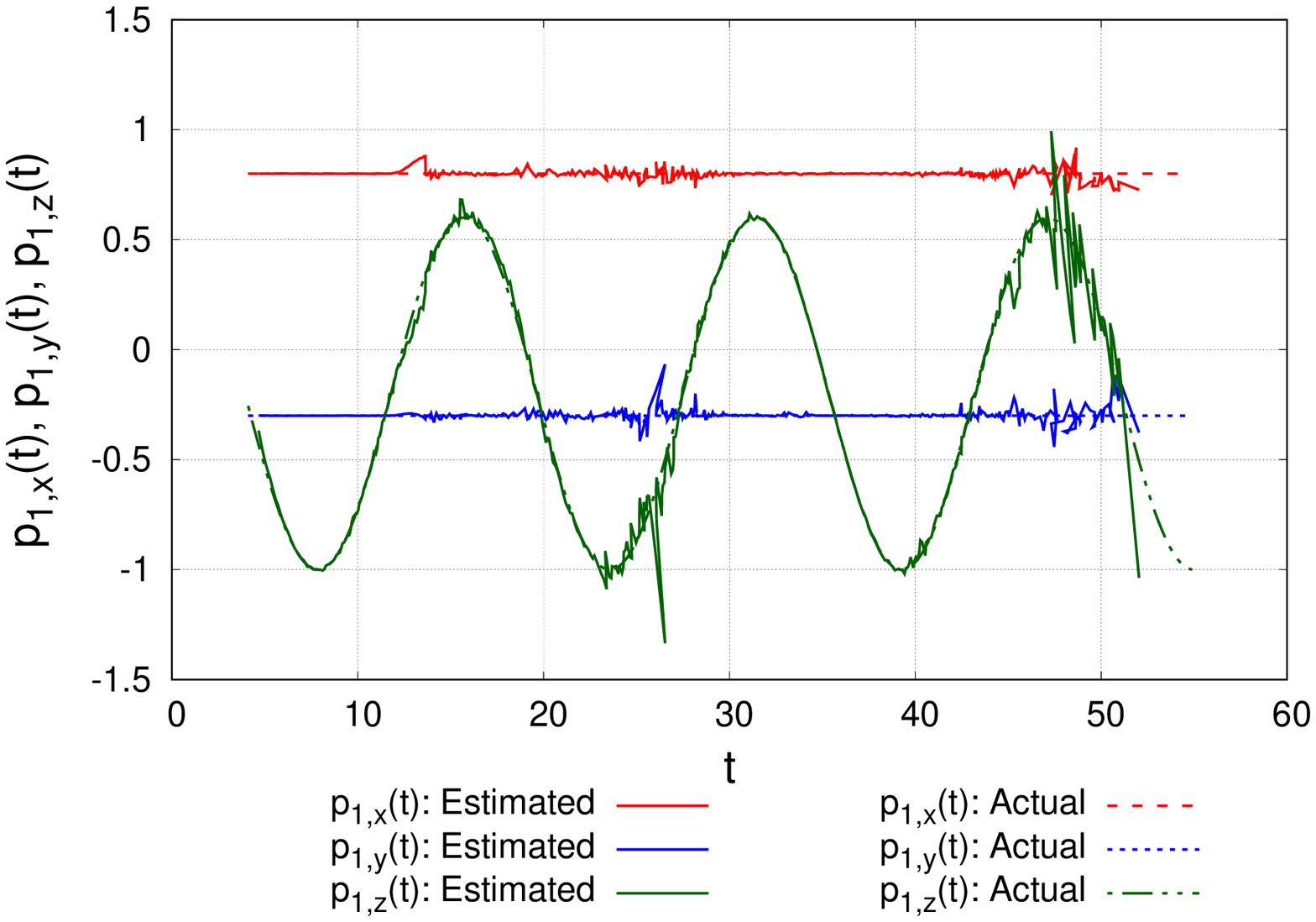}}}
  \subfloat[magnitude of source 1]{
  \resizebox*{8cm}{!}{\includegraphics{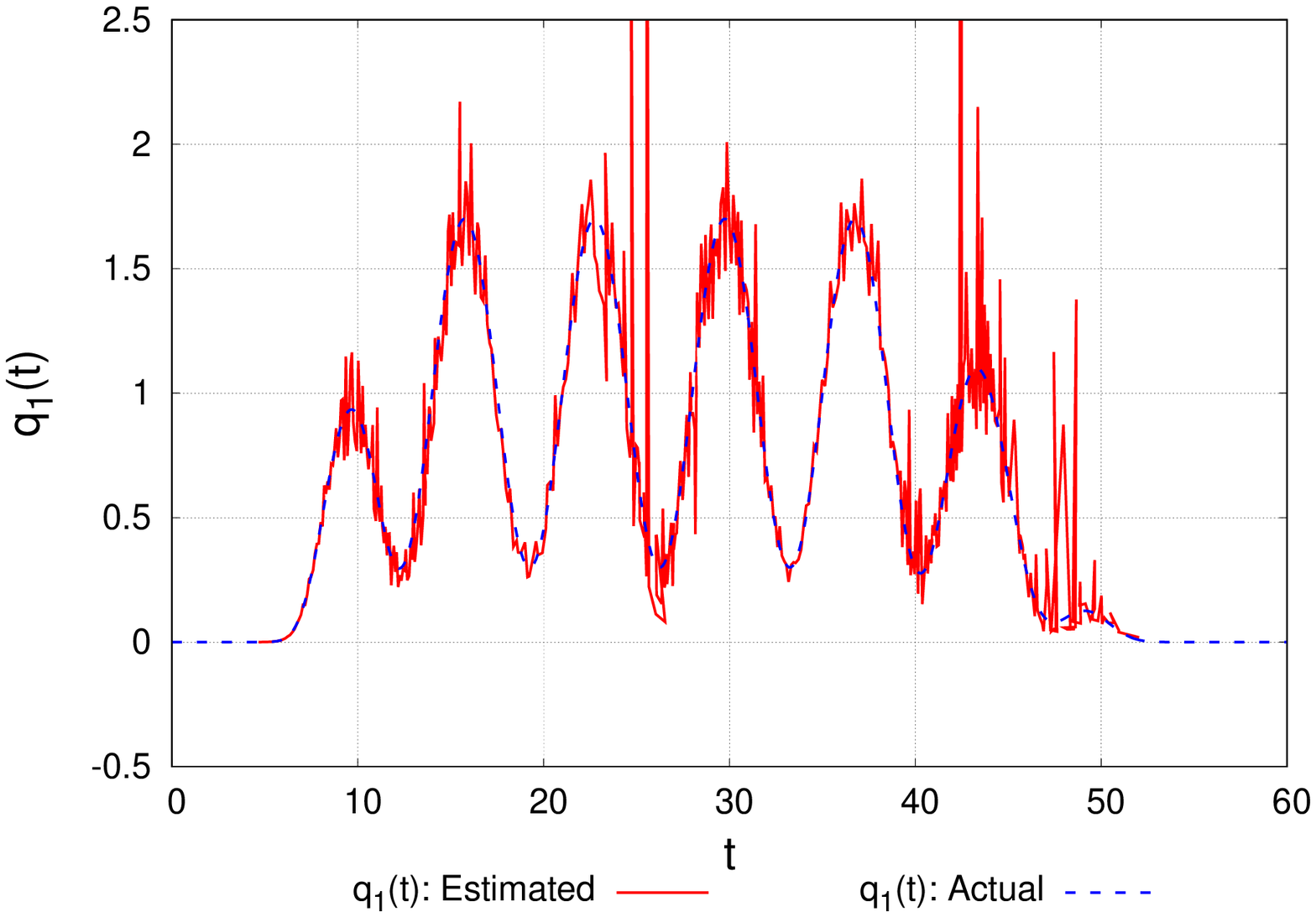}}}
  \\[-2ex]
  \subfloat[location of source 2]{
  \resizebox*{8cm}{!}{\includegraphics{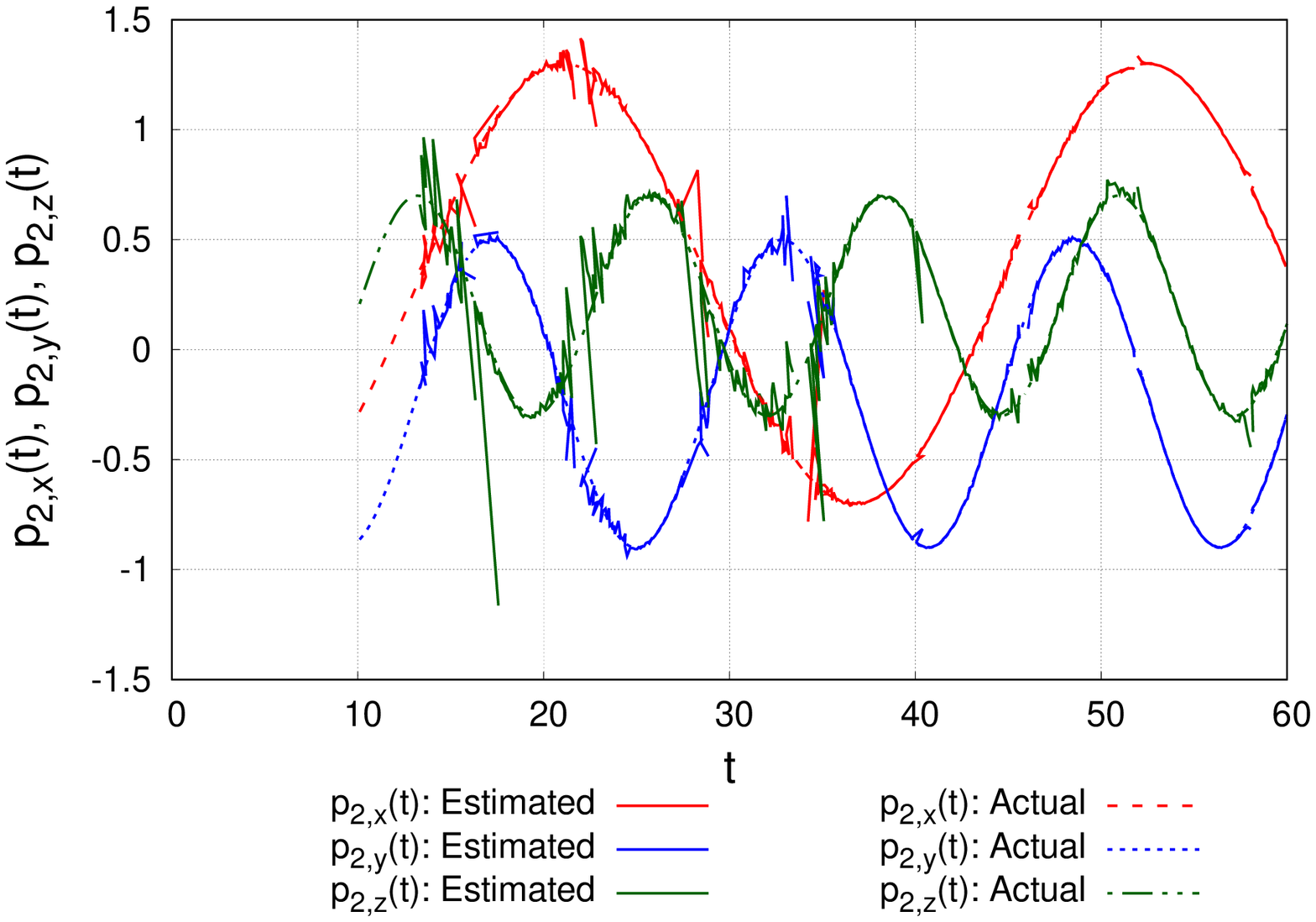}}}
  \subfloat[magnitude of source 2]{
  \resizebox*{8cm}{!}{\includegraphics{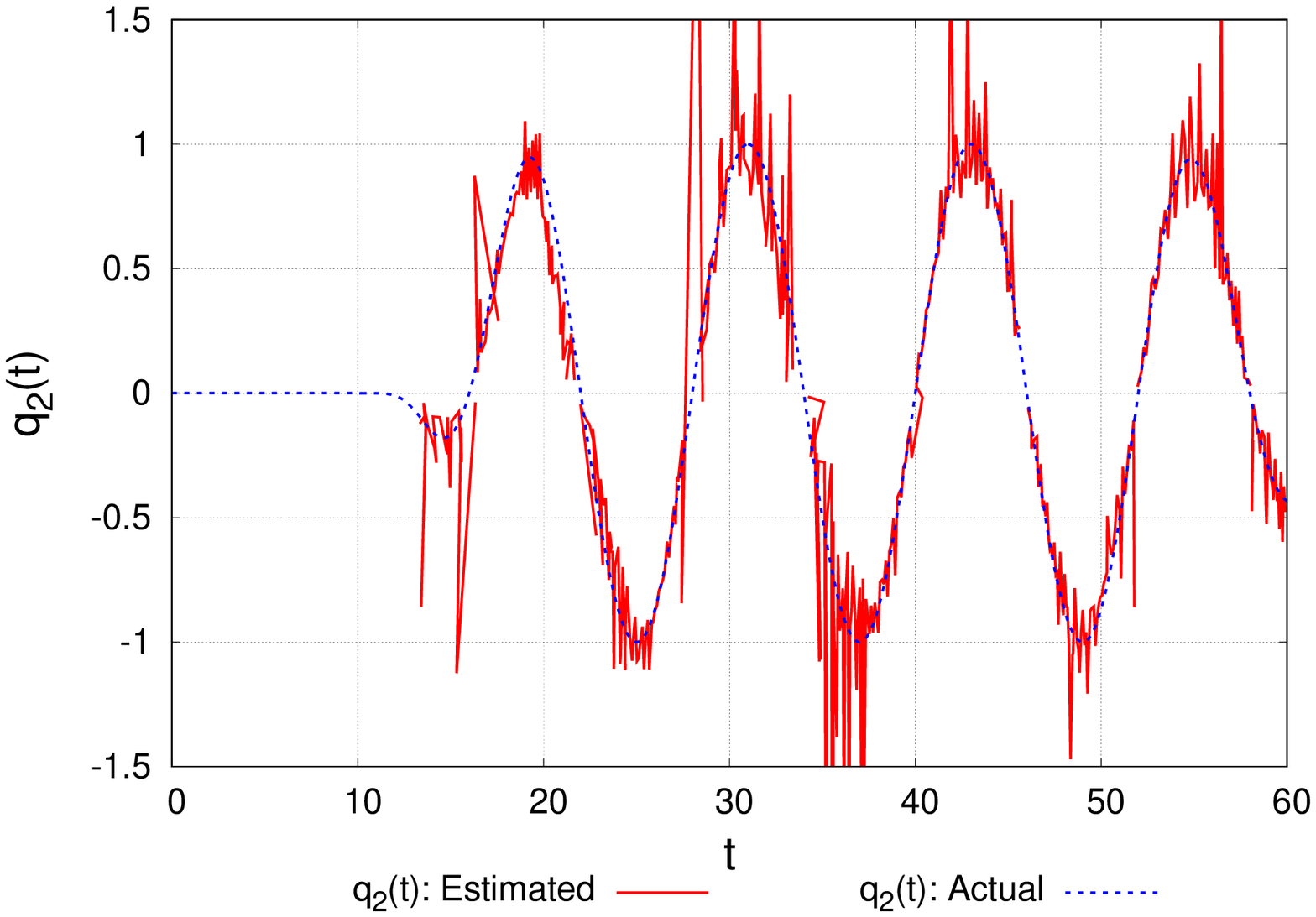}}}
 \\[-2ex]
  \subfloat[location of source 3]{
  \resizebox*{8cm}{!}{\includegraphics{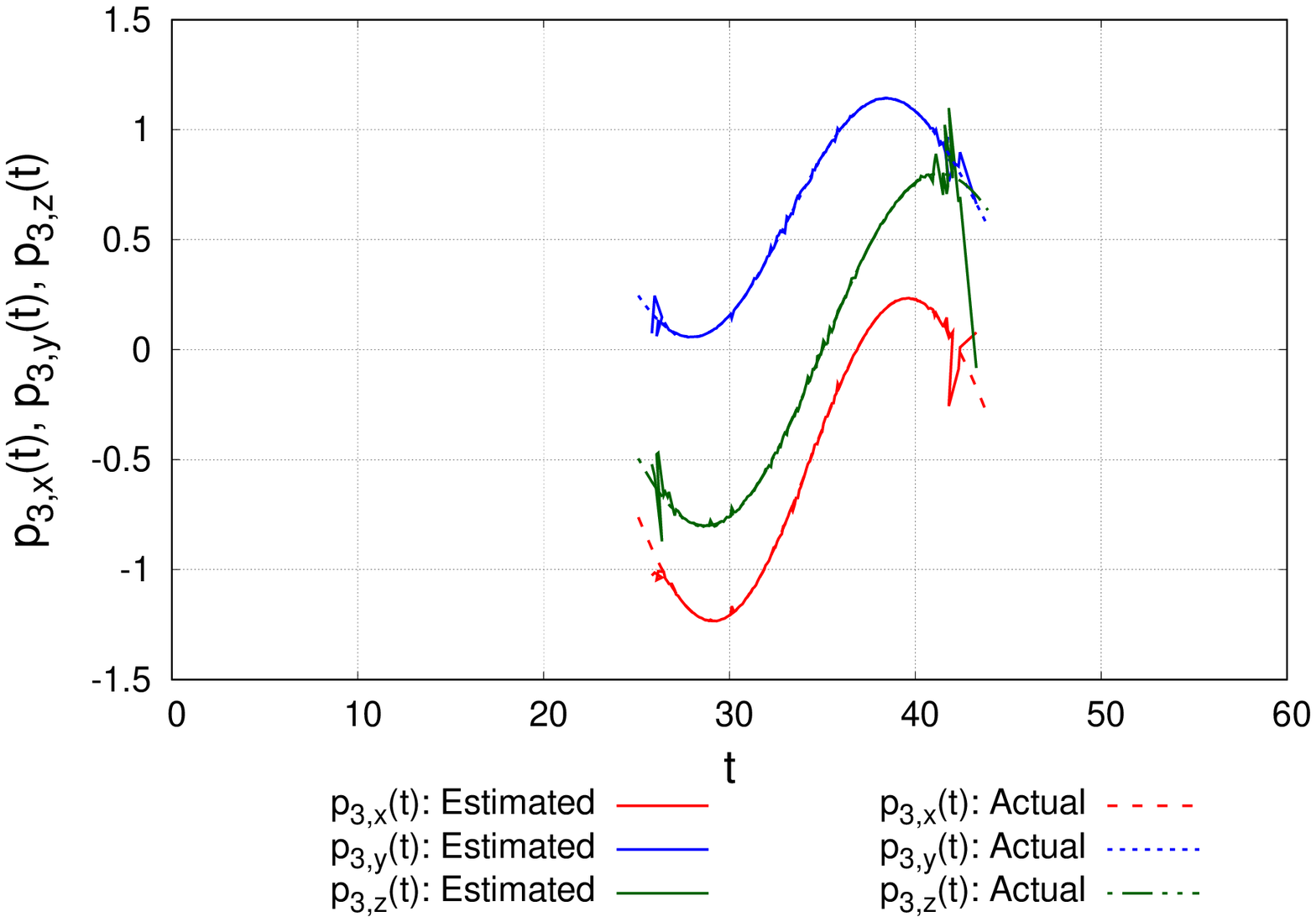}}}
  \subfloat[magnitude of sources 3]{
  \resizebox*{8cm}{!}{\includegraphics{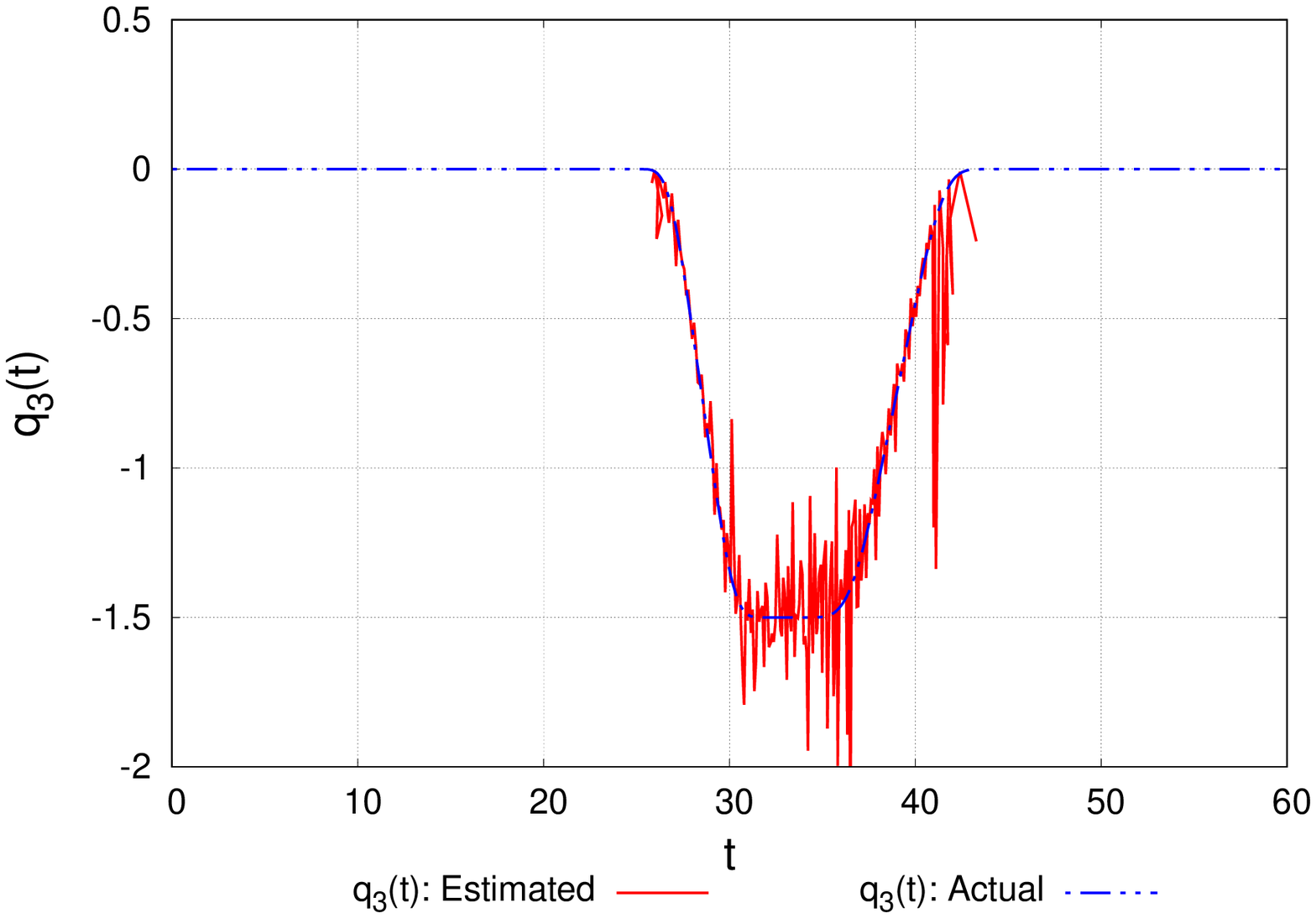}}}
\caption{Estimated locations and magnitudes of point sources for observation data with 0.5\% noise.}
\label{fig:reconstruction_results_loc_nl_0.5_point_source}
\end{figure}
\begin{figure}[ht]
\centering
  \subfloat[location of source 1]{
  \resizebox*{8cm}{!}{\includegraphics{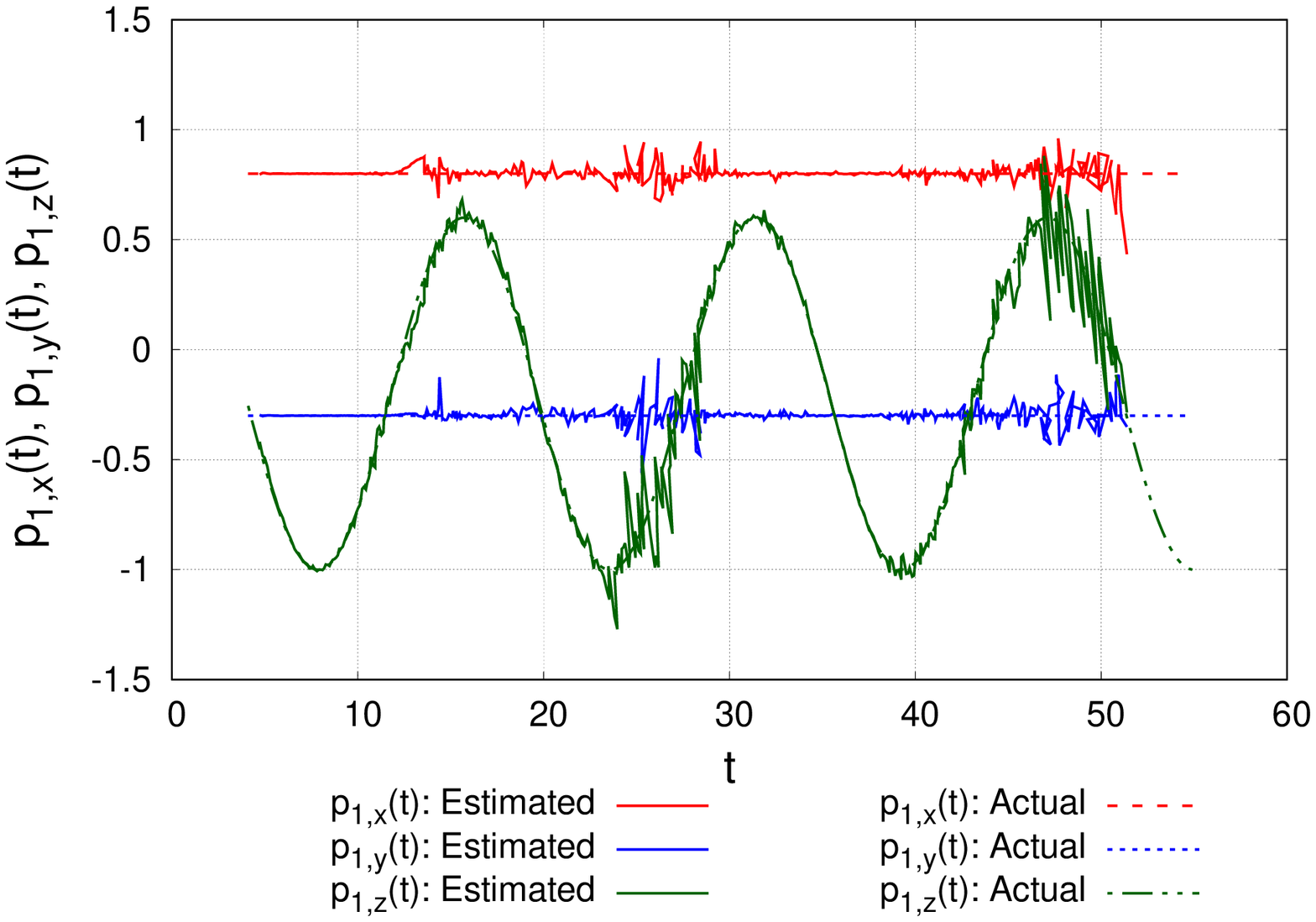}}}
  \subfloat[magnitude of source 1]{
  \resizebox*{8cm}{!}{\includegraphics{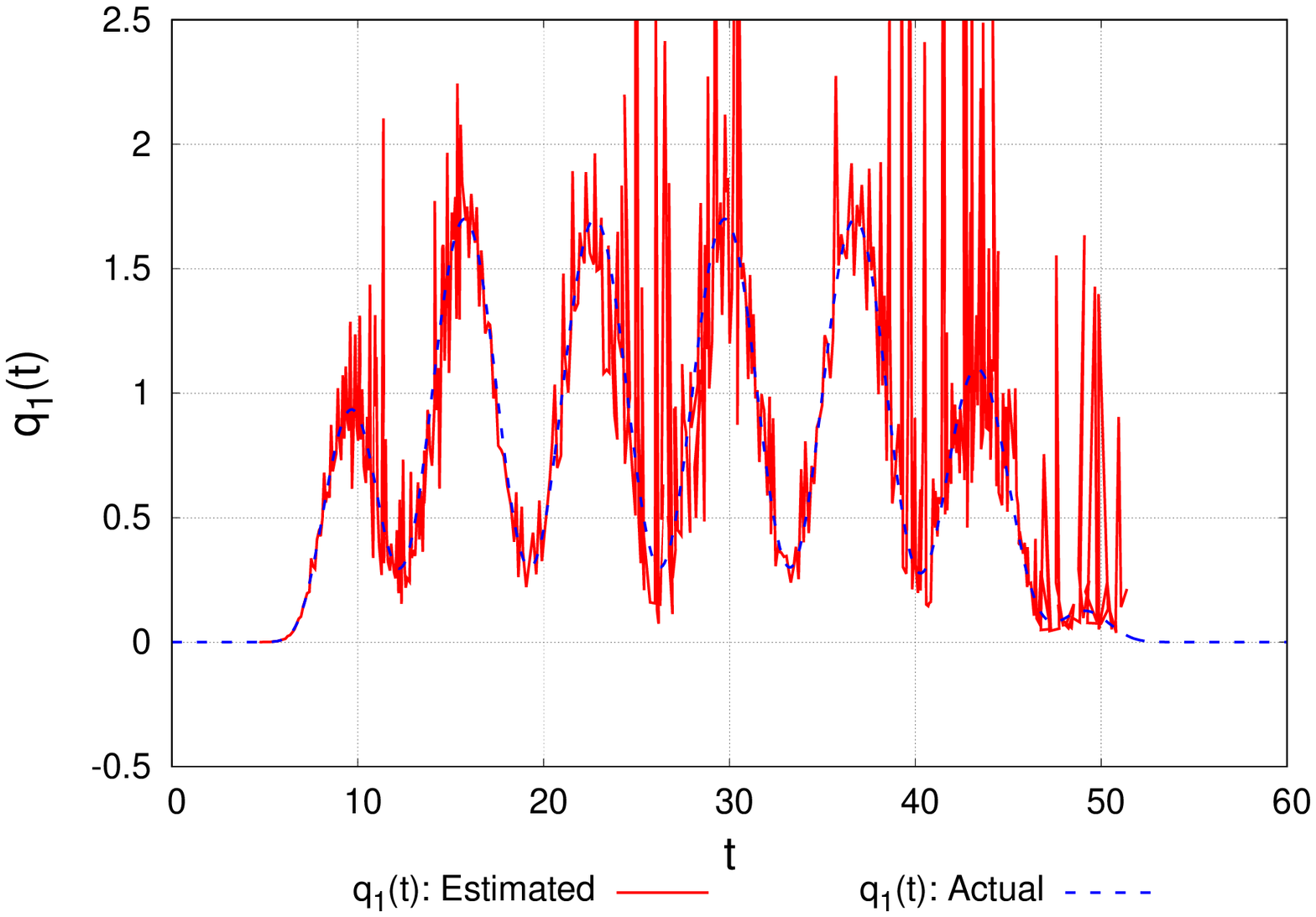}}}
  \\[-2ex]
  \subfloat[location of source 2]{
  \resizebox*{8cm}{!}{\includegraphics{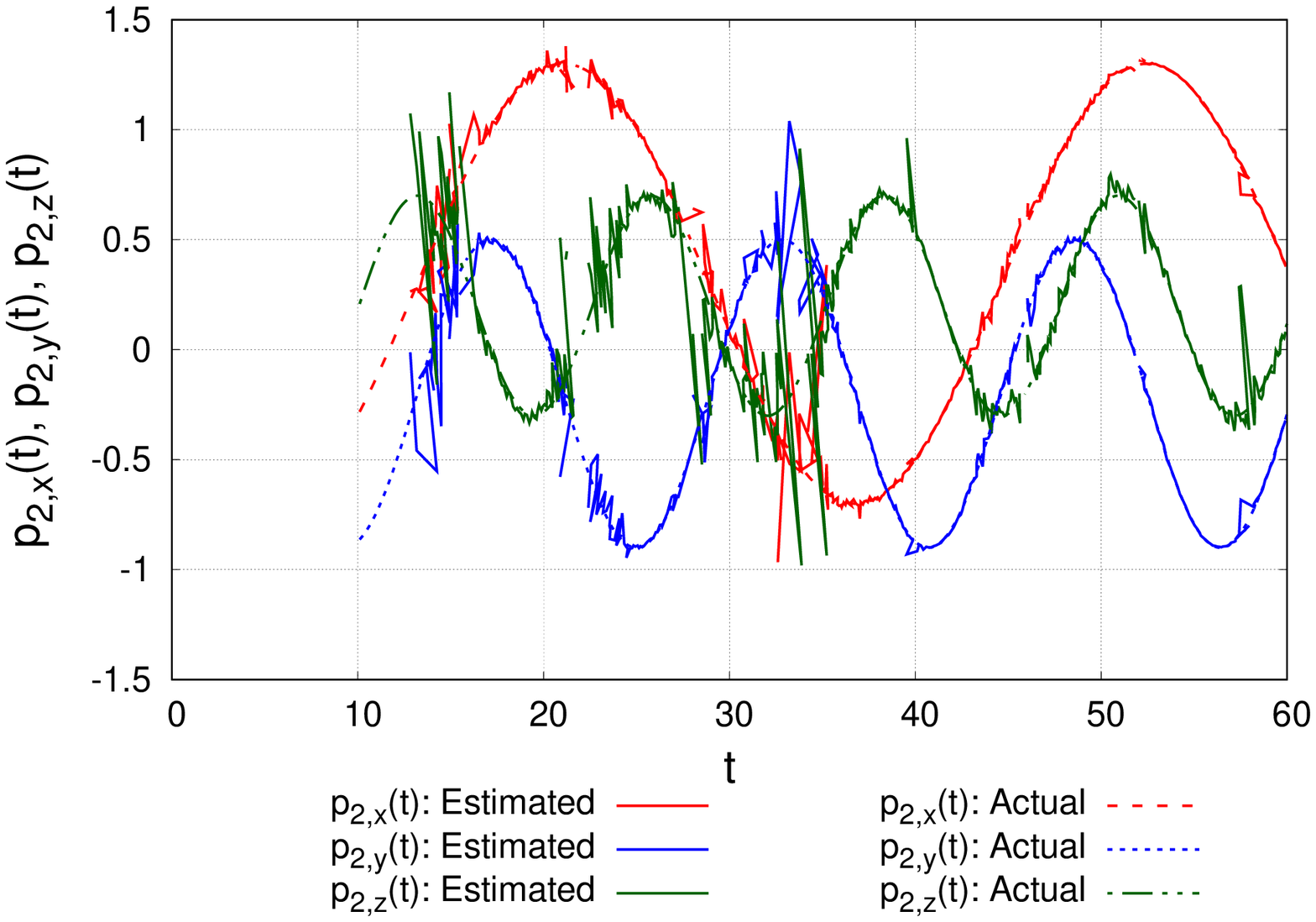}}}
  \subfloat[magnitude of source 2]{
  \resizebox*{8cm}{!}{\includegraphics{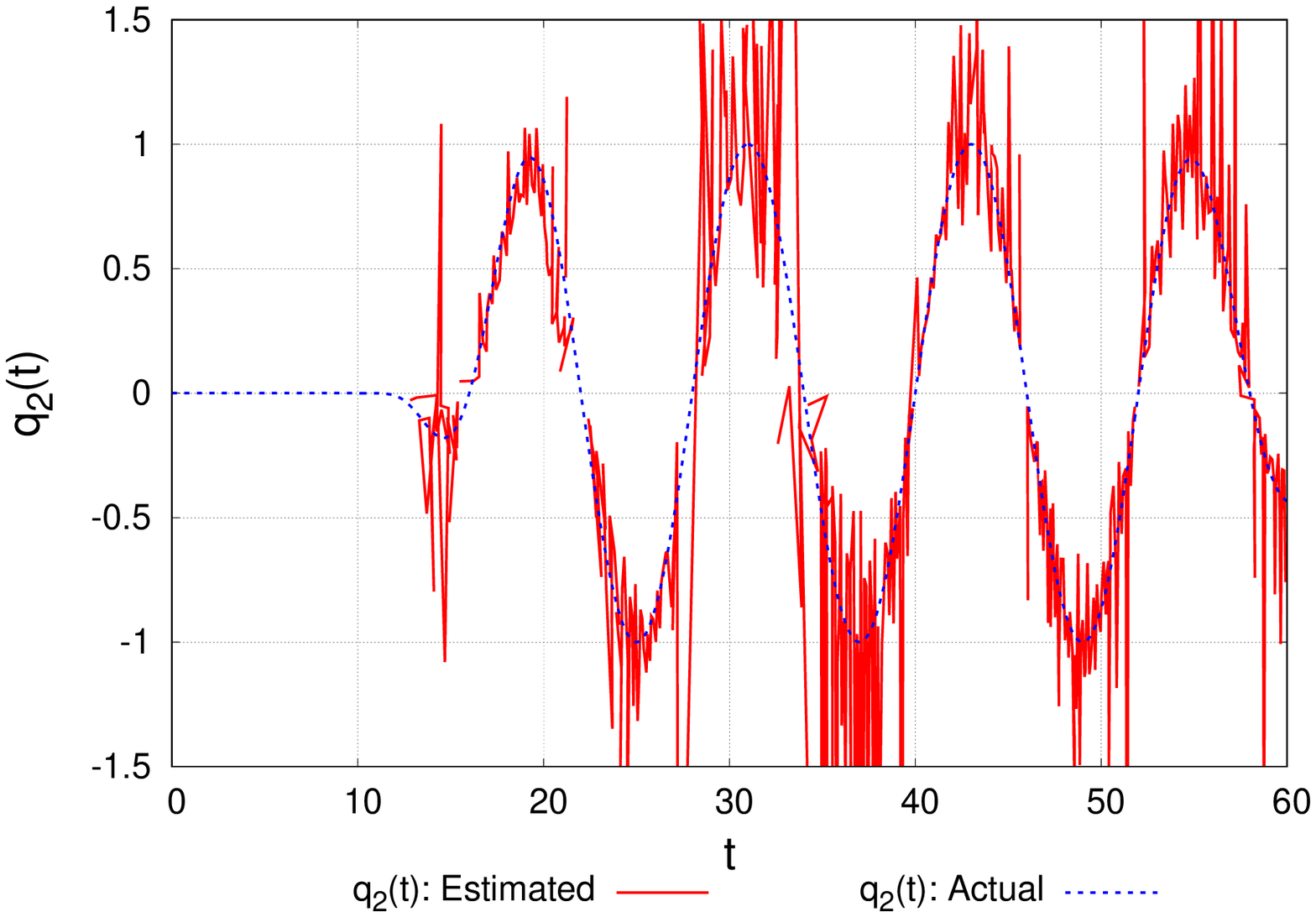}}}
   \\[-2ex]
  \subfloat[location of source 3]{
  \resizebox*{8cm}{!}{\includegraphics{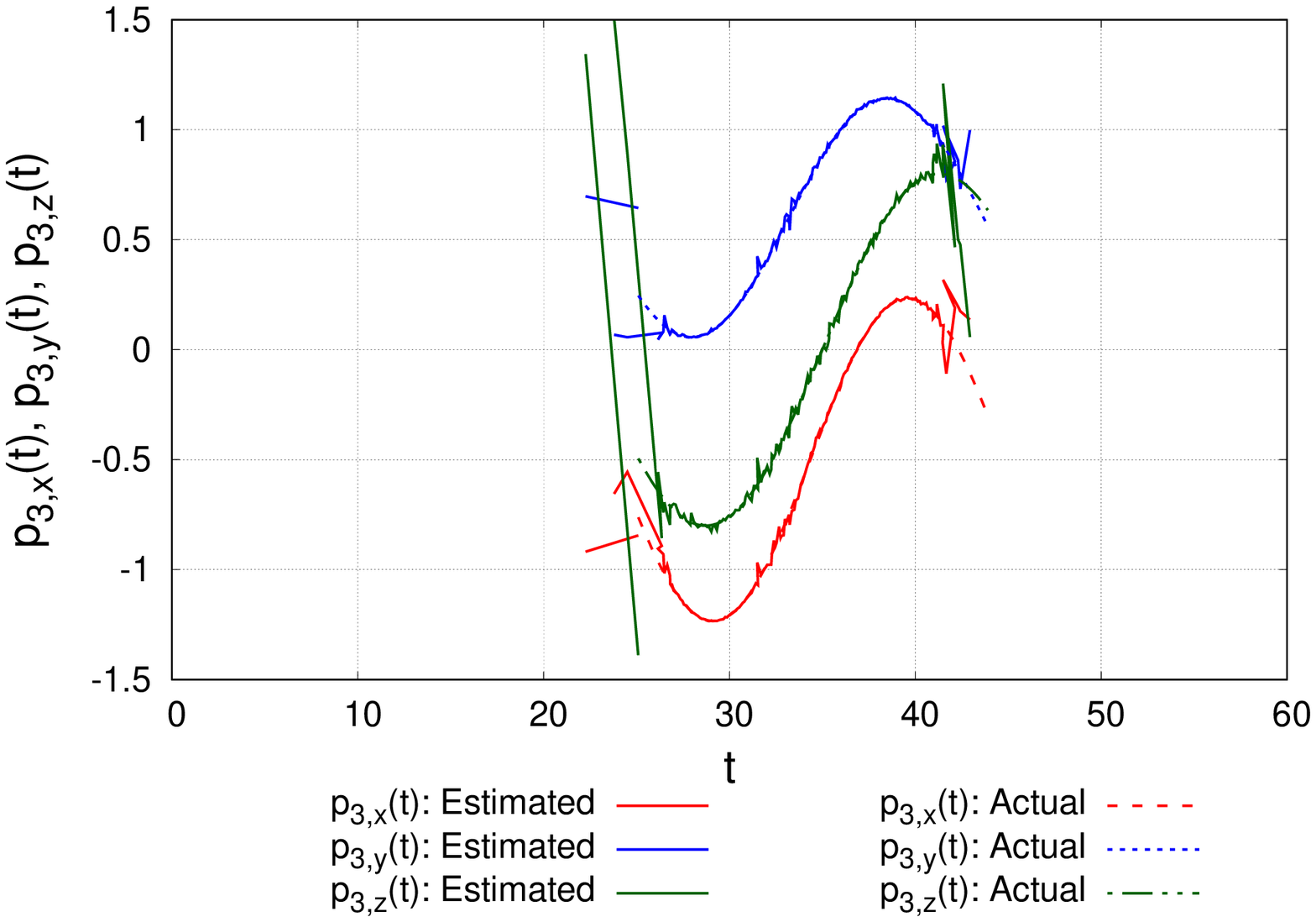}}}
  \subfloat[magnitude of source 3]{
  \resizebox*{8cm}{!}{\includegraphics{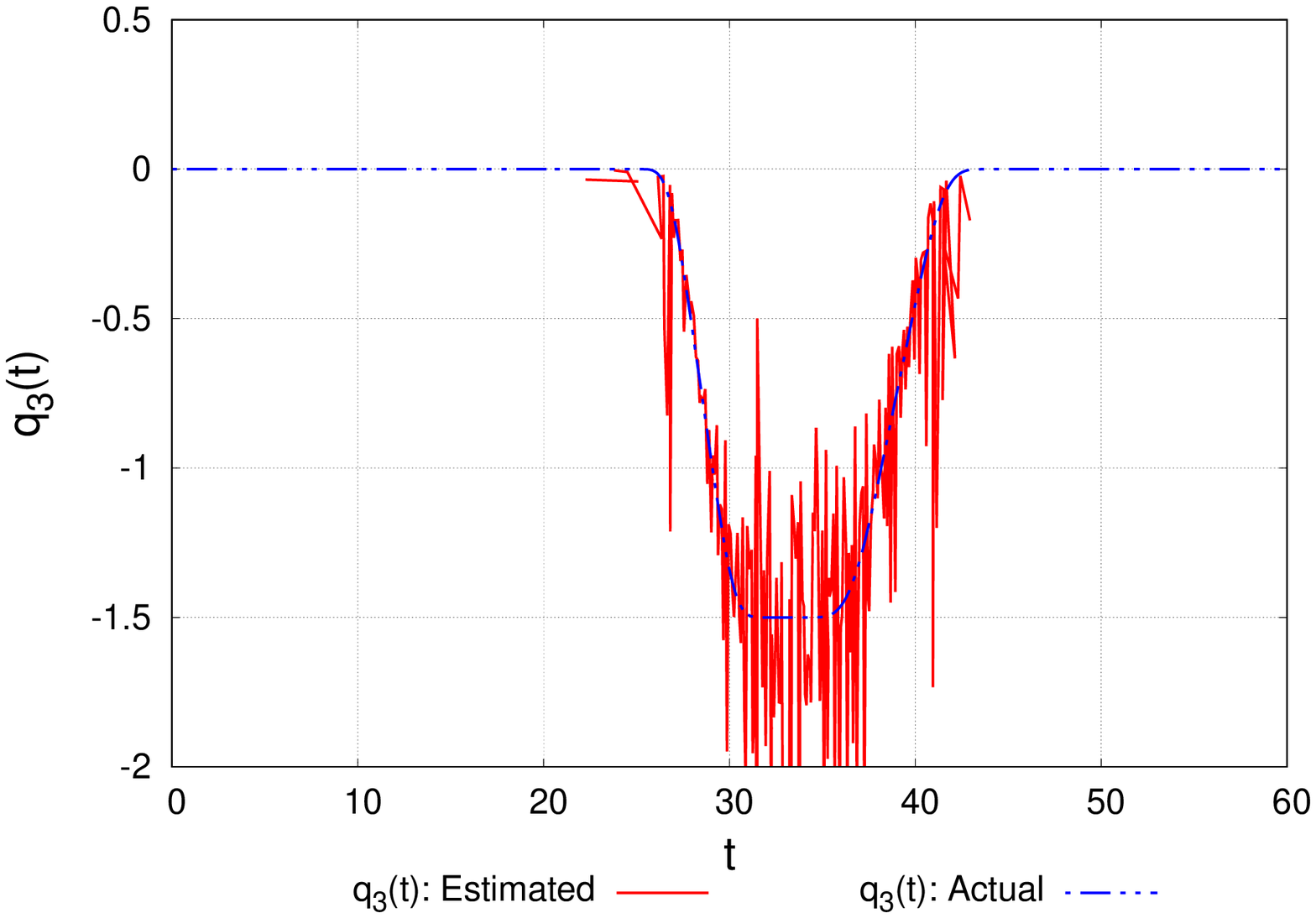}}}
\caption{Estimated locations and magnitudes of point sources for observation data with 1.0\% noise.}
\label{fig:reconstruction_results_loc_nl_1.0_point_source}
\end{figure}

Now, we show the reconstruction results of unknown parameters of moving point sources.
Figures
\ref{fig:reconstruction_results_loc_nl_0.0_point_source}-\ref{fig:reconstruction_results_loc_nl_1.0_point_source} display the reconstruction results of
locations and magnitudes for observation data without noise, and with 0.1\%, 0.5\% and 1.0\%  noise.
In these figures, 
we use the time variable $t$ for the horizontal axis 
by plotting $(t_k(\tau), \boldsymbol{p}_k(t_k(\tau)))$ and $(t_k(\tau), q_k(t_k(\tau)))$.
In Table \ref{table:Average_Error_point_source}, we show the average errors of estimated locations and magnitudes for observations
without noise, and with 0.1\%, 0.5\%, 1.0\%, and 5.0\% noise.
Here we define 
\begin{align*}
\mbox{average error of location }\boldsymbol{p}_k \equiv
& \displaystyle
\sqrt{\frac{1}{|I|}\int_{I} |\widehat{\boldsymbol{p}}_k(t_k(\tau)) -
 \boldsymbol{p}_k(t_k(\tau))|^2 {\mathrm d}\tau},\\
\mbox{average error of magnitude } q_k \equiv
& \displaystyle
\sqrt{\frac{1}{|I|}\int_{I}|\widehat{q}_k(t_k(\tau)) -
 q_k(t_k(\tau))|^2 {\mathrm d}\tau},
\end{align*}
where $\widehat{\boldsymbol{p}}_k$ and $\widehat{q}_k$ denote the estimated
location and magnitude of $k$-th point source, respectively, $I$ denotes the interval
in which we evaluate the average error, and $|I|$ denotes the measure of $I$.
\par
From results of Figures
\ref{fig:reconstruction_results_loc_nl_0.0_point_source}-\ref{fig:reconstruction_results_loc_nl_1.0_point_source}
and Table \ref{table:Average_Error_point_source}, we consider that our method works well if the
noise of observation data is smaller than 0.5\%.
We can observe that if noise exceeds 1\%, the influence of noise
becomes unignorable in reconstruction results.
Especially, estimations of magnitudes are highly affected by the
observation noise.
The reconstruction results become unreliable under the noise larger than 5\%.
\par
%
%
\begin{table}[p]
\tbl{The average errors of estimated locations and magnitudes in  each interval}
{\begin{tabular}{cccccc}
\toprule
 & $3.9 \leq \tau < 10.2$ & $10.2 \leq \tau < 24.6$ 
 & $24.6 \leq \tau < 44.6$ & $44.6 \leq \tau < 54.0$ & $54.0 \leq \tau < 60.0$ \\
 & $(K(\tau)=1)$ & $(K(\tau)=2)$ & $(K(\tau)=3)$ & $(K(\tau)=2)$ &
		     $(K(\tau)=1)$ \\ 
\midrule
(a) & \multicolumn{5}{c}{Without noise} \\
$\boldsymbol{p}_1$ & $1.7E-2$ & $3.9E-2$ & $2.2E-2$ & $8.1E-2$ & $-$ \\
$\boldsymbol{p}_2$ & $-$      & $3.6E-2$ & $4.0E-2$ & $3.1E-2$ & $1.3E-1$ \\
$\boldsymbol{p}_3$ & $-$      & $-$      & $1.3E-2$ & $-$      & $-$ \\
$q_1$      & $1.3E-3$ & $7.6E-2$ & $3.2E-2$ & $5.1E-2$ & $-$ \\
$q_2$      & $-$      & $1.1E-1$ & $5.0E-2$ & $3.2E-2$ & $1.6E-2$ \\
$q_3$      & $-$      & $-$      & $2.5E-2$ & $-$      & $-$ \\ 
\midrule
(b) & \multicolumn{5}{c}{With $0.1\%$ noise} \\
$\boldsymbol{p}_1$ & $1.7E-2$ & $3.9E-2$ & $2.3E-2$ & $8.3E-2$ & $-$ \\
$\boldsymbol{p}_2$ & $-$      & $5.0E-2$ & $3.2E-2$ & $3.6E-2$ & $1.5E-1$ \\
$\boldsymbol{p}_3$ & $-$      & $-$      & $3.8E-2$ & $-$      & $-$ \\
$q_1$      & $1.4E-2$ & $8.0E-2$ & $4.5E-2$ & $5.2E-2$ & $-$ \\
$q_2$      & $-$      & $1.2E-1$ & $5.6E-2$ & $3.8E-2$ & $2.7E-2$ \\
$q_3$      & $-$      & $-$      & $4.9E-2$ & $-$      & $-$ \\ 
\midrule
(c) & \multicolumn{5}{c}{With $0.5\%$ noise} \\
$\boldsymbol{p}_1$ & $1.9E-2$ & $4.7E-2$ & $6.3E-2$ & $1.9E-1$ & $-$ \\
$\boldsymbol{p}_2$ & $-$      & $2.2E-1$ & $1.7E-1$ & $4.1E-2$ & $3.7E-2$ \\
$\boldsymbol{p}_3$ & $-$      & $-$      & $8.8E-2$ & $-$      & $-$ \\
$q_1$      & $8.1E-2$ & $2.6E-1$ & $4.7E-1$ & $2.6E-1$ & $-$ \\
$q_2$      & $-$      & $2.1E-1$ & $6.6E-1$ & $1.4E-1$ & $2.0E-1$ \\
$q_3$      & $-$      & $-$      & $2.1E-1$ & $-$      & $-$ \\ 
\midrule
(d) & \multicolumn{5}{c}{With $1.0\%$ noise} \\
$\boldsymbol{p}_1$ & $1.6E-2$ & $7.9E-2$ & $8.8E-2$ & $2.1E-1$ & $-$ \\
$\boldsymbol{p}_2$ & $-$      & $2.7E-1$ & $2.4E-1$ & $5.9E-2$ & $1.1E-1$ \\
$\boldsymbol{p}_3$ & $-$      & $-$      & $2.5E-1$ & $-$      & $-$ \\
$q_1$      & $1.5E-1$ & $5.2E-1$ & $1.2E+0$ & $4.1E-1$ & $-$ \\
$q_2$      & $-$      & $2.7E-1$ & $1.0E+0$ & $4.5E-1$ & $6.0E-1$ \\
$q_3$      & $-$      & $-$      & $7.3E-1$ & $-$       & $-$ \\ 
\midrule
(e) & \multicolumn{5}{c}{With $5.0\%$ noise} \\
$\boldsymbol{p}_1$ & $5.7E-2$ & $2.1E-1$ & $2.2E-1$ & $7.1E-1$ & $-$ \\
$\boldsymbol{p}_2$ & $-$      & $6.0E-1$ & $5.0E-1$ & $2.5E-1$ & $1.9E-1$ \\
$\boldsymbol{p}_3$ & $-$      & $-$      & $2.9E-1$ & $-$      & $-$ \\
$q_1$      & $1.9E+0$ & $3.4E+0$ & $3.9E+0$ & $1.6E+0$ & $-$ \\
$q_2$      & $-$      & $2.1E+0$ & $3.7E+0$ & $2.9E+0$ & $3.6E+0$ \\
$q_3$      & $-$      & $-$      & $3.8E+0$ & $-$      & $-$ \\ 
\bottomrule
\end{tabular}}
\label{table:Average_Error_point_source}
\end{table}
%
%
%
%
%
%
%
%
%
%
%
%
%
%
%
%
%
%
%
\par
We can find that the reconstruction of location becomes bad for the source with small magnitude, e.\@g.\@ the result of source 2 around $t=12$ in
Figures \ref{fig:reconstruction_results_loc_nl_0.5_point_source}(b) and \ref{fig:reconstruction_results_loc_nl_1.0_point_source}(b), but we consider that these results are reasonable.
We can also find bad reconstruction results due to small magnitude
for source 2 in the interval $54.0 \leq \tau <
60.0$ for the observation without noise and 1.0\% noise in Table \ref{table:Average_Error_point_source}, and bad
identification of $K(\tau)$ around $\tau=12$ and $\tau=53$ in Figure
\ref{fig:behaviour_K_point_source}(b).
\vspace{2ex}
\par
\noindent
\textbf{Note 4.1.}
In Table \ref{table:Average_Error_point_source}, the error of
$\boldsymbol{p}_2$ in the interval $54.0 \leq
\tau \leq 60.0$ is concentrated in $57.6 \leq \tau \leq 57.8$ where $q_2(t_2(\tau)) \simeq 0$.
Except of these bad estimates, the average error of $\boldsymbol{p}_2$ in $54.0 \leq \tau <
60.0$ becomes $1.56E-2$ for the case  without noise and $1.45E-2$ for the
case with 0.1\% noise.
\vspace{2ex}
\par
\begin{figure}[ht]
\centering
\subfloat{
  \resizebox{10cm}{!}{\includegraphics{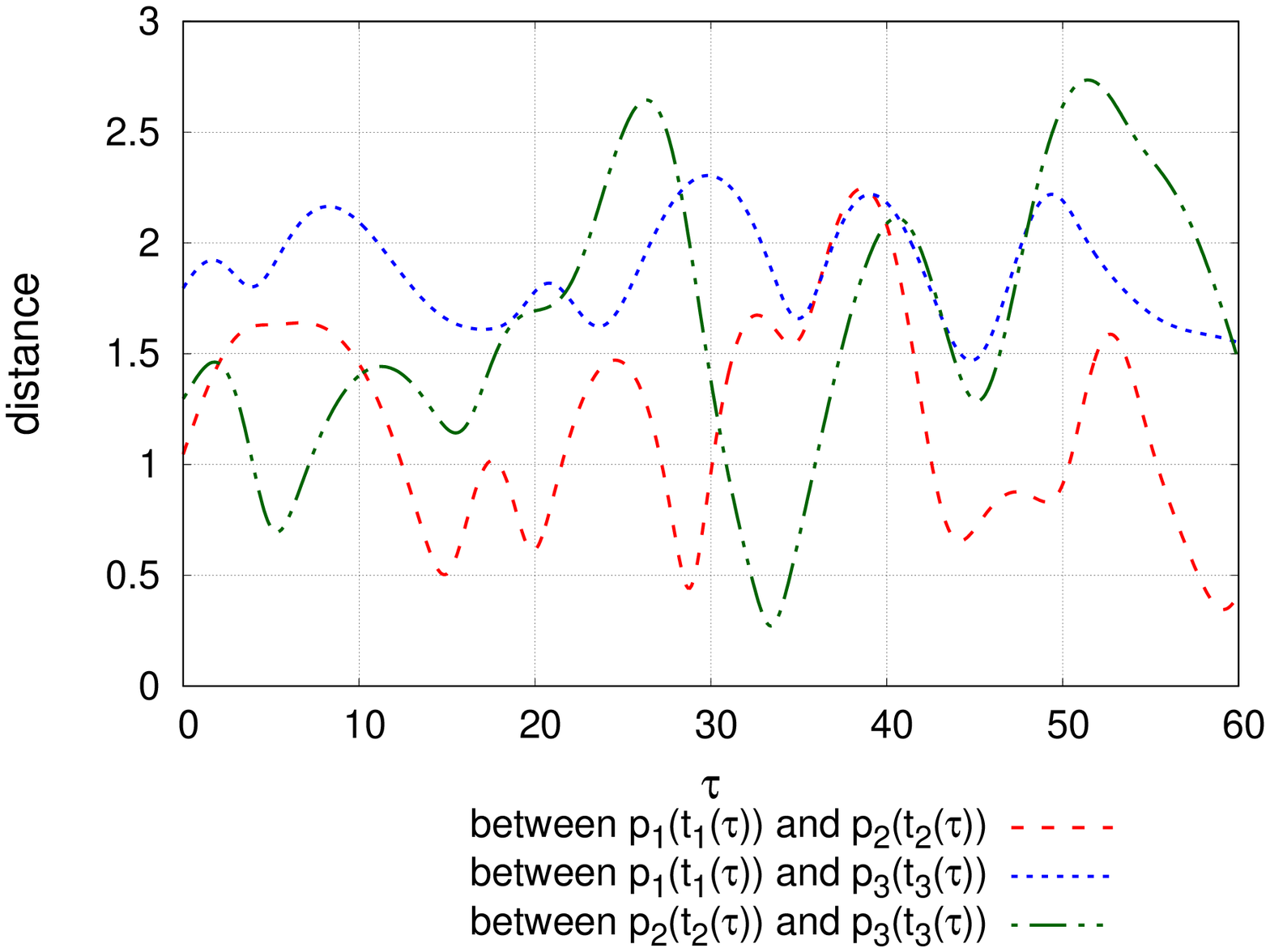}}}
\\[4ex]
  \caption{Behaviour of the distances between each two point sources.}
  \label{fig:distance_locations}
\end{figure}

In Figures
\ref{fig:reconstruction_results_loc_nl_0.5_point_source} and \ref{fig:reconstruction_results_loc_nl_1.0_point_source},
we can also find bad reconstruction results even for large magnitudes of point sources, e.\@g.\@ results for sources 2 and 3 at
around $t=30$.
To discuss the reason of these bad estimates, we display the behaviour
of distances between each two point sources in Figure \ref{fig:distance_locations}.
The results of Figures
\ref{fig:reconstruction_results_loc_nl_0.5_point_source}-\ref{fig:distance_locations}
suggest that if the arangement of point sources is well
separated, then our method works well and gives good reconstruction
results, however, as arrangement of point sources becomes closer, then 
it becomes more difficult to distinguish these sources, and our method gives bad estimations.
\par
\par
In Figure \ref{fig:behaviour_K_point_source}(b), we can find many wrong
identification results of $K(\tau)$, especially, in  $25 \leq \tau \leq
38$ where $K(\tau)$ is identified as 4.
However, the absolute value of magnitude of
reconstructed 4th source is very small
($<0.01$)  comparing to other sources, and so the 4th source can be
recognised as `noises' or `ghosts'.
\par
\vspace{2ex}
\par
\noindent
\textbf{Note 4.2.}
Since our reconstruction procedure works at every
moment $\tau$, we have to identify whether the point source
reconstructed at the moment $\tau$ belongs to which source
reconstructed at the previous moment $\tau-\Delta\tau$, especially,
when $K(\tau)$ changes.
To solve this problem, we measure distances between the reconstructed
sources at moments $\tau$ and $\tau-\Delta\tau$, and identify as the same
source which has the minimum distance.
\par
\vspace{2ex}
\par
Throughout of the numerical results, we can observe that the
observation noise affects to the estimation of magnitudes more heavily
than the one of locations.
This phenomena is mainly caused by the estimation of perturbation term $\xi_k(t_k(\tau))$.
As an alternative method,  we apply Step 4' in Note 3.5 to obtain
$dp_{k,z}(t_k(\tau))/d\tau$, and compare to
results by Step 4.
Here, we apply the central differences 
\begin{equation*}
 \frac{dp_{k,z}(t_k(\tau))}{d\tau}
\sim \frac{p_{k,z}(\tau+\Delta \tau)
-p_{k,z}(\tau-\Delta \tau) }{2\cdot\Delta\tau}.
\end{equation*}
to estimate $dp_{k,z}(t_k(\tau))/d\tau$, where $p_{k,z}(t_k(\tau))$ is
obtained in Step 3.
Table \ref{table:Average_Error_point_source_another_idea} shows the
average errors of estimated magnitudes using Step 4' instead of Step 4.
Comparing  Tables \ref{table:Average_Error_point_source} and
\ref{table:Average_Error_point_source_another_idea}, the errors are
almost the same in both methods for the case where the observation noise is smaller than
0.5\%. 
Hence we consider that both method can be used as an alternative if the
noise is small.
\par
%
\begin{table}[p]
\tbl{The average errors of estimated magnitudes in  each interval 
using Step 4' instead of Step 4.}
{\begin{tabular}{cccccc}
\toprule
 & $3.9 \leq \tau < 10.2$ & $10.2 \leq \tau < 24.6$ 
 & $24.6 \leq \tau < 44.6$ & $44.6 \leq \tau < 54.0$ & $54.0 \leq \tau < 60.0$ \\
 & $(K(\tau)=1)$ & $(K(\tau)=2)$ & $(K(\tau)=3)$ & $(K(\tau)=2)$ &
		     $(K(\tau)=1)$ \\ 
\midrule
(a) & \multicolumn{5}{c}{Without noise} \\
$q_1$      & $1.8E-2$ & $1.5E-1$ & $9.5E-2$ & $8.5E-2$ & $-$ \\
$q_2$      & $-$      & $1.5E-2$ & $4.3E-2$ & $5.5E-2$ & $5.0E-2$ \\
$q_3$      & $-$      & $-$      & $3.4E-2$ & $-$      & $-$ \\ 
\midrule
(b) & \multicolumn{5}{c}{With $0.1\%$ noise} \\
$q_1$      & $1.8E-2$ & $1.6E-1$ & $9.6E-2$ & $9.0E-2$ & $-$ \\
$q_2$      & $-$      & $5.0E-2$ & $6.0E-2$ & $5.8E-2$ & $5.4E-2$ \\
$q_3$      & $-$      & $-$      & $4.0E-2$ & $-$      & $-$ \\ 
\midrule
(c) & \multicolumn{5}{c}{With $0.5\%$ noise} \\
$q_1$      & $3.6E-2$ & $2.6E-1$ & $1.5E-1$ & $1.6E-1$ & $-$ \\
$q_2$      & $-$      & $2.0E-1$ & $2.6E-1$ & $8.8E-2$ & $8.8E-2$ \\
$q_3$      & $-$      & $-$      & $9.3E-2$ & $-$      & $-$ \\ 
\midrule
(d) & \multicolumn{5}{c}{With $1.0\%$ noise} \\
$q_1$      & $5.0E-2$ & $5.1E-1$ & $2.9E-1$ & $2.4E-1$ & $-$ \\
$q_2$      & $-$      & $4.2E-1$ & $6.0E-1$ & $1.4E-1$ & $1.3E-1$ \\
$q_3$      & $-$      & $-$      & $1.8E-1$ & $-$       & $-$ \\ 
\midrule
(e) & \multicolumn{5}{c}{With $5.0\%$ noise} \\
 $q_1$      & $2.5E-1$ & $1.2E+0$ & $8.0E-1$ & $5.5E-1$ & $-$ \\
 $q_2$      & $-$      & $1.6E+0$ & $1.8E+0$ & $6.2E-1$ & $7.0E-1$ \\
 $q_3$      & $-$      & $-$      & $6.2E-1$ & $-$      & $-$ \\ 
\bottomrule
\end{tabular}}
\label{table:Average_Error_point_source_another_idea}
\end{table}
\par
\subsection{Reconstruction of moving dipole sources}
Next, we demonstrate numerical experiments for reconstruction of moving dipole sources.
We consider the case where three dipole sources move in the domain $\Omega$.
The locations of dipoles are arranged as same as the previous example,
and moments of dipoles change as follows:
\begin{description}
\item[source 1.] ${\boldsymbol{m}}_1(t) = (m_{1,x}(t),\ m_{1,y}(t),\ 0)$ where
\begin{align*}
  m_{1,x}(t) &= \left\{ \begin{array}{ll}
                   0, & 0 \leq t < 2, \\[1ex]
                   \displaystyle 
                   \frac{1}{2} \eta\left(\frac{t-2}{7}\right) \cdot \cos\left(\frac{2\pi t}{10}\right),
                   & 2 \leq t < 9, \\[3ex]
                   \displaystyle 
                   \frac{1}{2} \cos\left(\frac{2\pi t}{10}\right),
                   & 9 \leq t < 32, \\[3ex]
                   \displaystyle 
                   \frac{1}{2}\left(1 -  \eta\left(\frac{t-32}{9}\right)\right)  \cdot 
		   \cos\left(\frac{2\pi t}{10}\right),
		   & 32 \leq t < 41, \\[3ex]
                   0, & 41 \leq t \leq 70,
                   \end{array}
                   \right. \\
  m_{1,y}(t) &= \left\{ \begin{array}{ll}
                   0, & 0 \leq t < 2, \\[1ex]
                   \displaystyle 
                   \frac{1}{2}\eta\left(\frac{t-2}{7}\right) \cdot \sin\left(\frac{2\pi t}{10}\right),
                   & 2 \leq t < 9, \\[3ex]
                   \displaystyle 
                   \frac{1}{2} \sin\left(\frac{2\pi t}{10}\right),
                   & 9 \leq t < 32, \\[3ex]
                   \displaystyle 
                   \frac{1}{2} \left(1 -  \eta\left(\frac{t-32}{9}\right)\right)  \cdot 
		   \sin\left(\frac{2\pi t}{10}\right),
		   & 32 \leq t < 41, \\[3ex]
                   0, & 41 \leq t \leq 70,
                   \end{array}
                   \right. 
                   %
\end{align*}
\item[source 2.] ${\boldsymbol{m}}_2(t) = (m_{2,x}(t),\ m_{2,y}(t),\ 0)$ where
\begin{align*}
  m_{2,x}(t) &= \left\{ \begin{array}{ll}
                   0, & 0 \leq t < 7, \\[2ex]
                   \displaystyle 
                   \eta\left(\frac{t-7}{5}\right) \cdot 
                   \left( \frac{3}{5} - \frac{1}{5} \cos
		    \left(\frac{2\pi(t-7)}{15}\right)\right)\cdot
                    \sin\left(-\frac{2\pi t}{11}\right),
                   & 7 \leq t < 12, \\[3ex]
                   \displaystyle 
                   \left(\frac{3}{5} - \frac{1}{5} \cos
		    \left(\frac{2\pi(t-7)}{15}\right)\right)\cdot
                    \sin\left(-\frac{2\pi t}{11}\right),
                   & 12 \leq t \leq 70,
                   \end{array}
                   \right. \\[2ex]
  m_{2,y}(t) &= \left\{ \begin{array}{ll}
                   0, & 0 \leq t < 7, \\[2ex]
                   \displaystyle 
                   \eta\left(\frac{t-7}{5}\right) \cdot 
                   \left( \frac{3}{5} - \frac{1}{5} \cos
		    \left(\frac{2\pi(t-7)}{15}\right)\right)\cdot
                    \cos\left(-\frac{2\pi t}{11}\right),
                   & 7 \leq t < 12, \\[3ex]
                   \displaystyle 
                   \left(\frac{3}{5} - \frac{1}{5} \cos
		    \left(\frac{2\pi(t-7)}{15}\right)\right)\cdot
                    \cos\left(-\frac{2\pi t}{11}\right),
                   & 12\leq t \leq 70.
                   \end{array}
                   \right. 
\end{align*}

\item[source 3.]${\boldsymbol{m}}_3(t) = (m_{3,x}(t),\ m_{3,y}(t),\ 0)$ where
\begin{align*}
  m_{3,x}(t) &= \left\{ \begin{array}{ll}
                   0, & 0 \leq t < 18, \\[1ex]
                   \displaystyle 
                   \eta\left(\frac{t-18}{8}\right) \cdot \cos\left(\frac{2\pi t}{8}\right),
                   & 18 \leq t < 26, \\[3ex]
                   \displaystyle 
                   \cos\left(\frac{2\pi t}{8}\right),
                   & 26 \leq t < 45, \\[3ex]
                   \displaystyle 
                   \left(1 -  \eta\left(\frac{t-45}{9}\right)\right)  \cdot 
		   \cos\left(\frac{2\pi t}{8}\right),
		   & 45 \leq t < 54, \\[3ex]
                   0, & 54 \leq t \leq 70,
                   \end{array}
                   \right. \\[2ex]
  m_{3,y}(t) &= \left\{ \begin{array}{ll}
                   0, & 0 \leq t < 18, \\[1ex]
                   \displaystyle 
                   \eta\left(\frac{t-18}{8}\right) \cdot \sin\left(\frac{2\pi t}{8}\right),
                   & 18 \leq t < 26, \\[3ex]
                   \displaystyle 
                   \sin\left(\frac{2\pi t}{8}\right),
                   & 26 \leq t < 45, \\[3ex]
                   \displaystyle 
                   \left(1 -  \eta\left(\frac{t-45}{9}\right)\right)  \cdot 
		   \sin\left(\frac{2\pi t}{8}\right),
		   & 45 \leq t < 54, \\[3ex]
                   0, & 54 \leq t \leq 70.
                   \end{array}
                   \right. 
\end{align*}
\end{description}
Hence $K(\tau)$ changes as
\begin{displaymath}
     K(\tau) = \left\{\begin{array}{ll}
           0, & 0 \leq \tau <2.4, \\[1ex]
           1, & 2.4 \leq \tau < 6.8,\ \ 54.8 \leq \tau \leq 70.0, \\[1ex]
           2, & 6.8 \leq \tau < 17.3,\ \ 40.2 \leq \tau < 54.8, \\[1ex]
           3, & 17.3 \leq \tau < 40.2.
\end{array}\right.
\end{displaymath}
We note that $\boldsymbol{m}_k \in C^4([0,70]),\ k=1,2,3$.
%
\par
\begin{figure}[ht]
\centering
\subfloat{
  \resizebox*{9cm}{!}{\includegraphics{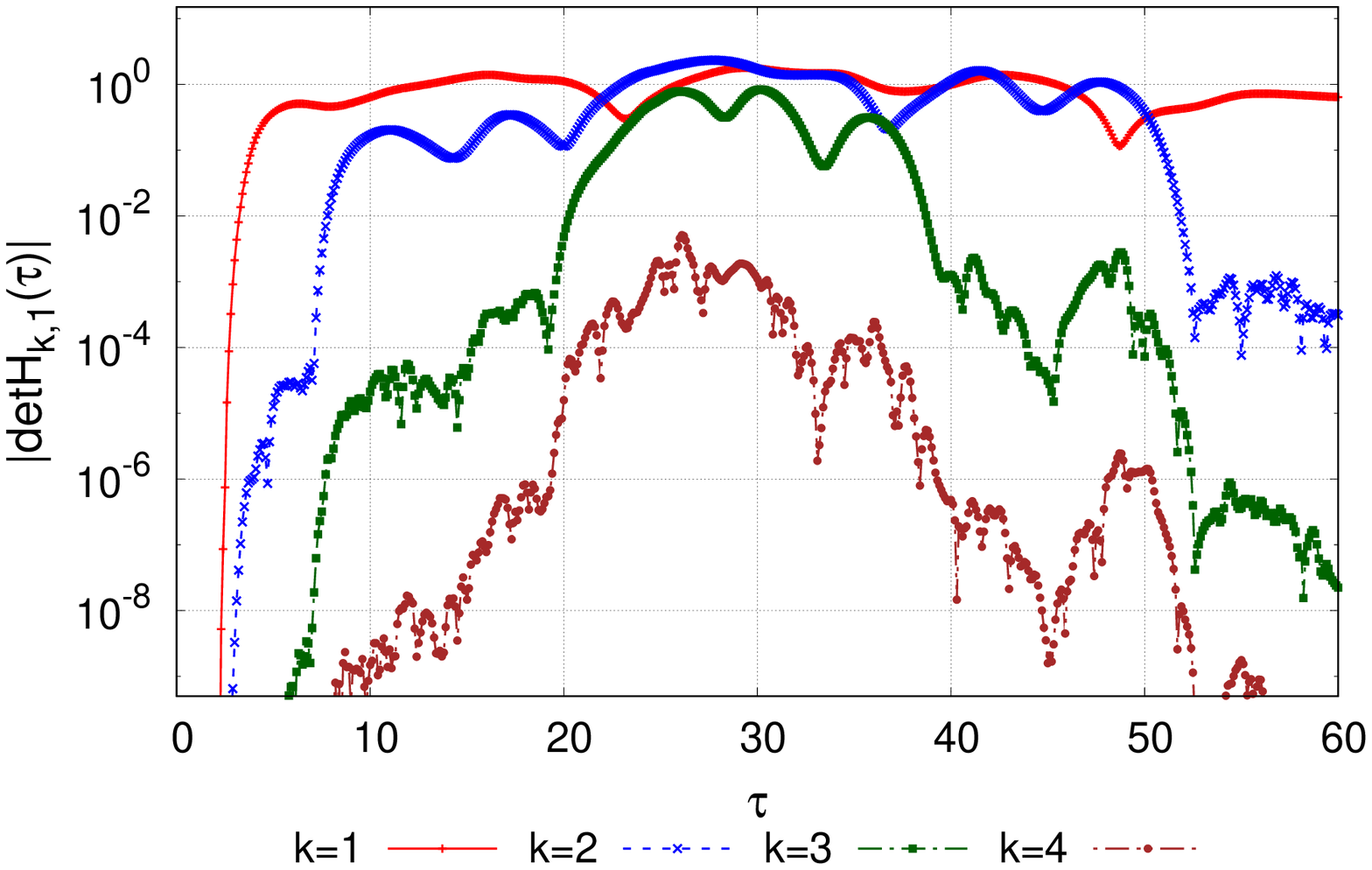}}}
  \\[4ex]
  \caption{Behaviour of $|\det H_{k,1}(\tau)|$ for observation data without noise: dipole source case.}
   \label{fig:behaviour_detH_dipole_source}
\end{figure}
\begin{figure}[ht]
\centering
  \subfloat[]{
  \resizebox*{9cm}{!}{\includegraphics{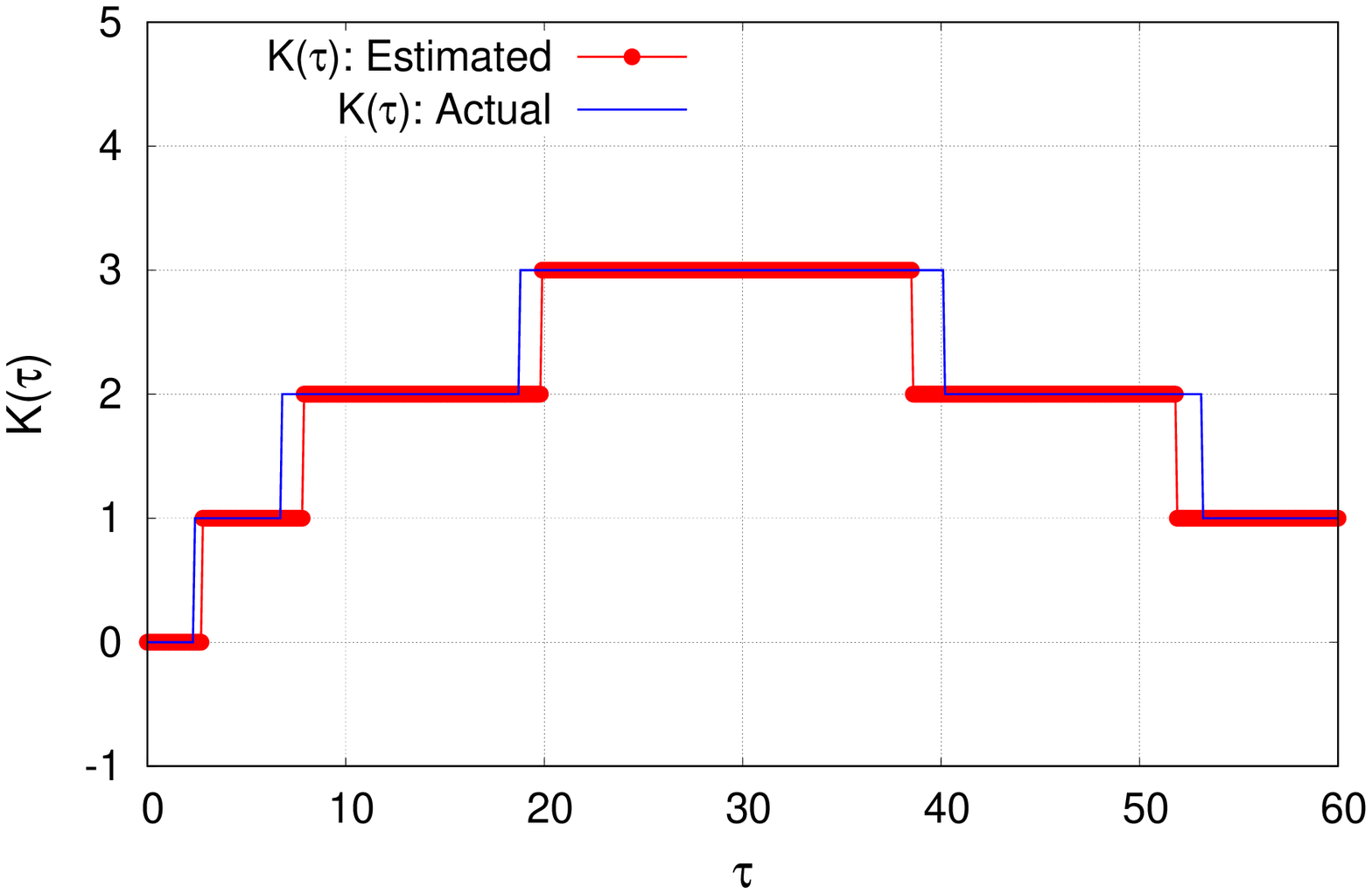}}}
  \\
  \subfloat[]{
  \resizebox*{9cm}{!}{\includegraphics{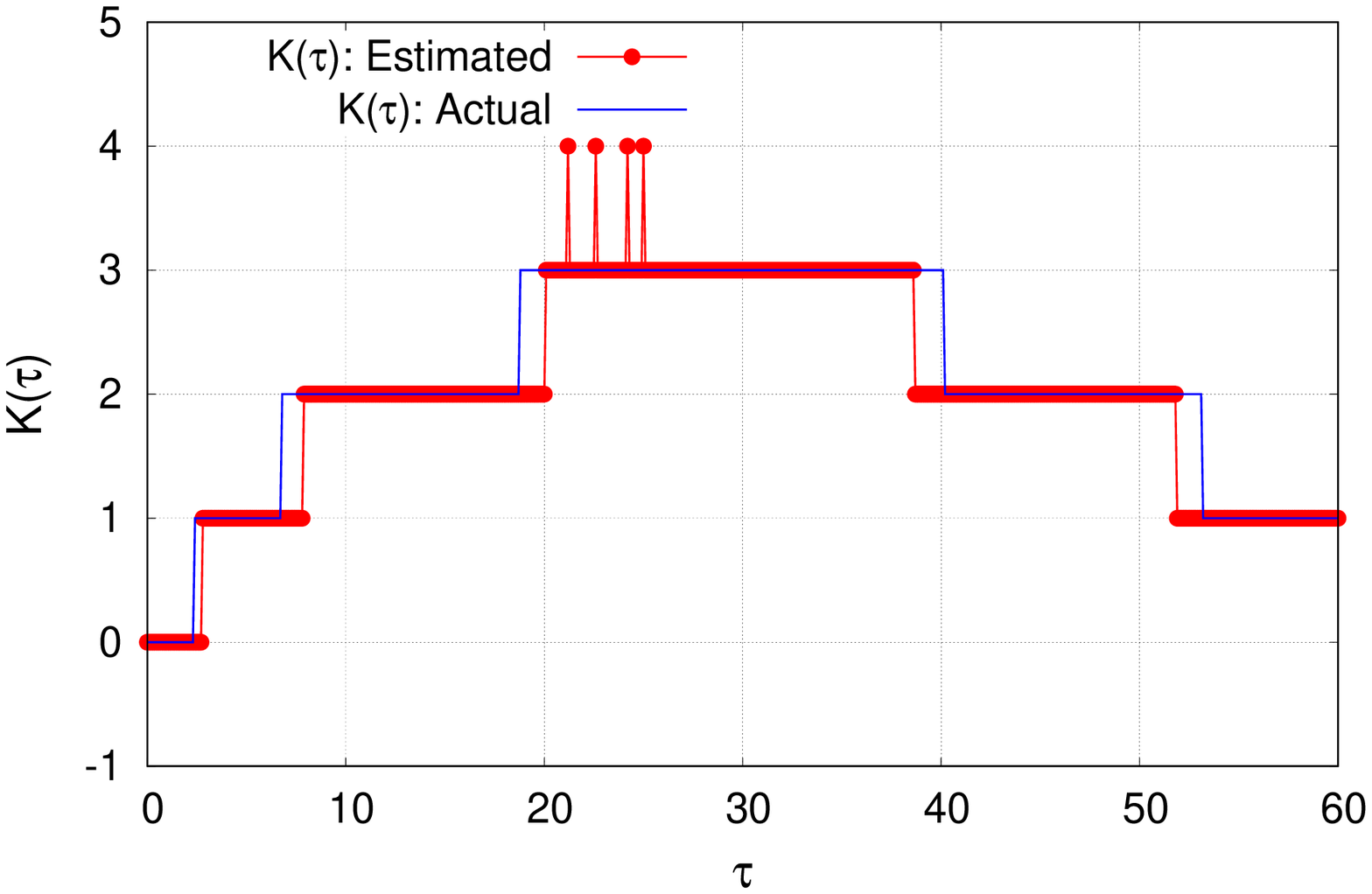}}}
  \caption{Behaviour of estimated number of dipole sources: (a) for
 observation data without noise, (b) with $0.5\%$ noise.}
  \label{fig:behaviour_K_dipole_source}
\end{figure}

We first discuss the identification of the number $K(\tau)$ of dipoles. 
Figure \ref{fig:behaviour_detH_dipole_source} shows the behaviour of $|\det H_{k,1}(\tau)|$ for $0 \leq \tau \leq
T_e$.
The behaviour of $|\det H_{k,1}(\tau)|$ is quite similar as for the one of
$|\det H_{k,0}(\tau)|$ for point
sources in Figure \ref{fig:behaviour_detH_point_source}, and hence we apply Algorithm 4.1 for the identification of
$K(\tau)$ using $H_{k,1}(\tau)$ instead of $H_{k,0}(\tau)$.
Figure \ref{fig:behaviour_K_dipole_source} displays the identification results of $K(\tau)$ from observations without noise and with 0.5\% noise.
Here, we set  the parameters $\varepsilon_0 = 1.0\times 10^{-4}$ and $\varepsilon_G =
2.5 \times 10^{-2}$ as same as for moving point sources.
From the result of Figure \ref{fig:behaviour_K_dipole_source}, 
we consider Algorithm 4.1  also works well for the identification of the
number of moving dipole sources.
\par
\begin{figure}[ht]
\centering
  \subfloat[location of source 1]{
  \resizebox*{8cm}{!}{\includegraphics{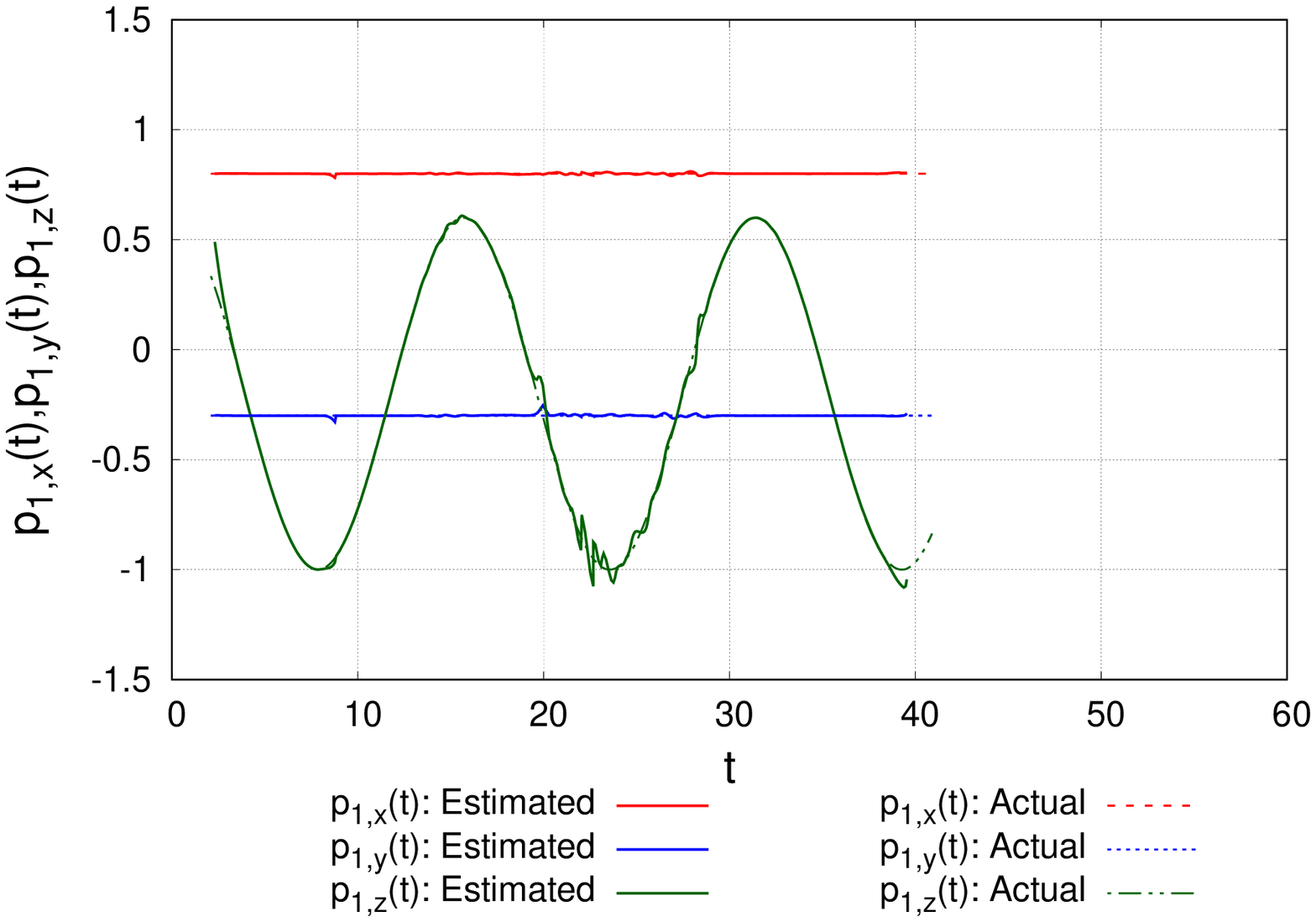}}}
  \subfloat[moment of sources 1]{
  \resizebox*{8cm}{!}{\includegraphics{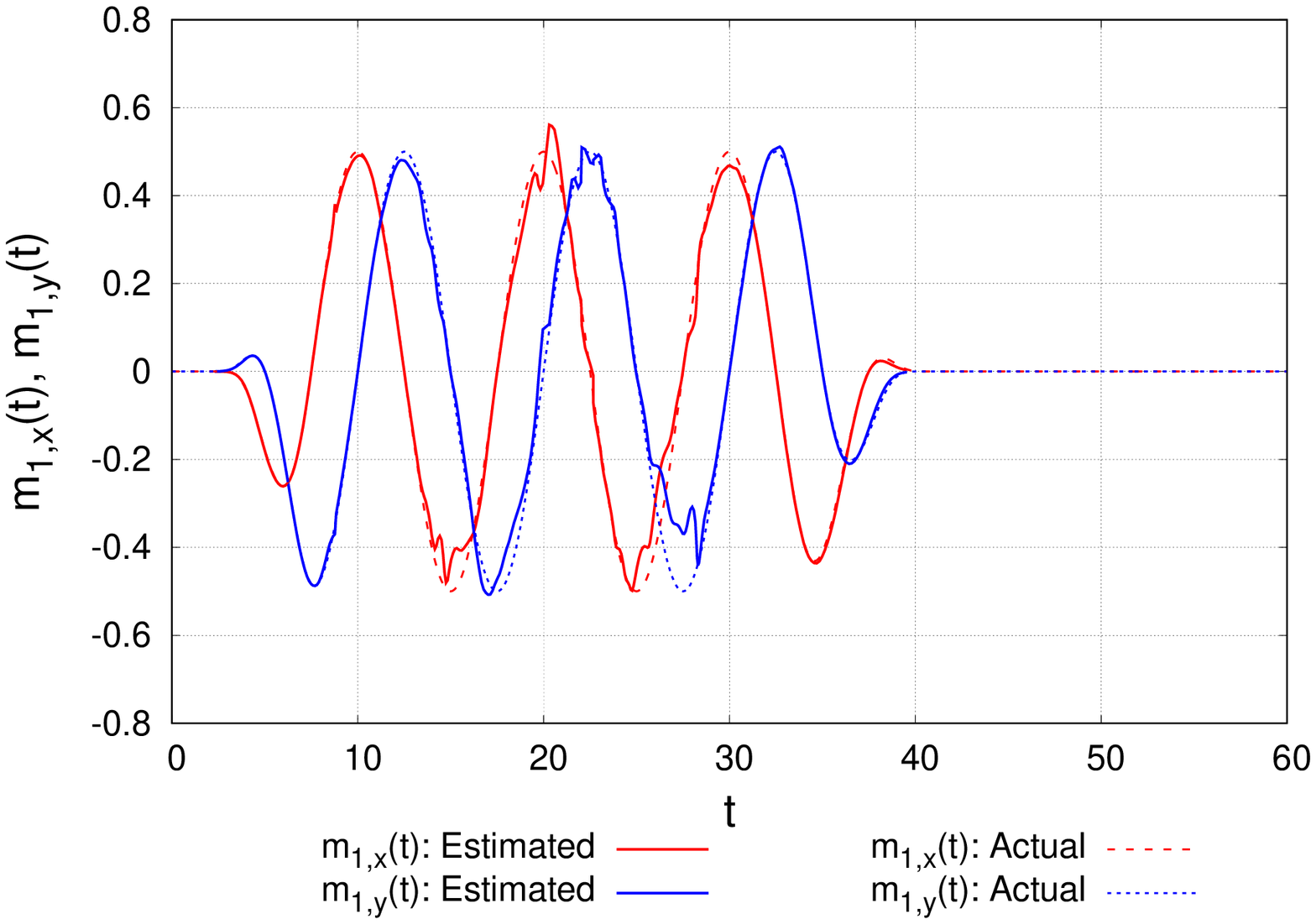}}}
\\[-2ex]
  \subfloat[location of source 2]{
  \resizebox*{8cm}{!}{\includegraphics{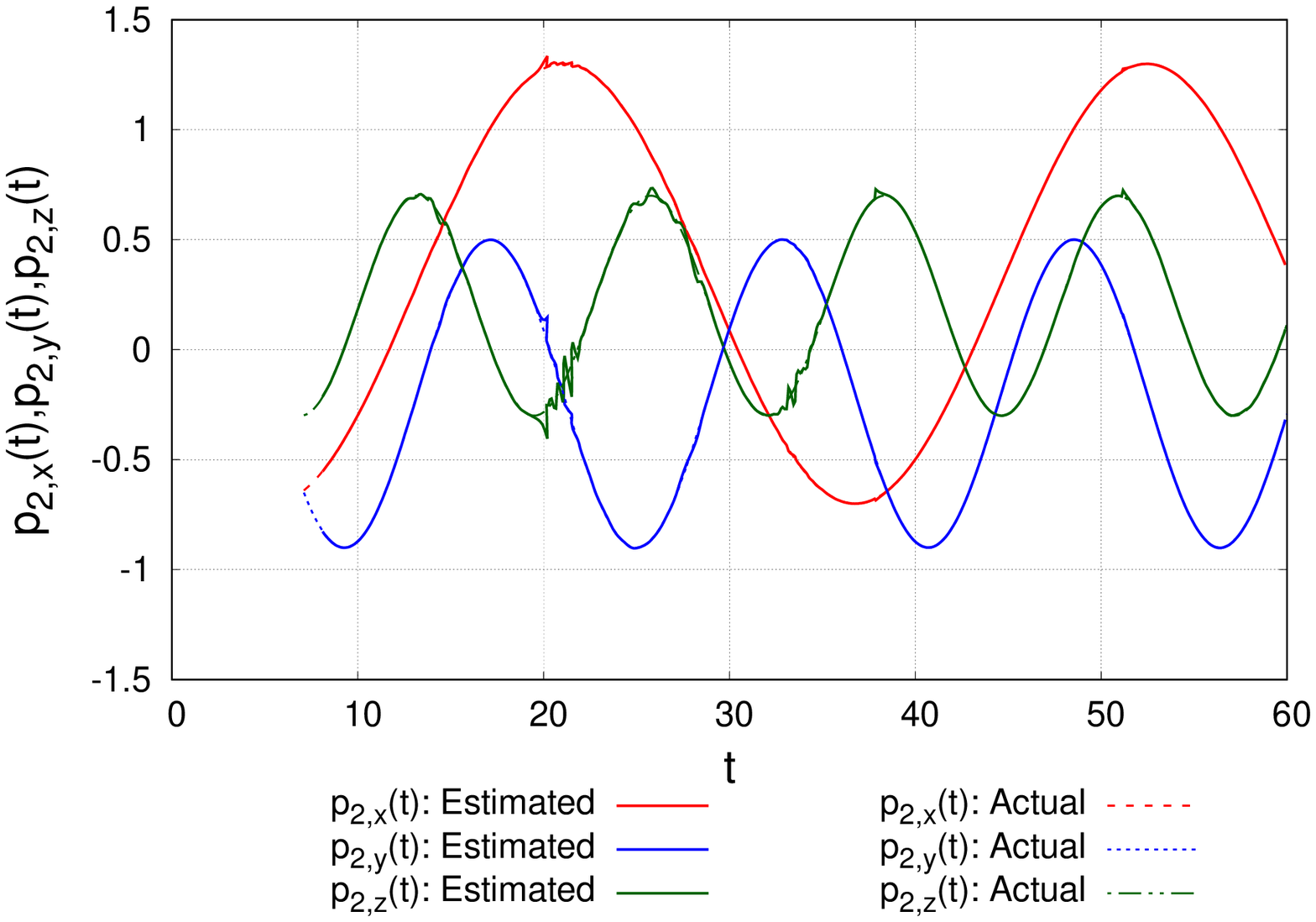}}}
  \subfloat[moment of source 2]{
  \resizebox*{8cm}{!}{\includegraphics{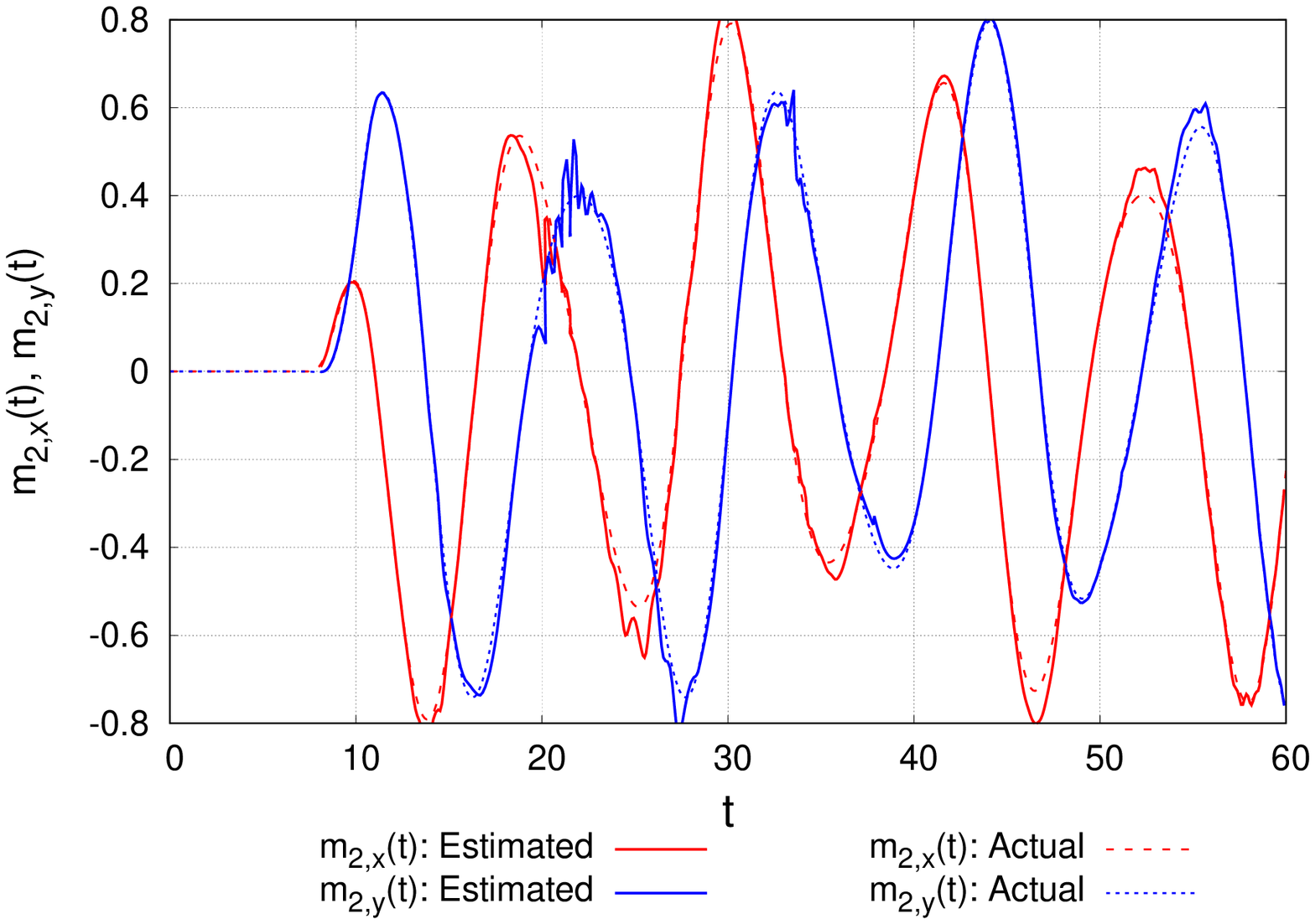}}}
\\[-2ex]
  \subfloat[location of source 3]{
  \resizebox*{8cm}{!}{\includegraphics{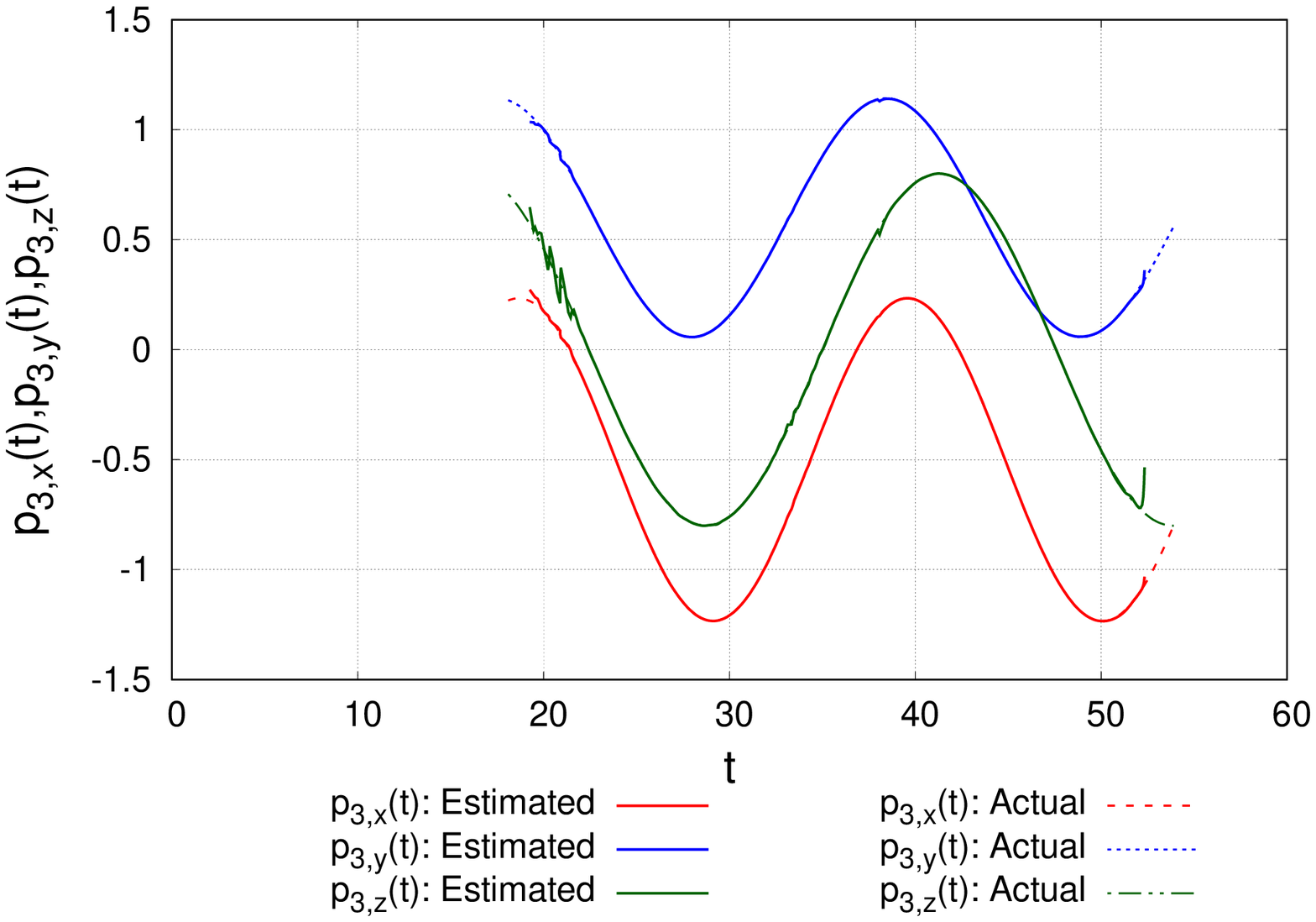}}}
  \subfloat[moment of source 3]{
  \resizebox*{8cm}{!}{\includegraphics{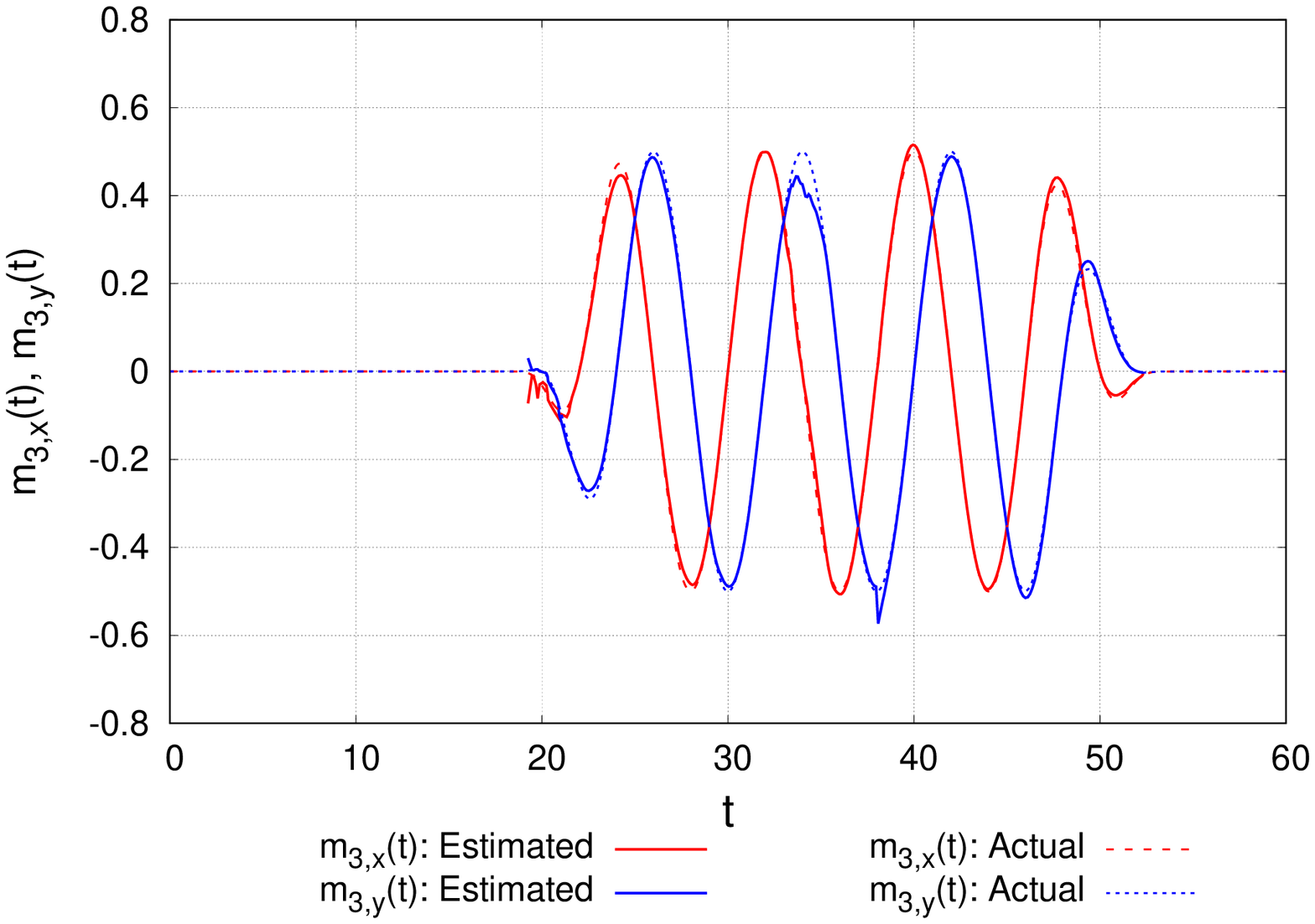}}}
 \caption{Estimated locations and moments of dipole sources for observation data without noise.}
\label{fig:reconstruction_results_loc_nl_0.0_dipole_source}
\end{figure}

\begin{figure}[ht]
\centering
  \subfloat[location of source 1]{
  \resizebox*{8cm}{!}{\includegraphics{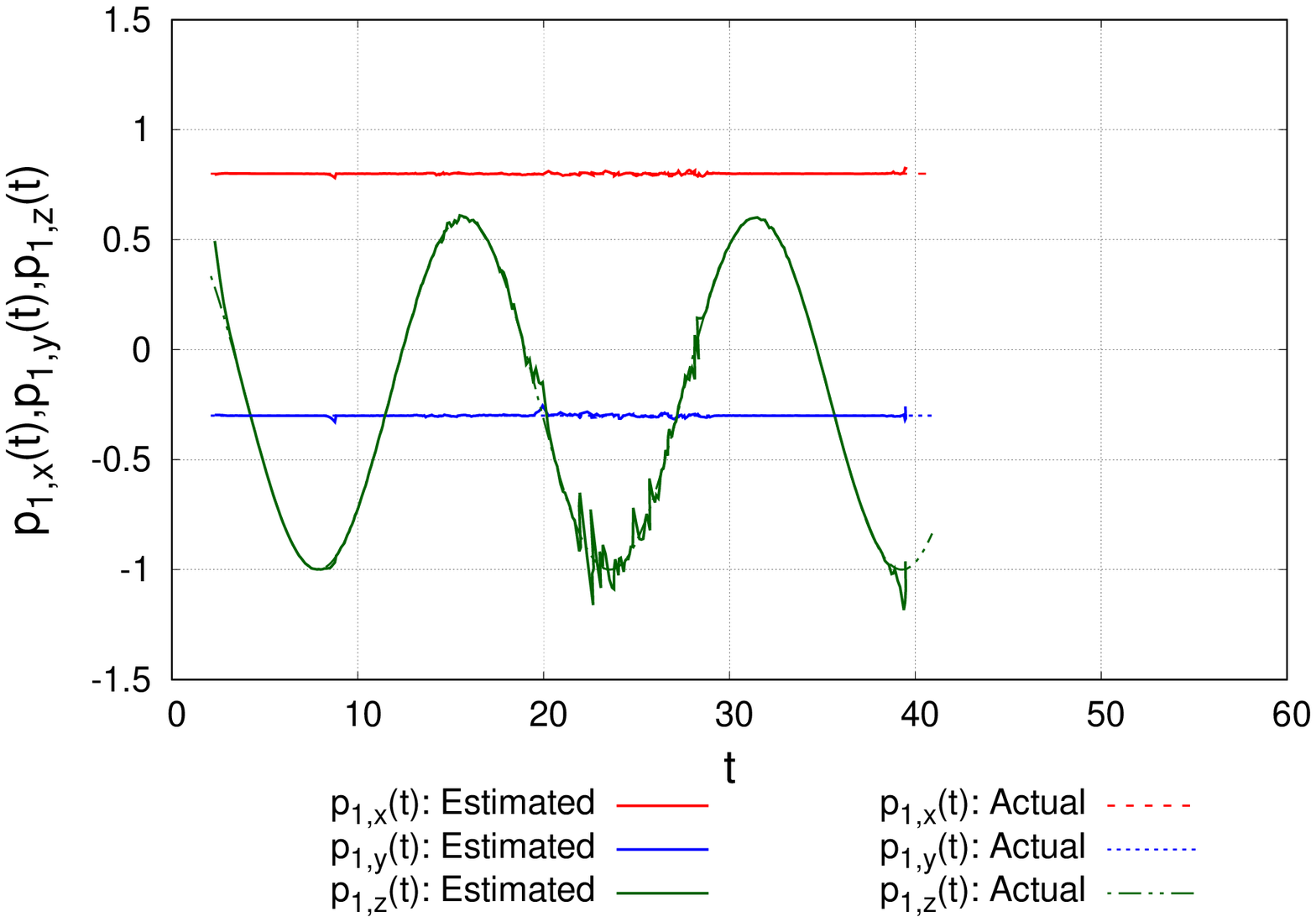}}}
  \subfloat[moment of sources 1]{
  \resizebox*{8cm}{!}{\includegraphics{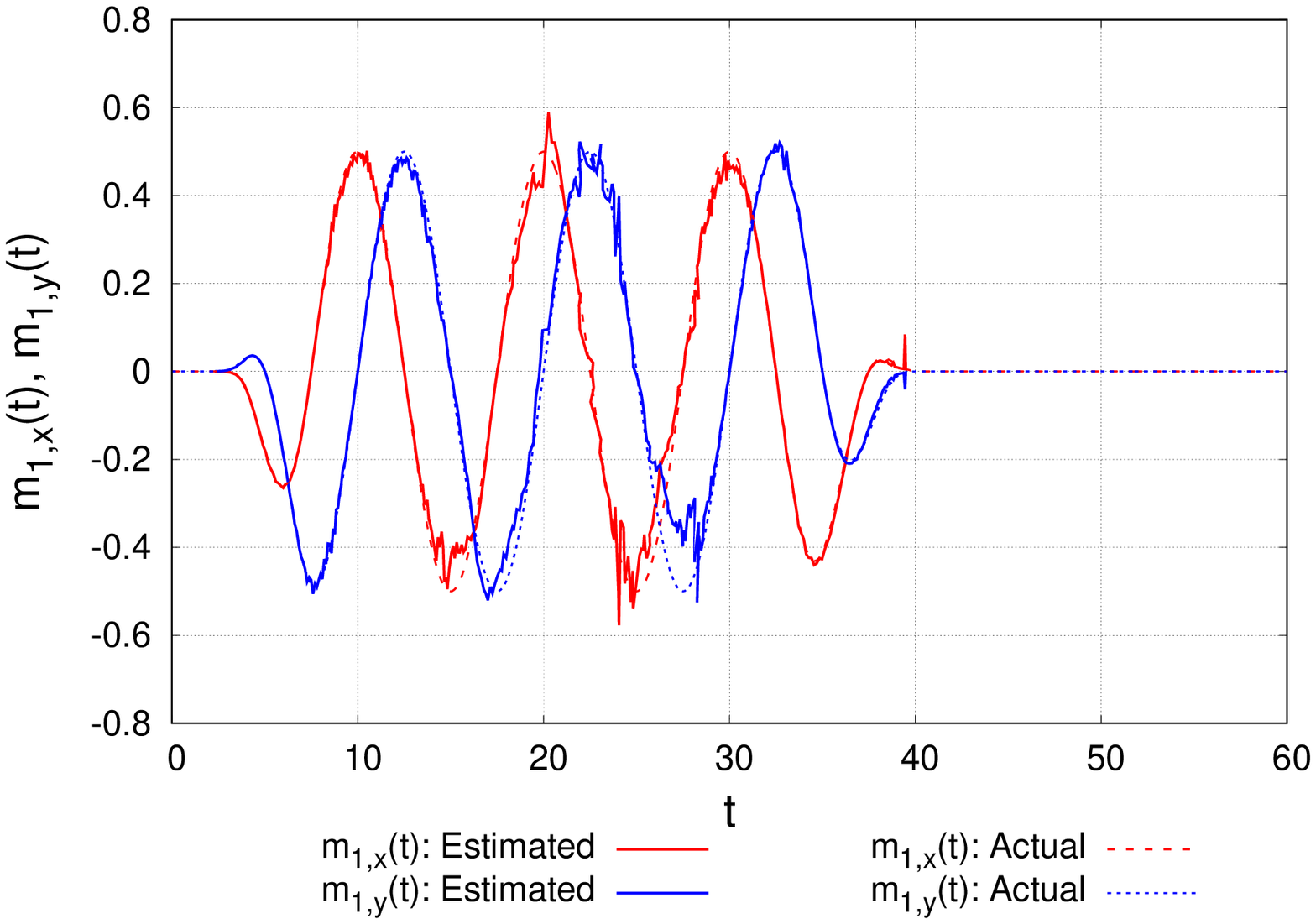}}}
\\[-2ex]
  \subfloat[location of source 2]{
  \resizebox*{8cm}{!}{\includegraphics{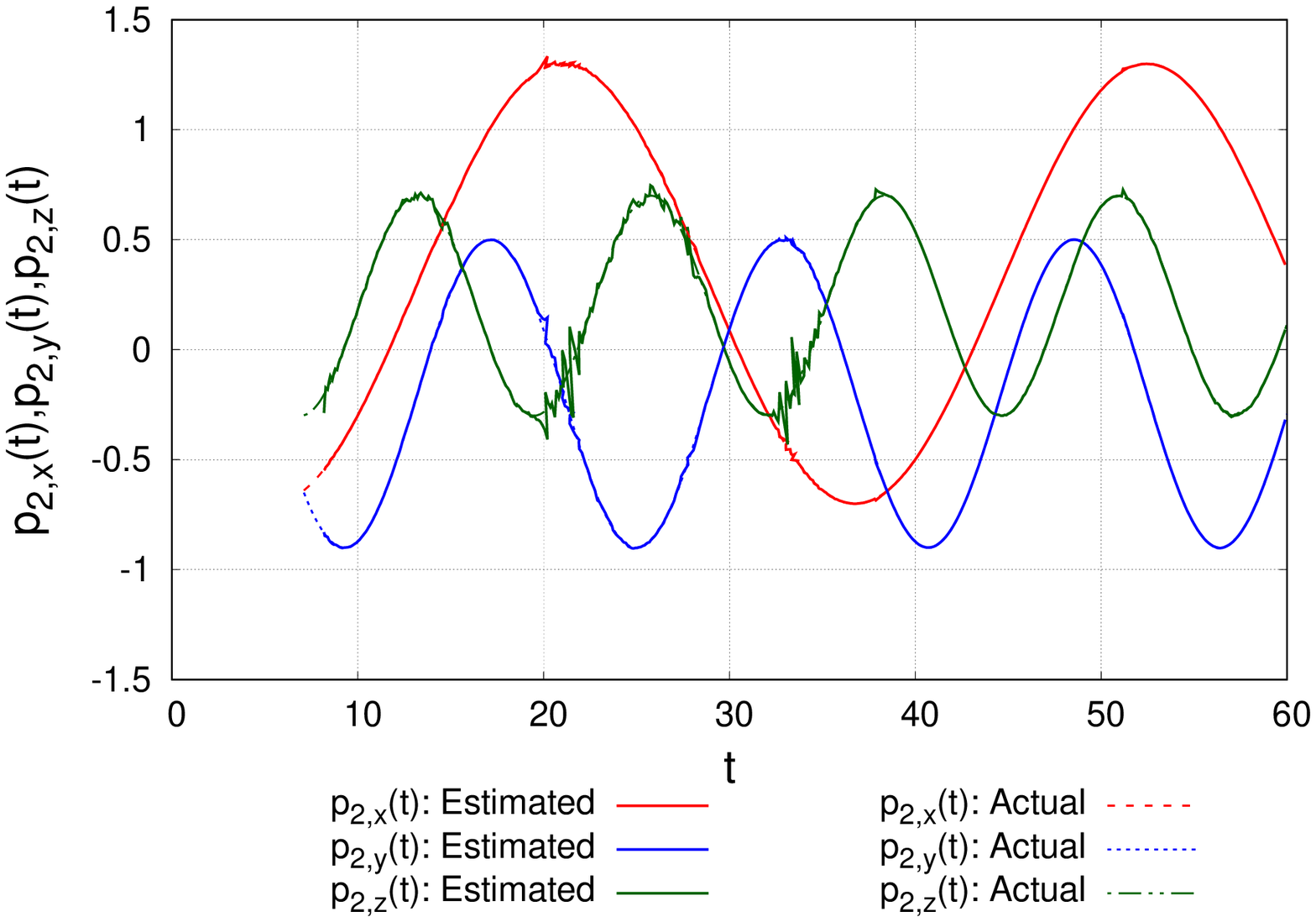}}}
  \subfloat[moment of source 2]{
  \resizebox*{8cm}{!}{\includegraphics{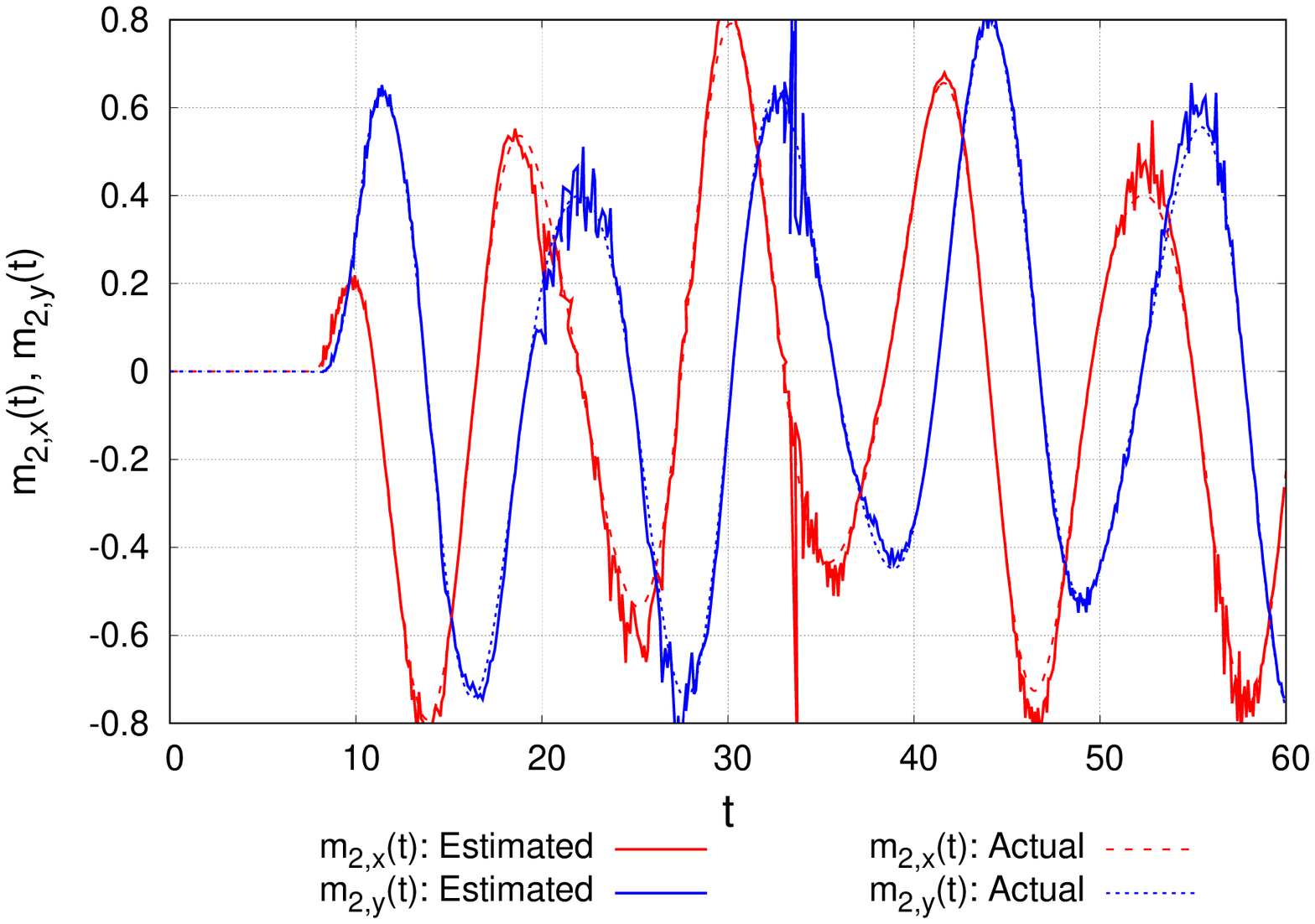}}}
\\[-2ex]
  \subfloat[location of source 3]{
  \resizebox*{8cm}{!}{\includegraphics{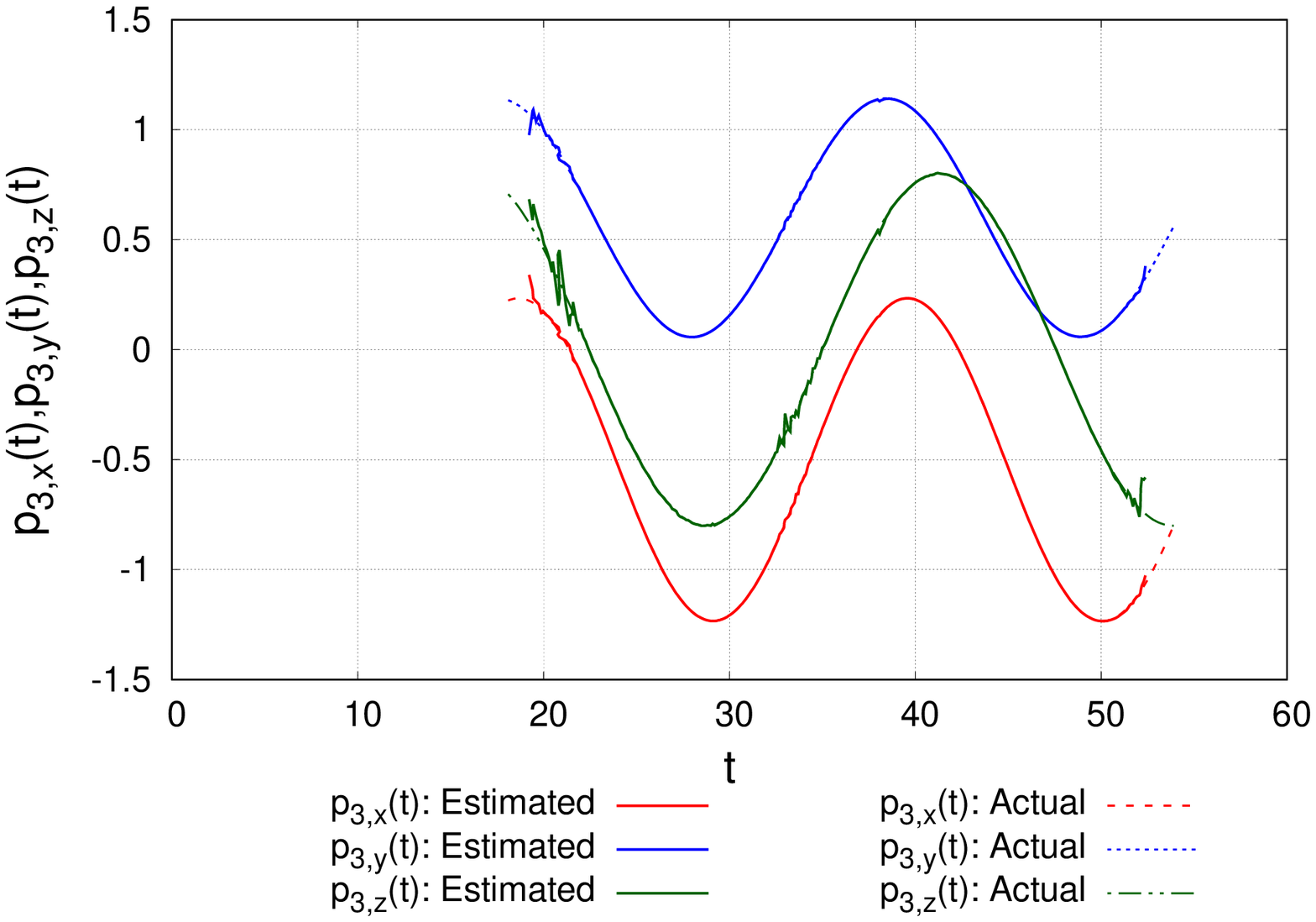}}}
  \subfloat[moment of source 3]{
  \resizebox*{8cm}{!}{\includegraphics{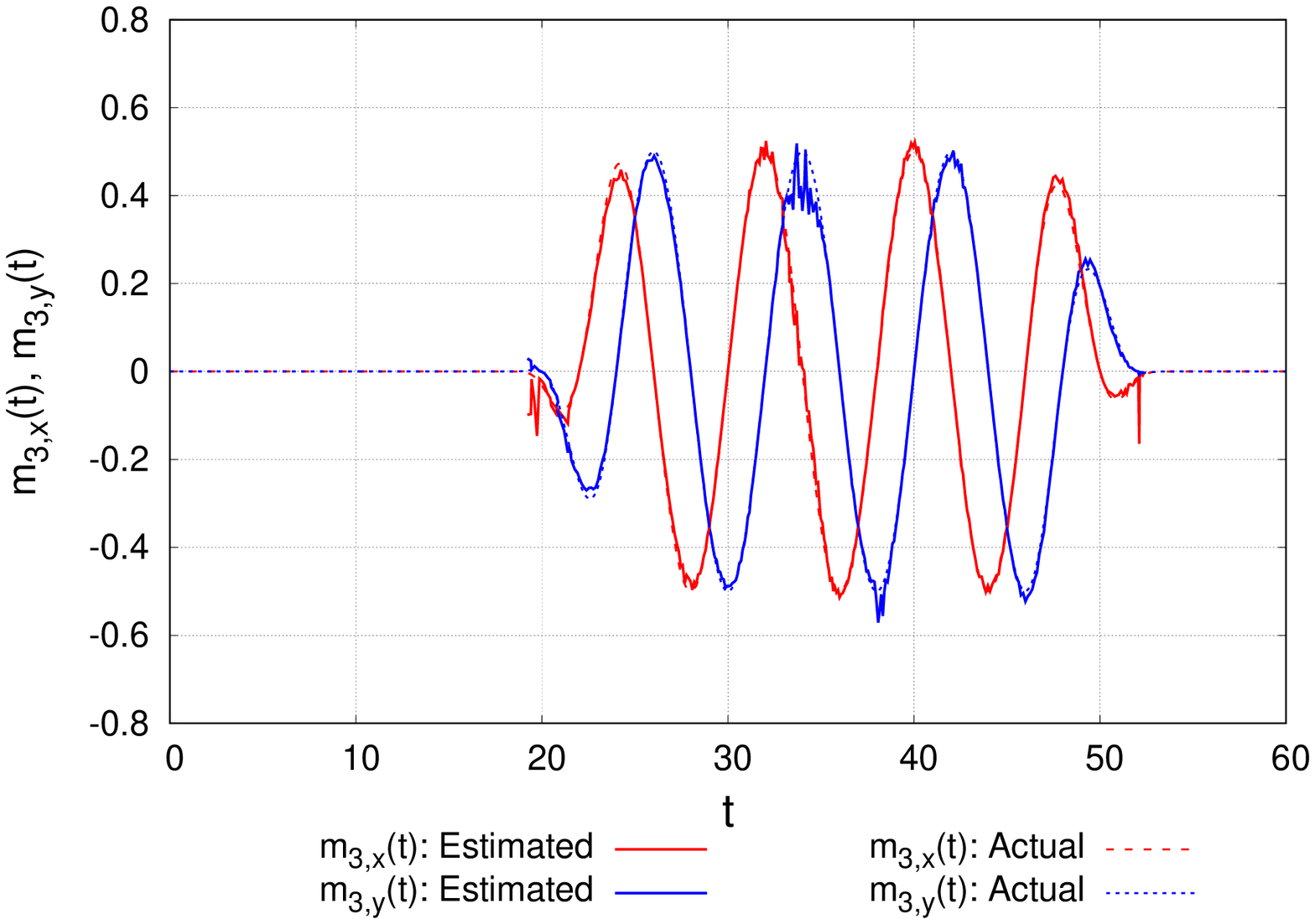}}}
 \caption{Estimated locations and moments of dipole sources for observation data with 0.1\% noise.}
\label{fig:reconstruction_results_loc_nl_0.1_dipole_source}
\end{figure}

\begin{figure}[ht]
\centering
  \subfloat[location of source 1]{
  \resizebox*{8cm}{!}{\includegraphics{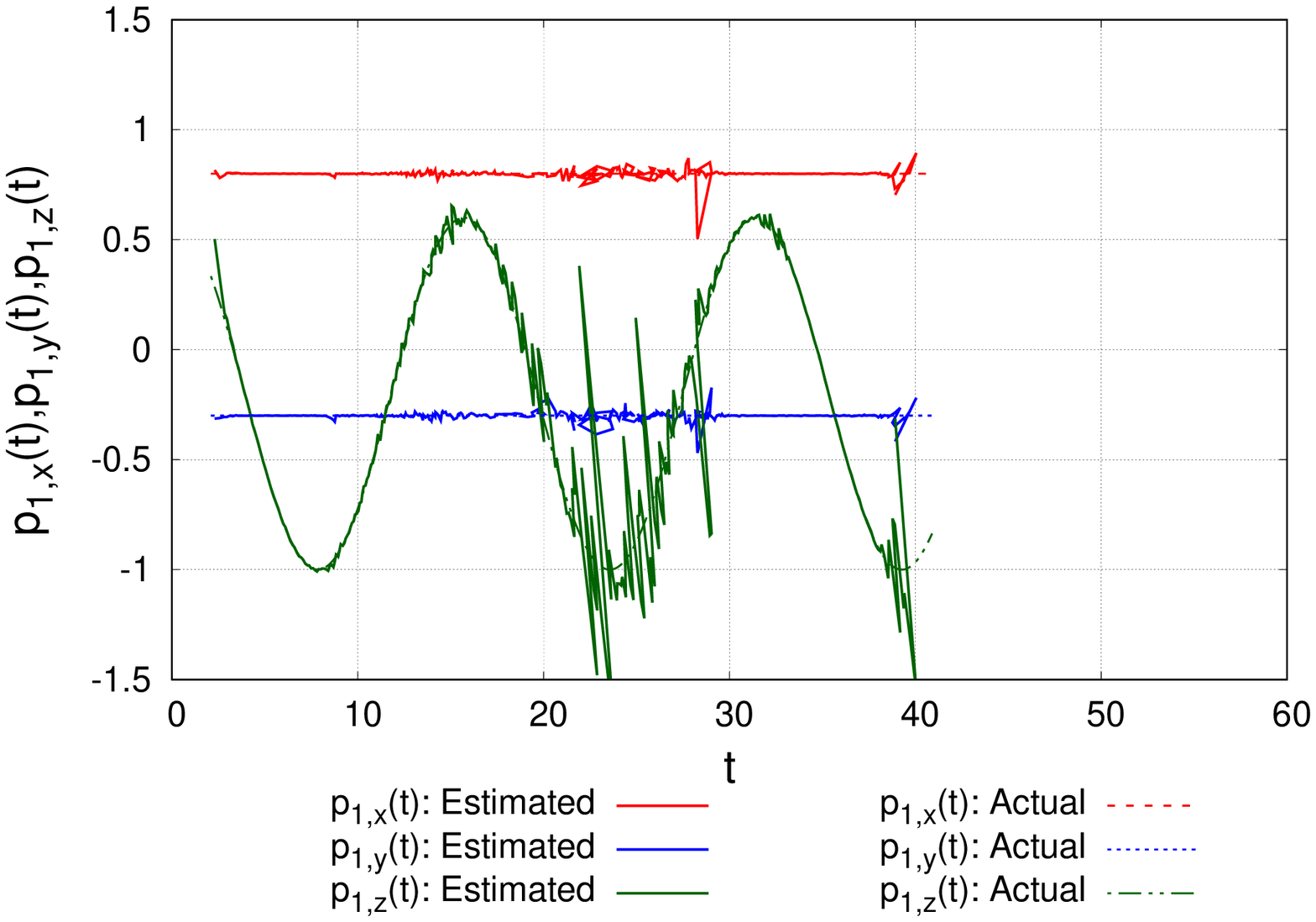}}}
  \subfloat[moment of sources 1]{
  \resizebox*{8cm}{!}{\includegraphics{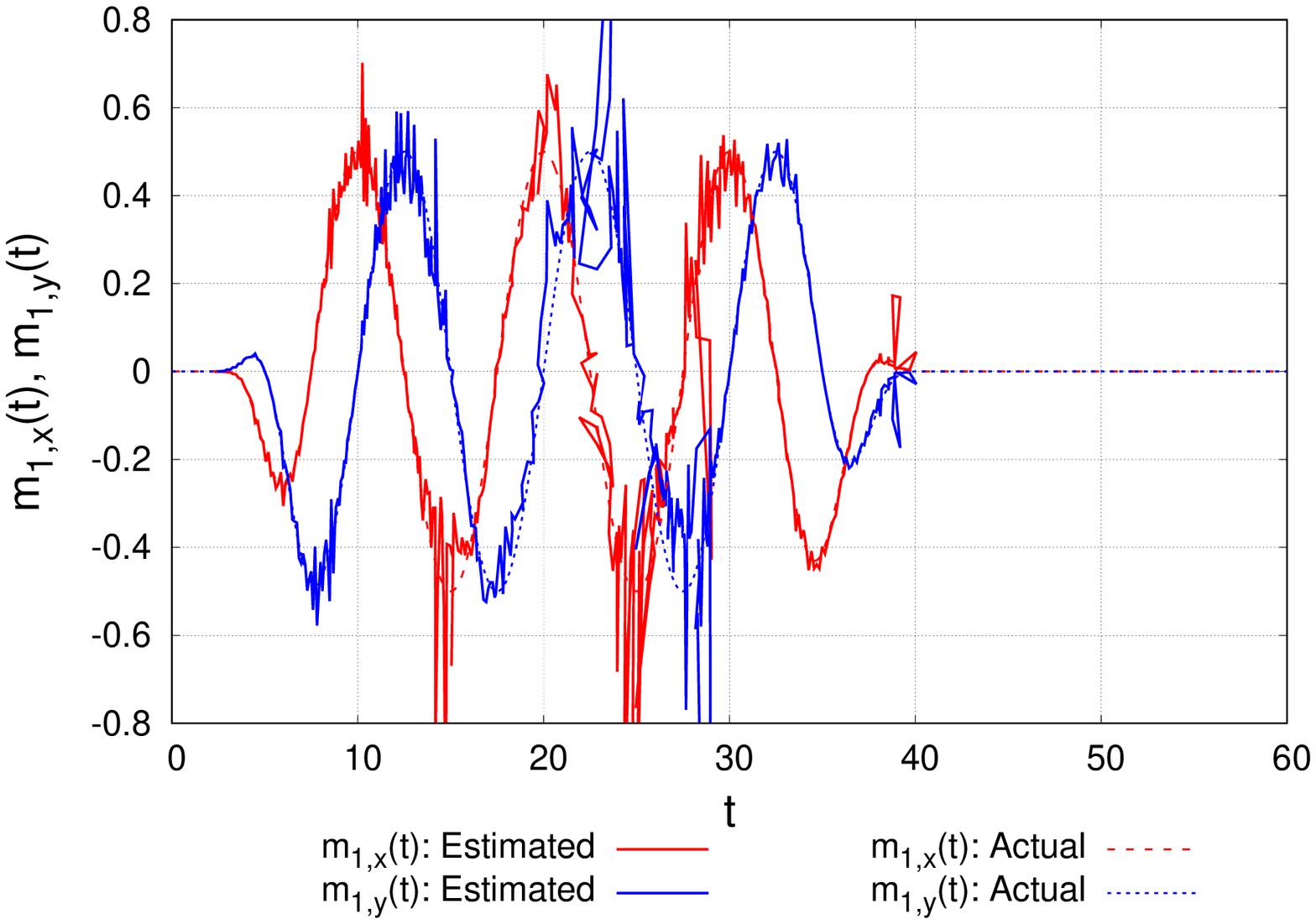}}}
\\[-2ex]
  \subfloat[location of source 2]{
  \resizebox*{8cm}{!}{\includegraphics{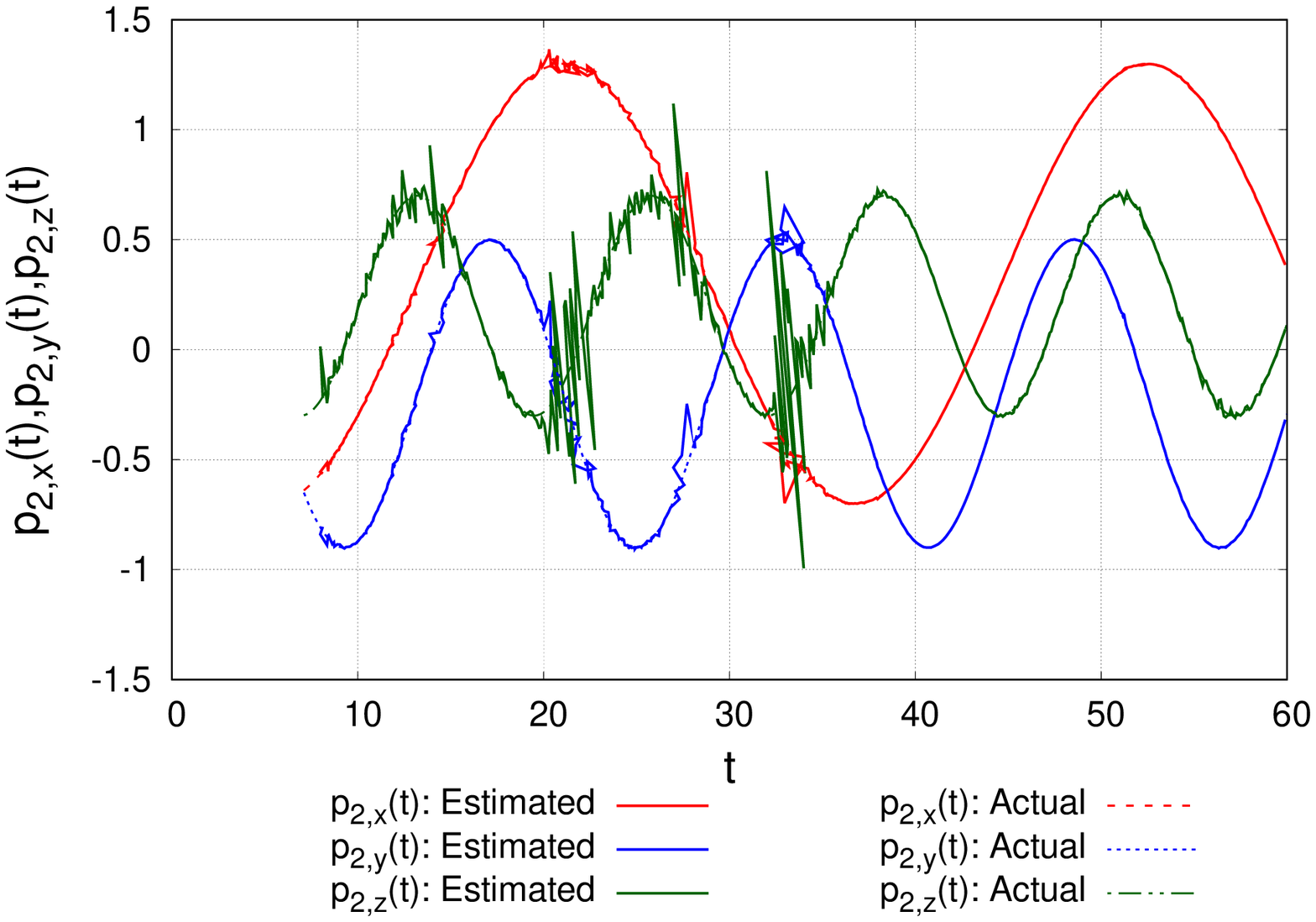}}}
  \subfloat[moment of source 2]{
  \resizebox*{8cm}{!}{\includegraphics{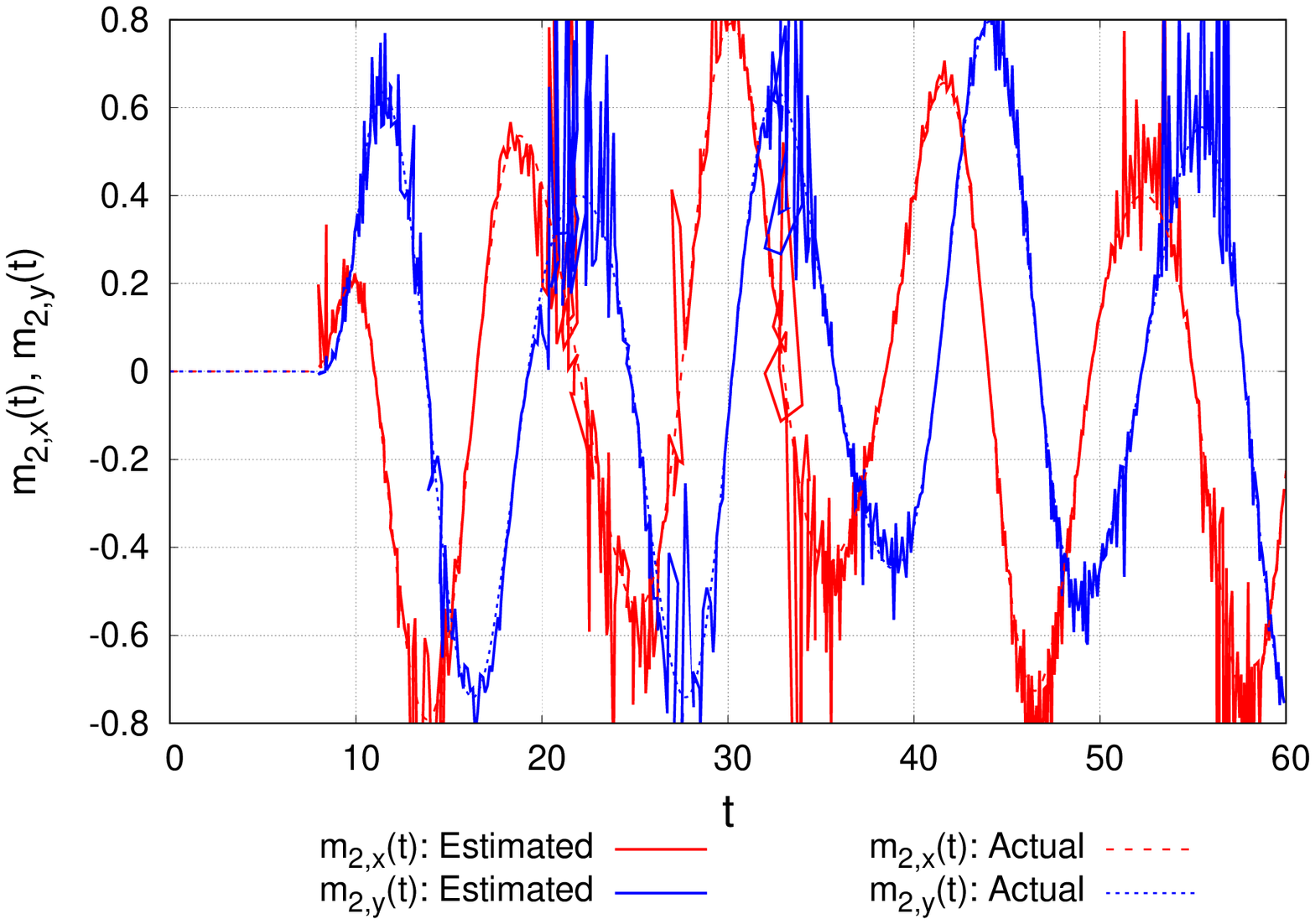}}}
\\[-2ex]
  \subfloat[location of source 3]{
  \resizebox*{8cm}{!}{\includegraphics{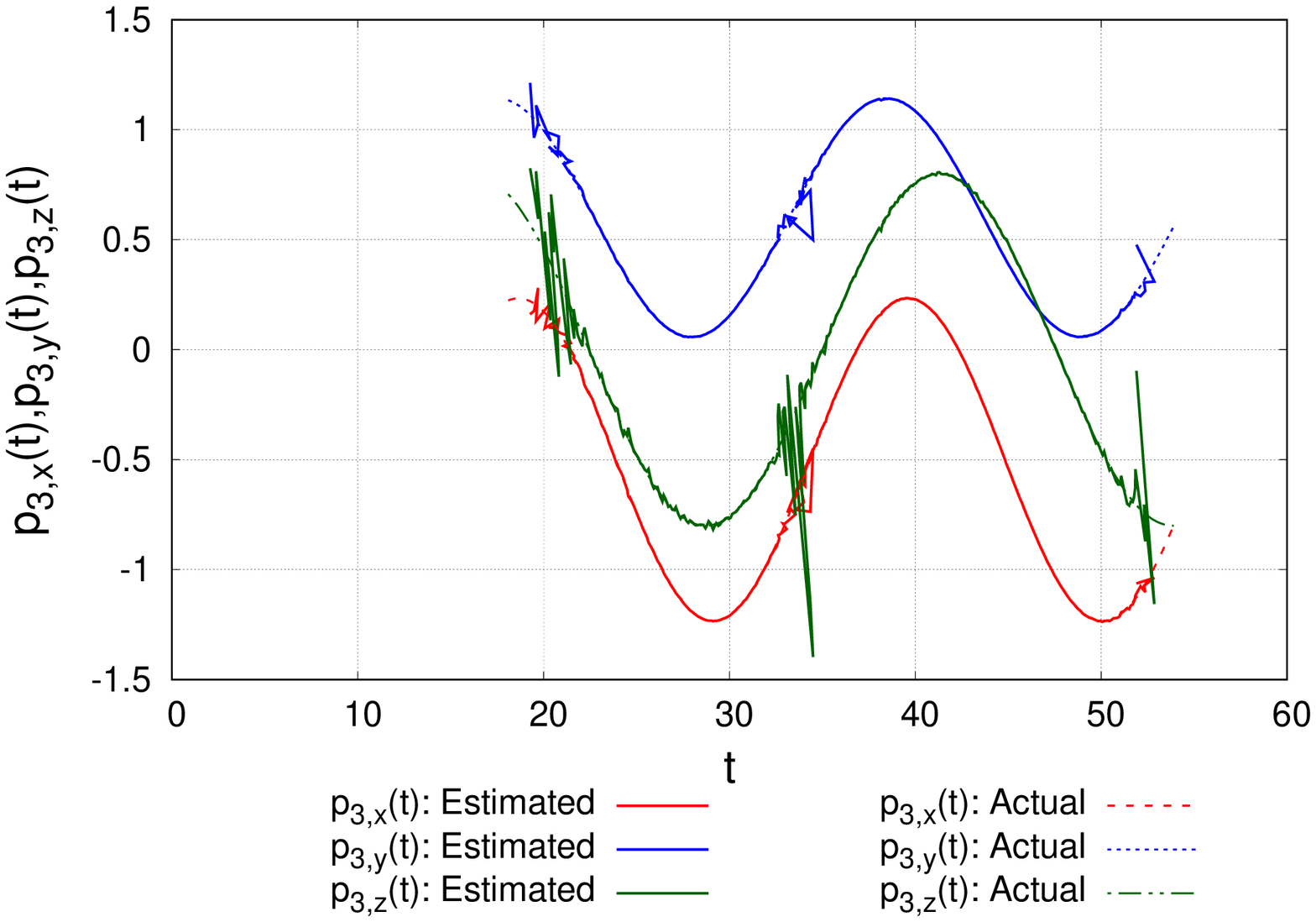}}}
  \subfloat[moment of source 3]{
  \resizebox*{8cm}{!}{\includegraphics{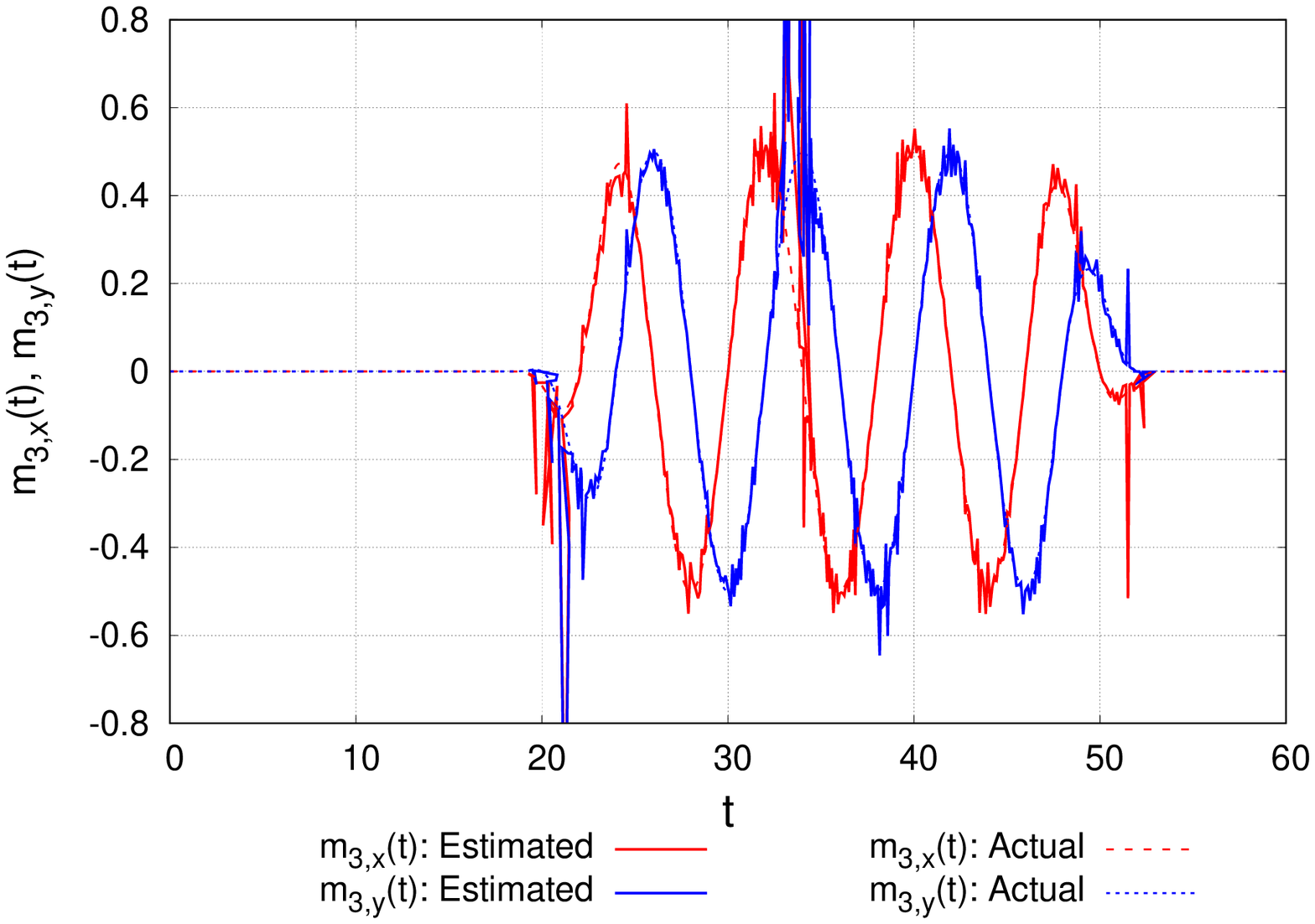}}}
 \caption{Estimated locations and moments of dipole sources for observation data with 0.5\% noise.}
\label{fig:reconstruction_results_loc_nl_0.5_dipole_source}
\end{figure}

\begin{figure}[ht]
\centering
  \subfloat[location of source 1]{
  \resizebox*{8cm}{!}{\includegraphics{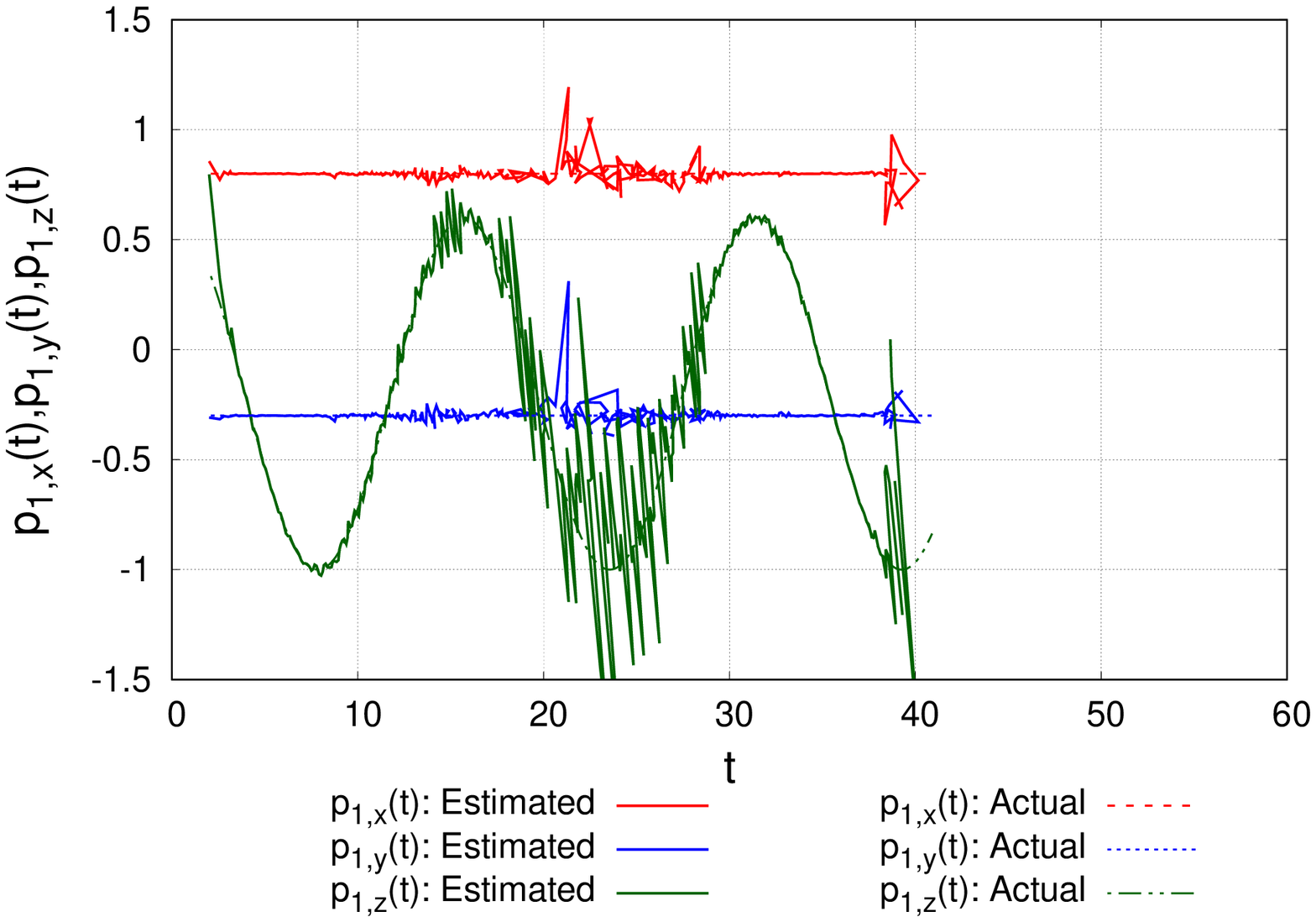}}}
  \subfloat[moment of sources 1]{
  \resizebox*{8cm}{!}{\includegraphics{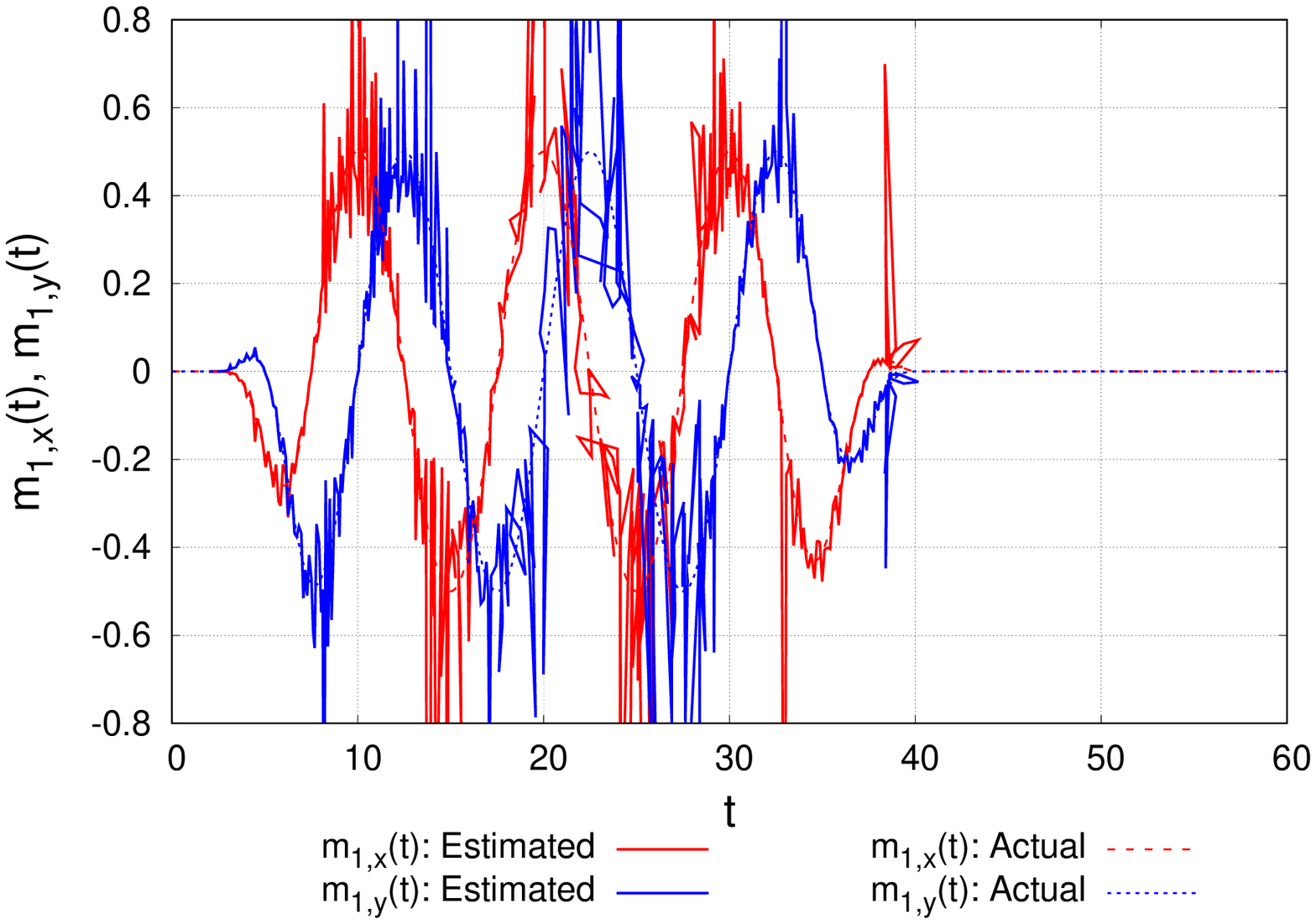}}}
  \\[-2ex]
  \subfloat[location of source 2]{
  \resizebox*{8cm}{!}{\includegraphics{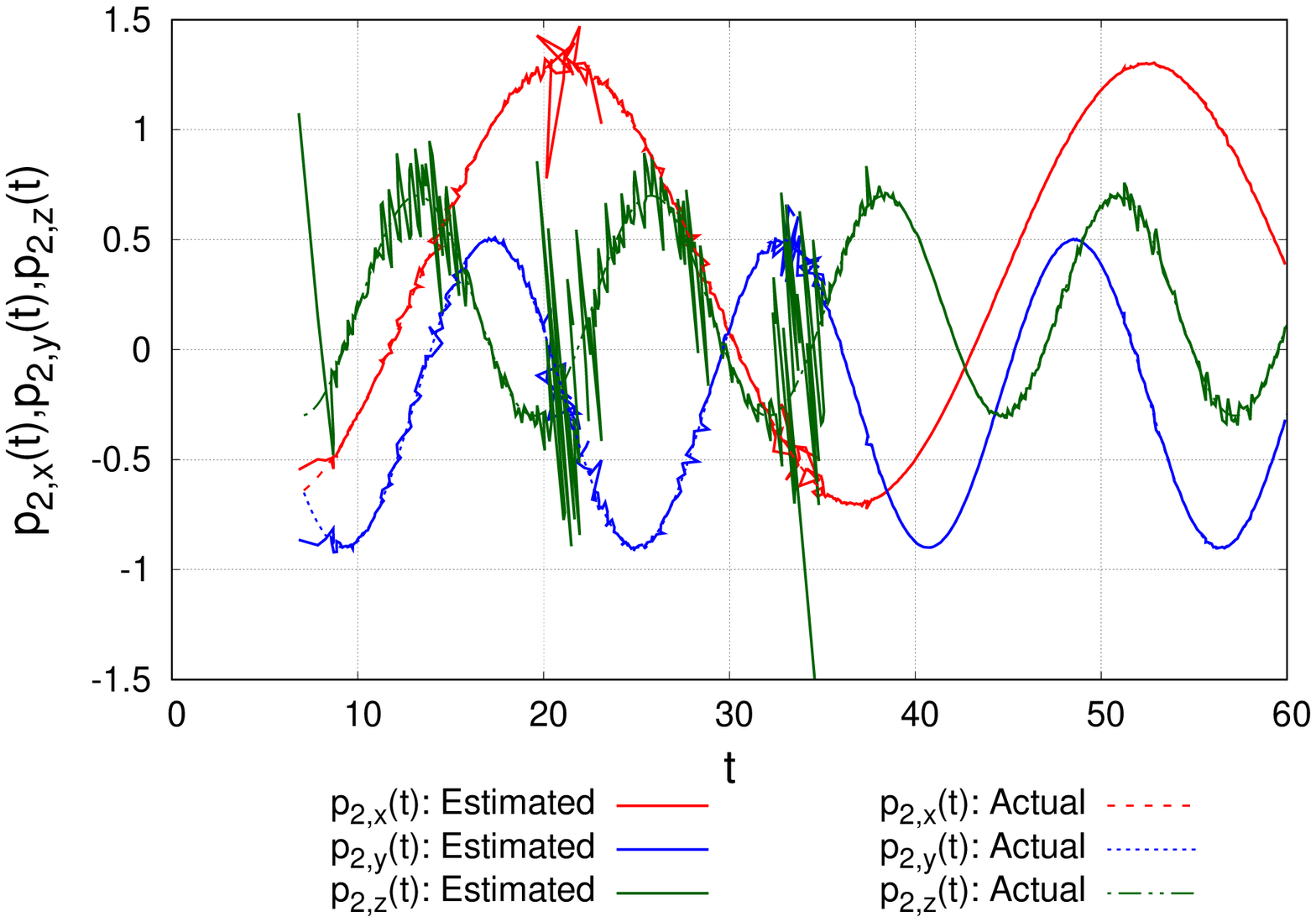}}}
  \subfloat[moment of source 2]{
  \resizebox*{8cm}{!}{\includegraphics{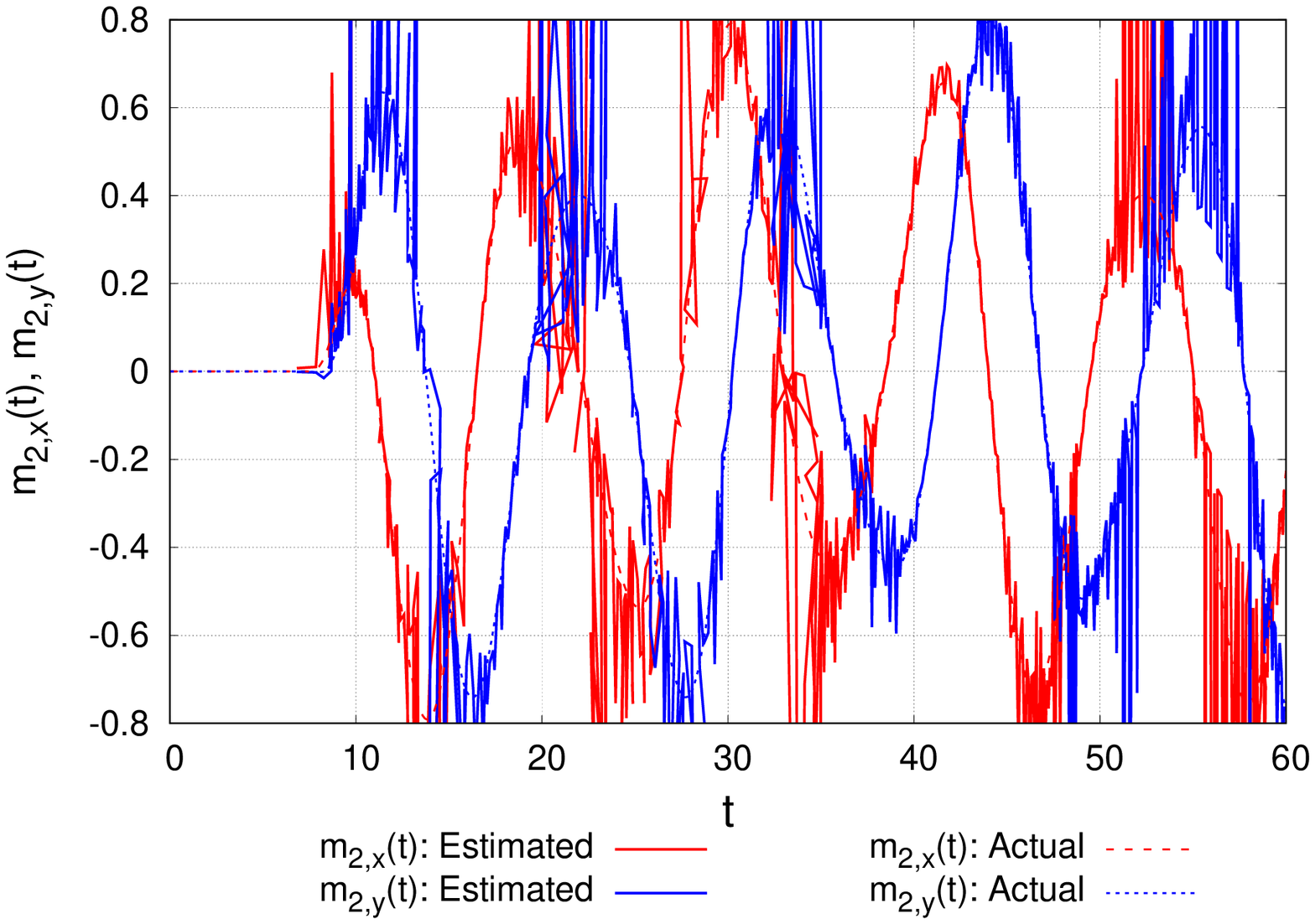}}}
  \\[-2ex]
  \subfloat[location of source 3]{
  \resizebox*{8cm}{!}{\includegraphics{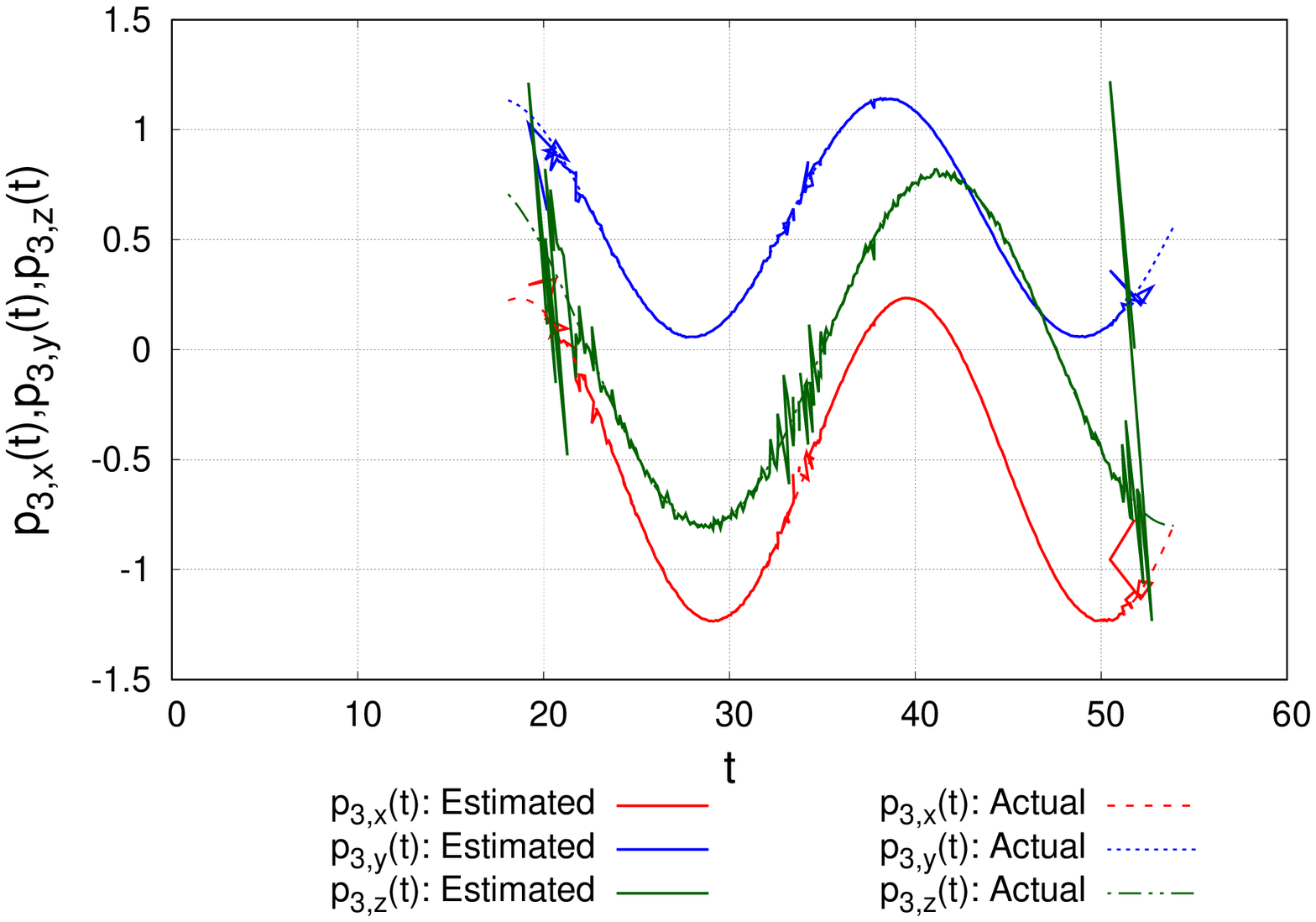}}}
  \subfloat[moment of source 3]{
  \resizebox*{8cm}{!}{\includegraphics{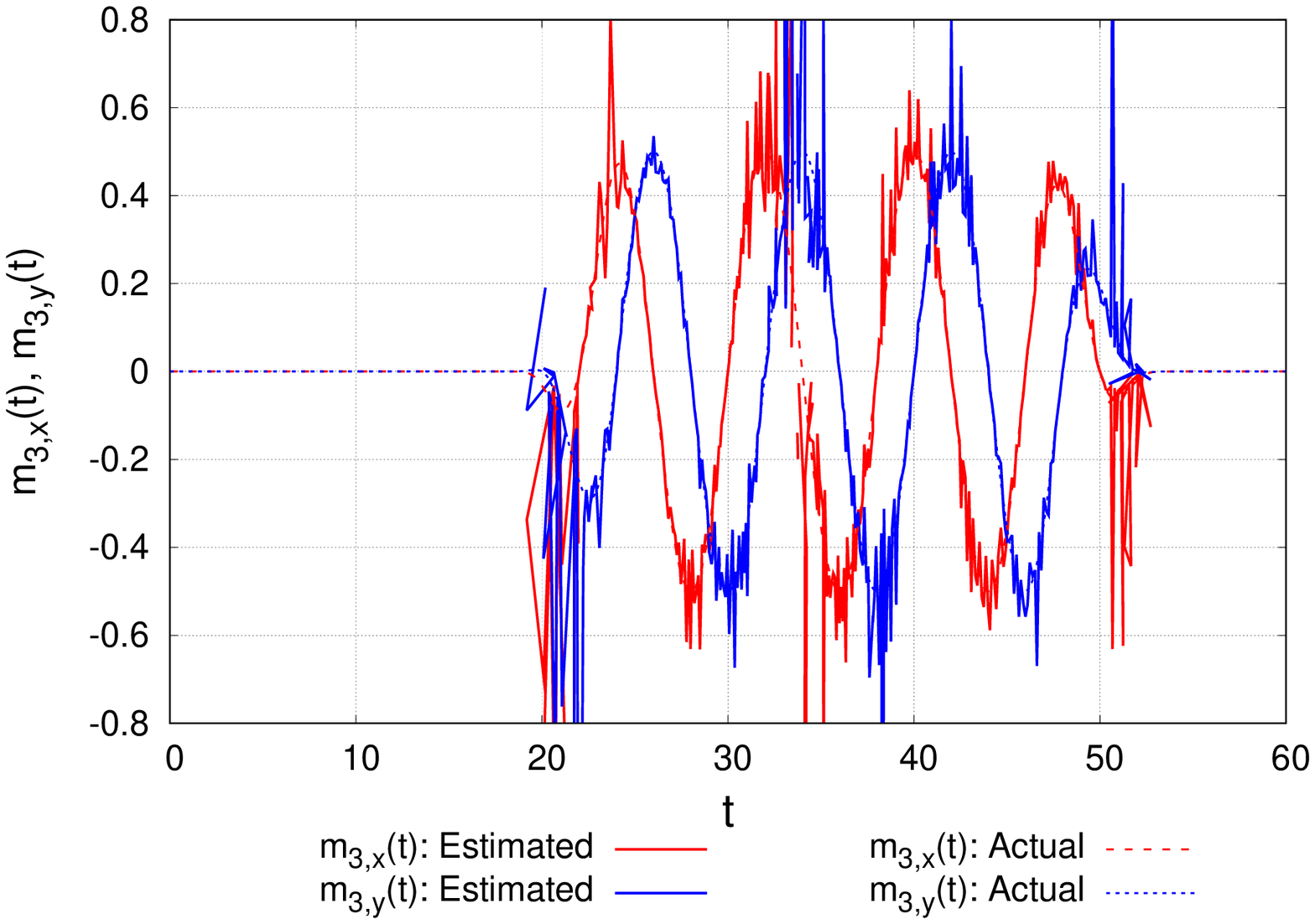}}}
 \caption{Estimated locations and moments of dipole sources for observation data with 1.0\% noise.}
\label{fig:reconstruction_results_loc_nl_1.0_dipole_source}
\end{figure}

Next, we show the reconstruction results for locations and moments of
 dipole sources.
Figures \ref{fig:reconstruction_results_loc_nl_0.0_dipole_source}-\ref{fig:reconstruction_results_loc_nl_1.0_dipole_source}
 display the reconstruction results 
 from observations without noise, and with 0.1\%, 0.5\% and 1.0\% noise.
Also in  
Table \ref{table:Average_Error_dipole_source}, we show the average errors of estimated
 locations and moments.
Here we define the average error of each moment by
\begin{displaymath}
\mbox{average error of moment } \boldsymbol{m}_k \equiv
 \displaystyle
\sqrt{\frac{1}{|I|}\int_{I}|\widehat{\boldsymbol{m}}_k(t_k(\tau)) -
 \boldsymbol{m}_k(t_k(\tau))|^2 {\mathrm d}\tau},
\end{displaymath}
where $\widehat{\boldsymbol{m}}_k$ denotes the estimated moment of $k$-th
dipole source.
Same as for the reconstruction of point sources, 
our method works well if the noise of observation data is smaller than 0.5\%.
The noise of observation becomes 1\%, the influence of noise becomes
unignorable, and our method can not give reliable estimates under 5\%
noise.
We can see that observation noises affect to the estimations of
 moments more heavily than the reconstruction of locations.
The reconstruction result also becomes bad when the
norm of moment is small, and when the arrangement of
dipoles is clustered.
\par
%
%

\begin{table}[p]
\tbl{The average errors of estimated locations and moments in each interval}
{\begin{tabular}{cccccc}
\toprule
 & $2.4 \leq \tau < 6.8$ & $6.8 \leq \tau < 17.3$ 
 & $17.3 \leq \tau < 40.2$ & $40.2 \leq \tau < 54.8$ & $54.8 \leq \tau < 60.0$ \\
 & $(K(\tau)=1)$ & $(K(\tau)=2)$ & $(K(\tau)=3)$ & $(K(\tau)=2)$ &
		     $(K(\tau)=1)$ \\
\midrule
(a) & \multicolumn{5}{c}{Without noise} \\
$\boldsymbol{p}_1$ & $5.6E-2$ & $6.9E-3$ & $3.7E-2$ & $-$      & $-$ \\
$\boldsymbol{p}_2$ & $-$      & $8.4E-3$ & $2.5E-2$ & $4.3E-3$ & $3.2E-3$ \\
$\boldsymbol{p}_3$ & $-$      & $-$      & $2.0E-2$ & $2.6E-2$ & $-$ \\
$\boldsymbol{m}_1$ & $1.6E-3$ & $3.6E-2$ & $6.7E-2$ & $-$      & $-$ \\
$\boldsymbol{m}_2$ & $-$      & $2.2E-2$ & $5.9E-2$ & $4.0E-2$ & $2.7E-2$ \\
$\boldsymbol{m}_3$ & $-$      & $-$      & $2.8E-2$ & $1.3E-2$ & $-$ \\ 
\midrule
(b) & \multicolumn{5}{c}{With $0.1\%$ noise} \\
$\boldsymbol{p}_1$ & $5.5E-2$ & $8.2E-3$ & $5.4E-2$ & $-$      & $-$ \\
$\boldsymbol{p}_2$ & $-$      & $1.5E-2$ & $4.3E-2$ & $4.6E-3$ & $4.1E-3$ \\
$\boldsymbol{p}_3$ & $-$      & $-$      & $3.3E-2$ & $2.9E-2$       & $-$ \\
$\boldsymbol{m}_1$ & $5.1E-3$ & $3.9E-2$ & $7.1E-2$ & $-$      & $-$ \\
$\boldsymbol{m}_2$ & $-$      & $3.2E-2$ & $1.5E-1$ & $4.7E-2$ & $5.6E-2$ \\
$\boldsymbol{m}_3$ & $-$      & $-$      & $3.4E-2$ & $2.1E-2$       & $-$ \\ 
\midrule
(c) & \multicolumn{5}{c}{With $0.5\%$ noise} \\
$\boldsymbol{p}_1$ & $5.5E-2$ & $2.6E-2$ & $2.1E-1$ & $-$      & $-$ \\
$\boldsymbol{p}_2$ & $-$      & $6.2E-2$ & $1.7E-1$ & $8.6E-3$ & $1.5E-2$ \\
$\boldsymbol{p}_3$ & $-$      & $-$      & $1.4E-1$ & $7.7E-2$ & $-$ \\
$\boldsymbol{m}_1$ & $2.0E-2$ & $1.1E-1$ & $3.4E-1$ & $-$      & $-$ \\
$\boldsymbol{m}_2$ & $-$      & $1.6E-1$ & $6.9E-1$ & $1.4E-1$ & $5.6E-1$ \\
$\boldsymbol{m}_3$ & $-$      & $-$      & $4.4E-1$ & $6.2E-2$  & $-$ \\ 
\midrule
(d) & \multicolumn{5}{c}{With $1.0\%$ noise} \\
$\boldsymbol{p}_1$ & $1.0E-1$ & $5.3E-2$ & $2.7E-1$ & $-$       & $-$ \\
$\boldsymbol{p}_2$ & $-$      & $1.9E-1$ & $2.5E-1$ & $2.2E-2$ & $3.2E-2$ \\
$\boldsymbol{p}_3$ & $-$      & $-$      & $1.7E-1$ & $2.1E-1$       & $-$ \\
$\boldsymbol{m}_1$ & $4.0E-2$ & $6.1E-1$ & $9.2E-1$ & $-$ & $-$ \\
$\boldsymbol{m}_2$ & $-$      & $1.2E+0$ & $9.0E-1$ & $8.4E-1$ & $2.2E+0$ \\
$\boldsymbol{m}_3$ & $-$      & $-$      & $5.9E-1$ & $1.6E-1$       & $-$ \\ 
\midrule
(e) & \multicolumn{5}{c}{With $5.0\%$ noise} \\
$\boldsymbol{p}_1$ & $9.5E-2$ & $2.7E-1$ & $6.8E-1$ & $-$ & $-$ \\
$\boldsymbol{p}_2$ & $-$      & $4.5E-1$ & $4.4E-1$ & $1.1E-1$ & $1.3E-1$ \\
$\boldsymbol{p}_3$ & $-$      & $-$      & $3.5E-1$ & $3.3E-1$       & $-$ \\
$\boldsymbol{m}_1$ & $1.2E+0$ & $2.2E+0$ & $2.9E+0$ & $-$ & $-$ \\
$\boldsymbol{m}_2$ & $-$      & $3.2E+0$ & $2.9E+0$ & $2.5E+0$ & $2.7E+0$ \\
$\boldsymbol{m}_3$ & $-$      & $-$      & $2.3E+0$ & $1.7E+0$       & $-$ \\ 
\bottomrule
\end{tabular}}
\label{table:Average_Error_dipole_source}
\end{table}
%
%
%
%
%
%
%
%
%
%
%
%
%
\par
Finally, we show the average of error of estimated moments by using Step
4' instead of Step 4 in Table \ref{table:Average_Error_dipole_source_another_idea}.
As same as for moving point source, the errors of estimation results are
almost the same in both methods if observation noise is smaller than
1.0\%.
\par
%
\begin{table}[p]
\tbl{The average errors of estimated moments in  each interval 
using Step 4' instead of Step 4.}
{\begin{tabular}{cccccc}
\toprule
 & $2.4 \leq \tau < 6.8$ & $6.8 \leq \tau < 17.3$ 
 & $17.3 \leq \tau < 40.2$ & $40.2 \leq \tau < 54.8$ & $54.8 \leq \tau < 60.0$ \\
 & $(K(\tau)=1)$ & $(K(\tau)=2)$ & $(K(\tau)=3)$ & $(K(\tau)=2)$ &
		     $(K(\tau)=1)$ \\
\midrule
(a) & \multicolumn{5}{c}{Without noise} \\
$\boldsymbol{m}_1$ & $1.7E-2$ & $2.4E-2$ & $1.4E-1$ & $-$      & $-$ \\
$\boldsymbol{m}_2$ & $-$      & $6.0E-2$ & $8.8E-2$ & $3.1E-2$ & $2.4E-2$ \\
$\boldsymbol{m}_3$ & $-$      & $-$      & $2.6E-2$ & $1.4E-2$ & $-$ \\ 
\midrule
(b) & \multicolumn{5}{c}{With $0.1\%$ noise} \\
$\boldsymbol{m}_1$ & $1.8E-2$ & $2.8E-2$ & $2.1E-2$ & $-$      & $-$ \\
$\boldsymbol{m}_2$ & $-$      & $8.0E-2$ & $1.6E-1$ & $3.3E-2$ & $3.1E-2$ \\
$\boldsymbol{m}_3$ & $-$      & $-$      & $4.9E-2$ & $1.5E-2$       & $-$ \\ 
\midrule
(c) & \multicolumn{5}{c}{With $0.5\%$ noise} \\
$\boldsymbol{m}_1$ & $2.1E-2$ & $8.8E-2$ & $8.5E-1$ & $-$      & $-$ \\
$\boldsymbol{m}_2$ & $-$      & $4.1E-1$ & $7.2E-1$ & $4.7E-2$ & $9.1E-2$ \\
$\boldsymbol{m}_3$ & $-$      & $-$      & $3.3E-1$ & $2.4E-2$  & $-$ \\ 
\midrule
(d) & \multicolumn{5}{c}{With $1.0\%$ noise} \\
$\boldsymbol{m}_1$ & $3.6E-2$ & $2.0E-1$ & $1.1E+0$ & $-$ & $-$ \\
$\boldsymbol{m}_2$ & $-$      & $8.4E-1$ & $8.4E-1$ & $1.0E-1$ & $2.0E-1$ \\
$\boldsymbol{m}_3$ & $-$      & $-$      & $2.3E-1$ & $5.1E-2$       & $-$ \\ 
\midrule
(e) & \multicolumn{5}{c}{With $5.0\%$ noise} \\
$\boldsymbol{m}_1$ & $1.1E-1$ & $8.7E-1$ & $2.0E+0$ & $-$ & $-$ \\
$\boldsymbol{m}_2$ & $-$      & $2.4E+0$ & $2.1E+0$ & $4.8E-1$ & $7.7E-1$ \\
$\boldsymbol{m}_3$ & $-$      & $-$      & $5.2E-1$ & $2.3E-1$       & $-$ \\ 
\bottomrule
\end{tabular}}
\label{table:Average_Error_dipole_source_another_idea}
\end{table}
\section{Conclusions}
%
We discuss a reconstruction method for
multiple moving point/dipole wave sources from boundary measurements.
We derive algebraic relations between the parameters of unknown sources and the reciprocity gap functionals for five sequences of
functions.
Based on these algebraic relations, we give a real-time reconstruction procedure for parameters of unknown
sources.
We examine our reconstruction procedure by numerical experiments.
Numerical results show that our procedure gives reliable reconstruction results for both moving
point and moving dipole sources when the observation noise
is smaller than 0.5\%.
However, the noise becomes 1\%, its influence becomes uninnorable, and
the reconstruction results become unreliable under 5\% noise.
Such bad influences can be found more heavily on the estimation results
of magnitudes/moments of sources.
\par
We need further discussions for more complicated cases, for example,
limited aperture cases, the case where the source term contains both
moving point and dipole sources, and so on. 
\par
\vspace{6ex}
%
%
%
\section*{Acknowledgements}
The author is grateful to Professor E. Nakaguchi in Tokyo Medical and 
Dental University for his useful discussions and helpful comments.
Also the author would like to acknowledge fruitful the discussions at 2017 IMI Joint Use Research
Program Workshop (II) "Practical inverse problems based on
interdisciplinary and industry-academia collaboration".
This work was supported in part by A3 Foresight Program ``Modeling and
Computation of Applied Inverse Problems'' and Grant-in-Aid for
Scientific Research (S) 15H05740 of Japan Society for the Promotion of
 Science.

\par
\vspace{6ex}
\section*{Appendix A. Derivation of
(\ref{eq:relation_RGf_parameter_point_source})-(\ref{eq:relation_RGj_parameter_point_source}) and 
(\ref{eq:relation_RGf_parameter_dipole_source})-(\ref{eq:relation_RGj_parameter_dipole_source})
}
Firstly, we show the derivation of (\ref{eq:relation_RGf_parameter_point_source})-(\ref{eq:relation_RGj_parameter_point_source}).
%
From the definition of $f_{n,\varepsilon}$, we have
\begin{align*}
 {\mathcal F}(f_n)(\tau) =&
 \displaystyle \lim_{\varepsilon \rightarrow +0}
  \sum_{k=1}^K \int_0^T q_k(t)  f_{n,\varepsilon}(t,\boldsymbol{p}_k(t);\tau)dt \\
=&
 \displaystyle \lim_{\varepsilon \rightarrow +0}
   \sum_{k=1}^K \int_0^T q_k(t) \cdot
  (p_{k,x}(t) + \mathrm{i} p_{k,y}(t))^n \cdot \eta_{\varepsilon}\left(t - \tau +
 							\frac{p_{k,z}(t)}{c}  \right) dt.
\end{align*}
Let $\displaystyle s = t +\frac{p_{k,z}(t)}{c}$ for each $k$, then
\begin{equation*}
  \frac{dt}{ds} = \left(1 - \frac{1}{c} \cdot \frac{d}{ds}(p_{k,z}(t_k(s)))\right) = \xi_k(s),
\end{equation*}
where $t_k(s)$ is defined in section 3.2, and we obtain
 \begin{equation}
 {\mathcal F}(f_n)(\tau) =
  \lim_{\varepsilon \rightarrow +0}
  \sum_{k=1}^K \int_{s_{0,k}}^{s_{1,k}} q_k(t_k(s))\xi_k(s) \cdot
  (p_{k,xy}(t_k(s)))^n \cdot\eta_{\varepsilon}(s  - \tau)
 ds,
  \label{eq:computation_RGf_point_source}
  \end{equation}
where $s_{0,k} = p_{k,z}(0)/c$ and $s_{1,k} = T + p_{k,z}(T)/c$.
Since $\eta_\varepsilon$ is the standard mollifier function, we
establish (\ref{eq:relation_RGf_parameter_point_source}).
\par
Next, from the definition of $g_{n,\varepsilon}$, we have
 \begin{align*}
 {\mathcal F}(g_n)(\tau) =&
 \displaystyle \lim_{\varepsilon \rightarrow +0}
  \sum_{k=1}^K \int_0^T q_k(t)  g_{n,\varepsilon}(t,\boldsymbol{p}_k(t);\tau)dt \\
=&
 \displaystyle -\lim_{\varepsilon \rightarrow +0}
  \sum_{k=1}^K \int_0^T q_k(t) \cdot\left. (p_{k,xy}(t))^n \cdot\eta^\prime_{\varepsilon}\left(t - \tau +
							\frac{z}{c}\right)\right|_{\boldsymbol{r}=\boldsymbol{p}_k(t)}dt \\
=&
 \displaystyle -\lim_{\varepsilon \rightarrow +0}
  \sum_{k=1}^K \int_0^T q_k(t) \cdot(p_{k,xy}(t))^n \cdot \eta^\prime_\varepsilon\left(t - \tau +
							\frac{p_{k,z}(t)}{c}\right)
  dt \\
=&
 \displaystyle -\lim_{\varepsilon \rightarrow +0}
  \sum_{k=1}^K \int_{s_{0,k}}^{s_{1,k}} q_k(t_k(s)) \xi_k(s) \cdot(p_{k,xy}(t_k(s)))^n
\cdot \eta^\prime_{\varepsilon}(s-\tau) 
  ds.
\end{align*}
Using integration by parts, it follows that
\begin{align}
 & {\mathcal F}(g_n)(\tau) \nonumber \\
= & -\displaystyle \lim_{\varepsilon \rightarrow +0}
  \sum_{k=1}^K  q_k(t_k(s_{1,k})) \xi_k(s_{1,k}) \cdot (p_{k,xy}(t_k(s_{1,k})))^n
 \cdot \eta_{\varepsilon}(s_{1,k}-\tau)  \nonumber \\
 & \displaystyle + \lim_{\varepsilon \rightarrow +0}
  \sum_{k=1}^K  q_k(t_k(s_{0,k})) \xi_k(s_{0,k}) \cdot (p_{k,xy}(t_k(s_{0,k})))^n \cdot \eta_{\varepsilon}(s_{0,k}-\tau) \nonumber \\
 &\displaystyle + \lim_{\varepsilon \rightarrow +0}
    \int_{s_{0,k}}^{s_{1,k}}
 \frac{d}{ds}\left(q_k(t_k(s))\xi_k(s) \cdot (p_{k,xy}(t_k(s)))^n
 \right) \cdot \eta_{\varepsilon}(s-\tau) ds  \nonumber \\
= &  \displaystyle \lim_{\varepsilon \rightarrow +0} \sum_{k=1}^K
    \int_{s_{0,k}}^{s_{1,k}} \frac{d}{ds}\left(q_k(t_k(s)) \xi_k(s)\right)
    \cdot (p_{k,xy}(t_k(s)))^n \cdot
    \eta_{\varepsilon}(s-\tau)   ds \nonumber \\
 &  + \displaystyle \lim_{\varepsilon \rightarrow +0} \sum_{k=1}^K
    \int_{s_{0,k}}^{s_{1,k}} n \cdot q_k(t_k(s))\xi_k(s) \cdot (p_{k,xy}(t_k(s)))^{n-1}
    \nonumber \\
 & \hspace{15ex} \times\frac{d}{ds} \left(p_{k,xy}(t_k(s))\right)
   \cdot \eta_{\varepsilon}(s-\tau)  ds.
   \label{eq:computation_RGg_point_source}
\end{align}
Here, we use $t'_k(s) = \xi_k(s),\  t_k(s_{0,k})= 0,\ t_k(s_{1,k}) = T$,
and 
$\eta_\varepsilon(s_{0,k}-\tau) = \eta_\varepsilon(s_{1,k}-\tau) = 0$.
Then, using similar computation as (\ref{eq:computation_RGf_point_source}), we
establish (\ref{eq:relation_RGg_parameter_point_source}).
\par
Also for ${\mathcal F}(h_n)(\tau),\ {\mathcal F}(i_n)(\tau)$ and
${\mathcal F}(j_n)(\tau)$, we can obtain the following form by 
changing variable $t$ to $s$ and integration by parts:
\begin{align*}
  & {\mathcal F}(h_n)(\tau)  \nonumber \\
= &  \displaystyle \lim_{\varepsilon \rightarrow +0} \sum_{k=1}^K
    2n \int_{s_{0,k}}^{s_{1,k}}
    q_k(t_k(s)) \xi_k(s) \cdot p_{k,z}(t_k(s)) \cdot (p_{k,xy}(t_k(s)))^{n-1}
    \cdot \eta_\varepsilon(s-\tau) ds \nonumber \\
 & + \displaystyle \lim_{\varepsilon \rightarrow +0} \sum_{k=1}^K
    \frac{1}{c}
    \int_{s_{0,k}}^{s_{1,k}}
    \frac{d}{ds}\left(q_k(t_k(s)) \xi_k(s) \cdot
 \overline{p_{k,xy}(t_k(s))} \cdot (p_{k,xy}(t_k(s)))^n \right) \cdot
    \eta_\varepsilon(s-\tau) ds \\[4ex]
 & {\mathcal F}(i_n)(\tau) \nonumber \\
= & \displaystyle \lim_{\varepsilon \rightarrow +0}  \sum_{k=1}^K 
    \int_{s_{0,k}}^{s_{1,k}} 
   \frac{d^2}{ds^2}\left(q_k(t_k(s))\xi_k(s) \cdot (p_{k,xy}(t_k(s)))^n 
		    \right) \cdot
   \eta_{\varepsilon}(s-\tau)  ds, 
\end{align*}
\begin{align*}
 & {\mathcal F}(j_n)(\tau) \\
= &  \displaystyle \lim_{\varepsilon \rightarrow +0} \sum_{k=1}^K
    2n \int_{s_{0,k}}^{s_{1,k}}
    \frac{d}{ds}\left(q_k(t_k(s)) \xi_k(s) \cdot p_{k,z}(t_k(s)) \cdot (p_{k,xy}(t_k(s)))^{n-1}
    \right) \eta_\varepsilon(s-\tau) ds \nonumber \\
 &  + \displaystyle \lim_{\varepsilon \rightarrow +0} \sum_{k=1}^K
    \frac{1}{c} \int_{s_{0,k}}^{s_{1,k}}
    \frac{d^2}{ds^2}\left(q_k(t_k(s)) \xi_k(s) \cdot
 \overline{p_{k,xy}(t_k(s))} \cdot(p_{k,xy}(t_k(s)))^n \right)
    \eta_\varepsilon(s-\tau) ds.
\end{align*}
Then, using similar computation as (\ref{eq:computation_RGg_point_source}), we
establish (\ref{eq:relation_RGh_parameter_point_source})-(\ref{eq:relation_RGj_parameter_point_source}).
\par
Next, we derive (\ref{eq:relation_RGf_parameter_dipole_source})-(\ref{eq:relation_RGj_parameter_dipole_source}) for moving dipole sources.
%
From the definition of $f_{n,\varepsilon}$ and the assumption that $m_{k,z}(t)\equiv 0$, we have
\begin{align*}
    {\mathcal F}(f_n)(\tau)
= & \displaystyle \lim_{\varepsilon \rightarrow +0}
  \sum_{k=1}^K \int_0^T \boldsymbol{m}_{k}(t)  \cdot \left. \nabla f_{n,\varepsilon}(t,\boldsymbol{r};\tau) \right|_{\boldsymbol{r}=\boldsymbol{p_{k}(t)}}dt \\
= &
 \displaystyle \lim_{\varepsilon \rightarrow +0}
   \sum_{k=1}^K \int_0^T n (m_{k,x}(t)+ \mathrm{i}m_{k,y}(t)) \cdot
  (p_{k,x}(t) + \mathrm{i} p_{k,y}(t))^{n-1} \\
 & \hspace{30ex}  \times \eta_{\varepsilon}\left(t - \tau + \frac{p_{k,z}(t)}{c}\right) dt. \\
= &
 \displaystyle \lim_{\varepsilon \rightarrow +0}
   \sum_{k=1}^K \int_{s_{0,k}}^{s_{1,k}} n \cdot m_{k,xy}(t_k(s))\xi_k(s) \cdot
   p_{k,xy}(t_k(s))^{n-1}\cdot
  \eta_{\varepsilon}(s - \tau) ds. 
\end{align*}
Since $\eta_\varepsilon$ is the standard mollifier function, we
establish (\ref{eq:relation_RGf_parameter_dipole_source}).
\par
For ${\mathcal F}(g_n)(\tau),\ {\mathcal F}(h_n)(\tau),\ {\mathcal
F}(i_n)(\tau),\ {\mathcal F}(j_n)(\tau)$, changing variable from $t$
to $s$ and integration by parts, we obtain the following form:
\begin{align*}
 & {\mathcal F}(g_n)(\tau) \nonumber \\
= & \displaystyle \lim_{\varepsilon \rightarrow +0}  \sum_{k=1}^K 
    n \int_{s_{0,k}}^{s_{1,k}} 
    \frac{d}{ds}\left(m_{k,xy}(t_k(s))\xi_k(s) \cdot (p_{k,xy}(t_k(s)))^{n-1} 
		    \right)
   \eta_{\varepsilon}(s-\tau)  ds,  \\[4ex]
%
 & {\mathcal F}(h_n)(\tau) \nonumber \\
= &  \displaystyle \lim_{\varepsilon \rightarrow +0} \sum_{k=1}^K
    2n(n-1) \int_{s_{0,k}}^{s_{1,k}}
    m_{k,xy}(t_k(s)) \xi_k(s) \cdot p_{k,z}(t_k(s)) \cdot
 (p_{k,xy}(t_k(s)))^{n-2} \cdot
     \eta_\varepsilon(s-\tau) ds \nonumber \\
 & + \displaystyle \lim_{\varepsilon \rightarrow +0} \sum_{k=1}^K
    \frac{1}{c}
    \int_{s_{0,k}}^{s_{1,k}}
    \frac{d}{ds}\left(\overline{m_{k,xy}(t_k(s))} \xi_k(s) \cdot
		 (p_{k,xy}(t_k(s)))^n \right)
    \eta_\varepsilon(s-\tau) ds \nonumber \\
 & + \displaystyle \lim_{\varepsilon \rightarrow +0} \sum_{k=1}^K
    \frac{n}{c}
    \int_{s_{0,k}}^{s_{1,k}}
    \frac{d}{ds}\left(m_{k,xy}(t_k(s)) \xi_k(s) \cdot
 \overline{p_{k,xy}(t_k(s))} \cdot
		 (p_{k,xy}(t_k(s)))^{n-1} \right)
    \eta_\varepsilon(s-\tau) ds  \\[4ex]
 & {\mathcal F}(i_n)(\tau) \nonumber \\
= & \displaystyle \lim_{\varepsilon \rightarrow +0}  \sum_{k=1}^K 
   n \int_{s_{0,k}}^{s_{1,k}} 
    \frac{d^2}{ds^2}\left(m_{k,xy}(t_k(s))\xi_k(s) \cdot (p_{k,xy}(t_k(s)))^{n-1} 
		    \right)
   \eta_{\varepsilon}(s-\tau)  ds, 
%
\end{align*}
\begin{align*}
 & {\mathcal F}(j_n)(\tau) \\
= &  \displaystyle \lim_{\varepsilon \rightarrow +0} \sum_{k=1}^K
    2n(n-1) \int_{s_{0,k}}^{s_{1,k}}
    \frac{d}{ds}\left(m_{k,xy}(t_k(s)) \xi_k(s) \cdot
 p_{k,z}(t_k(s)) \cdot (p_{k,xy}(t_k(s)))^{n-2}
    \right) \eta_\varepsilon(s-\tau) ds \nonumber \\
 &  + \displaystyle \lim_{\varepsilon \rightarrow +0} \sum_{k=1}^K
   \frac{1}{c} \int_{s_{0,k}}^{s_{1,k}}
    \frac{d^2}{ds^2}\left(\overline{m_{k,xy}(t_k(s))}\xi_k(s) \cdot
		     (p_{k,xy}(t_k(s)))^{n} \right)
    \eta_\varepsilon(s-\tau) ds \nonumber \\
 &  + \displaystyle \lim_{\varepsilon \rightarrow +0} \sum_{k=1}^K
    \frac{n}{c} \int_{s_{0,k}}^{s_{1,k}}
    \frac{d^2}{ds^2}\left(m_{k,xy}(t_k(s))\xi_k(s) \cdot
 \overline{p_{k,xy}(t_k(s))} \cdot 
		     (p_{k,xy}(t_k(s)))^{n-1} \right)
    \eta_\varepsilon(s-\tau) ds.
\end{align*}
Then, using similar computation as (\ref{eq:computation_RGg_point_source}), we
establish (\ref{eq:relation_RGg_parameter_dipole_source})-(\ref{eq:relation_RGj_parameter_dipole_source}).
\par
%
%
%
\section*{Appendix B. A proof of equation (\ref{eq:detV})}
In this appendix, we set $K(\tau) = K$ for simplicity, and omit the
argument $(t_j(\tau))$ in $\boldsymbol{p}_j(t_j(\tau))$.
\par
In $\det \tilde{V}(\tau)$, subtract 1st column from 2nd column, and $(K+1)$-th
column from $(K+2)$-th column, we obtain
\begin{equation}
\det \tilde{V}(\tau) = 
  \det \left(
        \begin{array}{cccccccccc}
        {\bm \psi}_1  & {\bm \psi}_2  - {\bm \psi}_1  & {\bm \psi}_3 & \cdots & {\bm \psi}_K &
        {\bm \psi}_1' & {\bm \psi}_2' - {\bm \psi}_1' & {\bm \psi}_3'& \cdots & {\bm \psi}_K'
\end{array}
\right).
 \label{eq:AppendixB_1}
\end{equation}
Since the 2nd and $(K+2)$-th columns of (\ref{eq:AppendixB_1}) have the factor $(p_{2,xy} -
p_{1,xy})$, we can factorise $\det \tilde{V}(\tau)$ as
\begin{equation*}
\det \tilde{V}(\tau) = (p_{2,xy} - p_{1,xy})^2 \cdot
  \det \left(
        \begin{array}{cccccccccc}
        {\bm \psi}_1  & {\bm \sigma}_2 & {\bm \psi}_3  & \cdots & {\bm \psi}_K &
        {\bm \psi}_1' & {\bm \sigma}_2'& {\bm \psi}_3' & \cdots & {\bm \psi}_K'
\end{array}
\right).
\end{equation*}
where the $k$-th components of 2$K$-vectors ${\bm \sigma}_2$ and ${\bm \sigma}_2'$ are expressed by
\begin{align*}
(\bm \sigma_2)_k &= \left\{
          \begin{array}{cc}
             0, & k = 1, \\
             1, & k = 2, \\[2ex]
             \displaystyle
             \sum_{l=0}^{k-2} (p_{2,xy})^{k-l-2} \cdot (p_{1,xy})^{l}, &
             k=3,4,\cdots 2K,
          \end{array}
\right. \nonumber \\[2ex]
(\bm \sigma_2')_k &= \left\{
          \begin{array}{cc}
             0, & k = 1,2,  \\
             2, & k = 3,    \\[2ex]
             \displaystyle
             (k-1)\sum_{l=0}^{k-3} (p_{2,xy})^{k-l-3} \cdot (p_{1,xy})^{l}, &
             k=4,5, \cdots 2K.
          \end{array}
\right.
 \end{align*}
\par
Also we can derive factor $(p_{2,xy} - p_{1,xy})$ two times
by the following two steps
\begin{itemize}
\item Subtract $(K +1)$-th column from the 2nd column, then 2nd column
      has the factor $(p_{2,xy} - p_{1,xy})$.
      Hence, we can factorise by $(p_{2,xy} - p_{1,xy})$.
\item After above step, subtract twice of 2nd column from $(K+2)$-th
      column, then $(K+2)$-th column has the factor $(p_{2,xy} -
      p_{1,xy})$.
\end{itemize}
Therefore, $\det \tilde{V}(\tau)$ has the factor $(p_{2,xy}-p_{1,xy})^4$.
\par
By similar computation, we can derive the factor $(p_{j,xy} - p_{k,xy})^4$
for any $j > k$.
It is easily see that the largest order of $p_{k,xy}$ in $\det
\tilde{V}(\tau)$ is $4(K-1)$ for any $k$,  then we can factorise $\det \tilde{V}(\tau)$ as
\begin{equation*}
  \det \tilde{V}(\tau) = C_K \cdot\prod_{j > k} (p_{j,xy} - p_{k,xy})^4,
\end{equation*}
where $C_K$ is a constant which depends only on $K$.
\par
Let us determine the constant $C_K$.
By replacing columns $(K-1)$ times and rows
$2(2K-2)$ times in $\det \tilde{V}(\tau)$, we obtain
\begin{align*}
& \det \tilde{V}(\tau) \\
= &  (-1)^{(K-1)} \\
   & \times \det \left(
        \begin{array}{cc|cccccccc}
        (p_{1,xy})^{2K-2} & (2K-2)(p_{1,xy})^{2K-3} & \multicolumn{8}{c}{\cdots} \\
        (p_{1,xy})^{2K-1}   & (2K-1)(p_{1,xy})^{2K-2}     &
	 \multicolumn{8}{c}{\cdots} \\[1ex]
         \hline 
	 & & \multicolumn{8}{c}{} \\[-2ex]
        \hat{\bm \psi}_1           & \hat{\bm \psi}'_1                  & \hat{\bm \psi}_2  &
        \hat{\bm \psi}_3 & \cdots  &\hat{\bm \psi}_K & \hat{\bm \psi}'_2  & \hat{\bm \psi}_3'& \cdots & \hat{\bm \psi}'_K
  \end{array}
\right), 
\end{align*}
where
\begin{equation*}
\hat{\bm \psi}_k = \left(
      \begin{array}{c}
      1 \\ 
      p_{k,xy} \\ 
     (p_{k,xy})^2 \\
     (p_{k,xy})^3 \\ 
      \vdots \\ 
     (p_{k,xy})^{2K-3}
      \end{array}
\right), \qquad
\hat{\bm \psi}'_k = \left(
      \begin{array}{c}
      0 \\ 1 \\ 
      2p_{k,xy} \\ 
      3(p_{k,xy})^2 \\ 
      \vdots \\ 
      (2K-3)(p_{k,xy})^{2K-4}
      \end{array}
\right).
\end{equation*}
Therefore the coefficient of the term $(p_{1,xy})^{4K-4}$ is given by
\begin{align}
 &   \displaystyle (-1)^{(K-1)} \cdot \left((2K-1)- (2K-2) \right) \cdot \det
   \left(
        \begin{array}{cccccccc}
        \hat{\bm \psi}_2   & \hat{\bm \psi}_3  & \cdots & \hat{\bm \psi}_K  &
	\hat{\bm \psi}'_2 & \hat{\bm \psi}'_3  & \cdots &  \hat{\bm \psi}'_K
  \end{array}
\right) \nonumber \\
 = & (-1)^{(K-1)} \cdot C_{K-1}\cdot \prod_{j>k\geq 2}(p_{j,xy} - p_{k,xy})^4.
\label{eq:deriviation_CK}
\end{align}
Applying equation (\ref{eq:deriviation_CK}) inductively, we obtain
\begin{equation*}
  C_K = (-1)^{(K-1)} C_{K-1} = \cdots = (-1)^{(K-1) + (K-2) + \cdots + 2} \cdot C_2.
\end{equation*}
Here, 
\begin{equation*}
  \det \tilde{V}_2 = 
\det \left(
     \begin{array}{cccc} 1 & 1 & 0 & 0 \\
                         p_{1,xy} & p_{2,xy} & 1 & 1 \\
                         (p_{1,xy})^2 & (p_{2,xy})^2 & 2 p_{1,xy} & 2p_{2,xy}  \\
                         (p_{1,xy})^3 & (p_{2,xy})^3 & 3(p_{1,xy})^2 & 3(p_{2,xy})^2
     \end{array}
     \right)
     = (-1)\cdot(p_{2,xy} - p_{1,xy})^4,
\end{equation*}
we have $C_2 = -1$.
Then $C_K= (-1)^{K(K-1)/2}$, and we complete the proof of (\ref{eq:detV}).
\par
%


\end{document}